\documentclass[a4paper,12pt]{article}
\usepackage{amsmath}
\usepackage{enumitem}
\usepackage{bm}
\usepackage{cite}
\usepackage{makecell}
\usepackage{mathtools}
\usepackage{amssymb, amsbsy, amstext}
\usepackage{srcltx}
\usepackage{amscd,amsfonts,color}
\usepackage[figuresright]{rotating}
\usepackage{url}

\input amssym.def
\input amssym.tex

\newtheorem{lem}{Lemma}[section]%
\newtheorem{theorem}[lem]{Theorem}%
\newtheorem{prop}[lem]{Proposition}%
\newtheorem{rem}[lem]{Remark}%
\newtheorem{conj}[lem]{Conjecture}%

\def\a{\alpha} \def\b{\beta} \def\g{\gamma} \def\d{\delta} \def\e{\varepsilon}
  \def\t{\tau}  \def\k{\kappa}
\def\th{\theta} \def\ld{\lambda} 
\def\Del{\Delta}  
\def\si{\Sigma} \def\O{\Omega} \def\G{\Gamma}
\def\De{\Delta}
\def\D{{\rm D}}

\def\di{\bigm|} \def\lg{\langle} \def\rg{\rangle}

\def\Alt{\hbox{\rm Alt}}
\def\PSL{\hbox{\rm PSL}}  
\def\Aut{\hbox{\rm Aut\,}}  \def\Syl{\hbox{\rm Syl}}
 \def\soc{\hbox{\rm soc}} \def\Fix{\hbox{\rm Fix }}
  \def\mod{\hbox{\rm mod }}
  
 \def\PG{\hbox{\rm PG}} \def\PGL{\hbox{\rm PGL}}
  \def\GL{\hbox{\rm GL}}  \def\P\GL{\hbox{\rm P\GL}}
 \def\SL{\hbox{\rm SL}} \def\FF{{\hbox{\sf F\kern-.43emF}}}
\def\cal{\mathcal}  
\def\i{\imath}
\def\j{\jmath}

\def\PSigmaU{{\rm P\Sigma U}}
\def\PGammaU{{\rm P\Gamma U}}
\def\Sym{\hbox{\rm Sym}}

\def\N{\hbox{\rm N}}
\def\C{\hbox{\rm C}}
\def\Z{\hbox{\rm Z}}
\def\A{\hbox{\rm A}}
\def\soc{\hbox{\rm soc}}
\def\SS{\hbox{\rm S}}

\def\Sz{\hbox{\rm Sz}}
\def\Ree{\hbox{\rm Ree}}
\def\GU{\hbox{\rm GU}}
\def\PGU{\hbox{\rm PGU}}
\def\PSU{\hbox{\rm PSU}}
\def\PGammaU{{\rm P\Gamma U}}
\def\SO{\hbox{\rm SO}}
\def\SU{\hbox{\rm SU}}
\def\PSO{\hbox{\rm PSO}}
\def\PGO{\hbox{\rm PGO}}
\def\PSigmaU{{\rm P\Sigma U}}

\def\det{{\rm det}}
\def\mat {{\rm mat}}
\def\iso {{\rm iso}}
\def\rank {{\rm rank}}

\def\ZZ{\mathbb{Z}} 
    \def\HH{\mathcal{H}}

\def\nd{\mathrel{\bigm|\kern-.7em/}} 
 \def\f{\noindent}
\def\qed{\hfill $\Box$} \def\demo{\f {\bf Proof}\hskip10pt}

\usepackage{titlesec}

\titleformat{\paragraph}
{\normalfont\normalsize\bfseries}{\theparagraph}{1em}{}
\titlespacing*{\paragraph}
{0pt}{3.25ex plus 1ex minus .2ex}{1.5ex plus .2ex}

\begin{document}

\begin{center} {\bf\large The Burness-Giudici Conjecture on
  \\Primitive Groups with Socle $\PSU(3,q)$}
\vskip 3mm
{\sc Huye Chen}\\
{\footnotesize
School of Mathematics, Guangxi University,   Nanning 530004, P. R. China}\\

{\sc Shaofei Du}\\
{\footnotesize
School of Mathematical Sciences, Capital Normal University, Beijing 100048, P. R. China}\\

{\sc Weicong Li}\\
{\footnotesize
School of Sciences, Great Bay University, Dongguan 523000, P. R. China}\\

\end{center}

\renewcommand{\thefootnote}{\empty}
  \footnotetext{E-mail addresses:  chenhy280@gxu.edu.cn (H. Chen), dushf@mail.cnu.edu.cn (S. Du),
  liweicong@gbu.edu.cn (W. Li).}

\begin{abstract}
Let $G$ be a transitive permutation group on a set $\O$, and suppose
$G_\a\cap G_\b=1$ for some distinct $\a, \b\in\O$. The Saxl graph $\Sigma(G)$ of $(G, \O)$ is defined as the graph with vertex set $\O$, where two vertices $\a', \b'$ are adjacent if and only if $G_{\a'}\cap G_{\b'}=1$.
Burness and Giudici conjectured that for every primitive permutation group $G$, its Saxl graph has the property that any two vertices share a common neighbor.

We focus on proving the conjecture for all primitive groups $G$ whose socle is a simple group of  Lie-type of rank $1$; that is, $\soc(G)\in \{\PSL(2,q),\PSU(3,q), \Ree(q),$ \\ $\Sz(q)\}$.
The case $\soc(G)=\PSL(2,q)$ has been treated in two earlier papers. The purpose of the present paper is to settle the case $\soc(G)=\PSU(3,q)$.
To finsh this work, we draw on methods from abstract- and permutation- group theory, finite unitary geometry,  probabilistic approach, number theory (employing Weil's bound), and, most importantly, algebraic combinatorics, which provides us some key ideas.
 \end{abstract}
\section{Introduction}
 A {\em base} for a finite permutation group $G$ on a set $\Omega$ is a subset $\Delta\subset \Omega$ whose pointwise  stabilizer $G_{(\Del)}$ is trivial.  This  concept  is a natural generalization of basis for  vector spaces. The {\em base size} $b(G)$ of $G$  is the minimal cardinality $|\Delta|$, where $\Del$ goes  over all bases for $G$, while  a  base $\Delta '$ with $|\Delta '|=b(G)$  is called a {\em base size set}. With a history dating to the 19th century, these concepts remain vital to contemporary algebra and combinatorics; see \cite{Bailey-Cameron,LS} and \cite[$\S5$]{B} for overviews and \cite{G,S} for applications.

Primitive permutation groups with a base of size $1$ are cyclic groups $\ZZ_p$, which makes the family of finite primitive groups with a base of size 2 particularly significant due to their notably small base.
Saxl's ``base-two programme", which aims to classify finite primitive permutation groups with base size $2$, has inspired a great deal of research. To study the bases of permutation groups $G$ with $b(G)=2$, Burness and Giudici \cite{BG} introduced the {\em Saxl graph}, which is defined in terms of these groups. Let $G$ be a permutation group on $\Omega$ with $b(G)=2$. The vertex set of a Saxl graph $\Sigma(G)$ (denoted by $\Sigma$, simply) of $G$ is just
$\Omega$, while two vertices are adjacent if and only if they form a base.
  In \cite{BG}, the authors discussed  the valency, connectivity, hamiltonicity and the independence number of $\Sigma(G)$, and   proposed the following conjecture.
\vskip -10pt
\begin{conj}$($BG-Conjecture$)$ \label{BG1}
If $G$ is a finite primitive permutation group with $b(G)=2$, then every pair of vertices in its Saxl graph $\Sigma(G)$ has a common neighbor.
\end{conj}
\vskip -5pt

 From now on, we assume that $G\leq \Sym(\O)$ acts transitively on $\O$. Fixing a point $\alpha\in \Omega$, we call the orbits of $G_\a$ the {\em suborbits} of $G$ relative to $\a$; the trivial suborbit is $\{\a\}$ itself.
   Now $\{\alpha,\beta\}$ is a base for $G$ with $b(G)=2$ if and only if $G_{\alpha}$ acts regularly on the suborbit containing $\b$.
   It follows that the neighborhood $\si_1(\a)$ of $\alpha$ in
  $\Sigma(G)$ is the union $\G$ of all regular suborbits of $G$ relative to $\a$.
Using a probabilistic approach, Burness and Giudici \cite{BG} verified the conjecture for a range of almost simple groups. For instance, it holds when $G=S_n$ or $A_n$ (with $n>12$) and the point stabilizer is primitive. Their results also cover certain primitive groups of diagonal type and twisted wreath products of sufficiently large order, which have been extended recently by Huang (see Sections 5.8 and 5.9 of \cite{HPHD}).
Furthermore, by employing the Magma database, they verified the conjecture for all primitive groups of degree at most $4095$. The paper \cite{BG} also confirms the conjecture for many sporadic simple groups.
Burness and Huang \cite{BH} later proved the conjecture for almost simple primitive groups with soluble point stabilisers; additionally, they studied it for primitive groups of product type in a separate paper \cite{BH2}.
Lee and Popiel \cite{LP} have recently proved the conjecture for most affine-type groups with sporadic point stabilizers. The theory was subsequently extended with the introduction of generalized Saxl graphs
and Saxl hypergraphs in
\cite{FH, LP1}.

Resolving this conjecture is widely recognized as a hard problem, with very few cases having been established for groups of Lie type. It is therefore sensible to focus first on rank one groups. At present, there appear to be no viable ideas for extending the approach to higher ranks, especially considering that a full classification of base-two primitive almost simple groups remains out of reach.
In this series of papers \cite{Chen-Du, Chen-Du-Li1, Chen-Du2}, we have established an effective and general method, evolving from the initial study in \cite{Chen-Du}. Our procedure consists of the following steps.
First, we determine the valency of the Saxl graph. If the valency is large enough (analogous to the $\Ree(q)$ and $\Sz(q)$ cases in \cite{Chen-Du2}, or to specific instances in \cite{Chen-Du-Li1}), we then employ a computational approach combining probabilistic arguments with a precise determination of suborbits.
  When the valency is smaller than half the number of vertices, we instead construct an associated bipartite graph between the stabilizer subgroups and the set of involutions in $G$ (see Sections 3-6 of the present paper or \cite[Section 3.1]{Chen-Du}).
  Furthermore, certain special cases can be resolved by directly analyzing their underlying geometric actions (cf. \cite[Section 3.3]{Chen-Du}).
  It is important to note that cases where the valency is less than half the number of vertices constitute the most challenging part, typically requiring a combination of geometric actions, suborbit calculations, and probabilistic methods.

 Our research aims to tackle this conjecture in case of groups of Lie-type of rank one.
  This class of groups contains four families, with  respective socles $\PSL(2,q)$,   $\PSU(3,q)$, $\Ree(q)$ and $\Sz(q)$.
Having established the conjecture for groups with socle $\PSL(2,q)$ in \cite{BH,Chen-Du},
  we now address the case of $\PSU(3,q)$ in this paper, that is
\begin{theorem} \label{main} Let $G$ be a primitive permutation group with socle $\PSU(3,q)$, where $q=p^m\ge 7$ for an odd prime $p$. Suppose a point-stabilizer
  $M$ of $G$ contains $\PSO(3,q)$. Then $b(G)=2$ and either (i) $G=\PSU(3,q)$ and $M=\SO(3,q)$; or (ii) $G=\PSU(3,q)\ZZ_{m_1}$ and  $M=\SO(3,q)\ZZ_{m_1}$ $(m_1\mid 2m$ and  $m_1\ne  1)$. Further, the BG-Conjecture is valid for  case (i);
  and for   $\sqrt{q}\geq 17$ or $45$,  if $3\nmid (q+1)$  and $3\di (q+1)$, respectively in case (ii).
  \end{theorem}

Remark that for the reason of the length of the paper, we are not able to  improve the lower bound  $\sqrt{q}\geq 17$ (resp. $45$) for small exceptions in case (ii) of Theorem~\ref{main}.
 Hopefully,  these  small cases may be verified directly with a sufficiently powerful computer.

Combining Theorem~\ref{main} and the main results of \cite{BH,Chen-Du,Chen-Du2}, we get
\begin{theorem} \label{main'} Let $G$ be  a primitive  group whose socle is a simple group of  Lie-type of rank $1$ such that $b(G)=2$.
 Then  the BG-Conjecture is valid, except for $G=\PSU(3,q)\ZZ_{m_1}$ and  $M=\SO(3,q)\ZZ_{m_1}$ $(m_1\mid 2m$ and  $m_1\ne  1)$, where
    $\sqrt{q}\le 17$ or $45$,  if $3\nmid (q+1)$  and $3\di (q+1)$, respectively.
  \end{theorem}

\vskip 3mm
\f {\bf  Outline  of the proof of Theorem~\ref{main}}: From now on, let $G$ be an almost simple group with the socle $T=\PSU(3,q)$, so that $T\leq G\leq \Aut(T)$, where $q=p^m$.
 Let $\a=M$ be a maximal subgroup of $G$ containing $M_0=\PSO(3,q)$, up to conjugacy. Consider the primitive permutation representations of $G$ on the right cosets $[G:M]$.
 Then $q\ge7$ and either: (i) $G=T$ and $M=M_0$; or (ii) $G=T\ZZ_{m_1}$ and  $M=M_0\ZZ_{m_1}$, where $m_1\mid 2m$ and $m_1\ne 1$.
For the first case, we exhibit a common  neighbor of $\a $ and $\a^g$ in the Saxl graph for any $g\in G$ via a direct geometric construction in Section 2.
However, we are not able to extend  this geometric construction   to the second case, which necessitates a new approach developed in Section 3.
Consequently, the proof of Theorem~\ref{main} is obtained by combining Theorems~\ref{main1} and \ref{main2}.
Moreover, to prove Theorem~\ref{main2} crucially depends on estimates of $|\Gamma_r(\alpha)|$ and $|\Gamma_{nr}(\alpha) \cap \Gamma_{nr}(\alpha')|$, developed in Sections 4 and 5 respectively. Recall that $\Gamma_r(\alpha)$ (resp. $\Gamma_{nr}(\alpha)$) is the union of all regular (resp. non-regular) suborbits of $G$ relative to $\alpha$.
\vskip -7pt
\section{$G=\PSU(3,q)$ and $M=\PSO(3,q)$}\label{geo}
 Theorem~\ref{main1} is the main result of this section, which will be proved in   Section~\ref{geoproof}.  To do that,
   a geometrical description of the action of $G$ is presented in
   Section~\ref{geometrical description}.
    \vskip -10pt\begin{theorem} \label{main1}  Let $G=\PSU(3,q)$ and $M=\PSO(3,q),$ a maximal subgroup of $G$, where $q=p^m\ge 7$ for an odd prime $p$.
Consider the primitive right multiplication action  of $G$ on the set of right cosets $\Omega =[G:M]$.  Then $b(G)=2$ and the BG-Conjecture holds.
  \end{theorem}
\subsection{A geometrical description of  the action}\label{geometrical description}
As usual, let $V:=V(3,q^2)$ and   $\PG(2,q^2)$ denote the $3$-dimensional row vector space  and  the  $2$-dimensional projective space on the finite field  $\FF_{q^2}.$
 Every point $P$ in   $\PG(2,q^2)$ is written as $\lg(x_0,x_1,x_2)\rg$ and for  any   subset $S\subseteq V$, set   $\PG(S):=\{\lg s\rg\mid s\in S\}$ and in particular,
 $\PG(S)$ is a line if $S$ is a 2-dimensional subspace of $V$.
  For  a subfield  $\FF_{q'}$  of $\FF_{q^2}$ and $v_{\imath}\in V$ with $1\leq \imath\leq k$, set $\lg v_1,v_2,\cdots,v_k\rg_{\FF_{q'}}=\{ \sum_{\i=1}^k a_\i v_\i\mid a_{\imath}\in\FF_{q'}, 1\leq \imath\leq k\}.$ Throughout this paper,  $\lg v_1,v_2,\cdots,v_k\rg_{\FF_{q^2}}$ is simply denoted by $\lg v_1,v_2,\cdots,v_k\rg$.

 From now on, let $V$ be a vector space equipped with a Hermitian form $( \cdot,\cdot)$. A {\it Baer subline} (resp. {\it subplane}) of $\PG(1,q^2)$ (resp. $\PG(2,q^2)$) is the orbit of $\PG(1,q)$ (resp. $\PG(2,q)$)  under the group $\PGL(2,q^2)$ and (resp. $\PGL(3,q^2)$).
Following the definition in \cite[p.~311]{suzuki}, three points $\lg\alpha\rg$, $\lg\beta\rg$, and $\lg\gamma\rg$ in $\PG(2,q^2)$ or $\PG(2,q)$ are said to form a {\it frame} if they are non-collinear.
Consequently, the corresponding vectors $\alpha$, $\beta$, and $\gamma$ are linearly independent in $V(3, q^2)$ and thus form a basis.
Let $\hat{g}$ denote the image in $\PGU(3,q)$ of an element $g\in\GU(3,q)$.
Throughout this paper, let $W_0=\lg e_1, e_2, e_3\rg_{\FF_q}$, where $e_1, e_2$ and $e_3$  are unit vectors so that
$\mathcal{B}_0:=\{\PG(W_0)^{\hat{g}}:=\PG(W_0^g)\mid \hat{g}\in G=\PSU(3,q)\}$, where $G_{\PG(W_0)}=\PSO(3,q)\cong \SO(3,q)\cong \PGL(2,q)$. Then the right multiplication of $\PSU(3,q)$ on the right cosets  $[\PSU(3,q):\PSO(3,q)]$ is equivalent to its  natural action on its orbit $\mathcal{B}_0.$
Furthermore, let $\mathcal{B}:=\{\PG(W_0)^{\hat{g}}:=\PG(W_0^g)\mid g\in \GU(3,q)\}$.
 For a notational convenience, we denote each element $\mathrm{PG}(W) \in \mathcal{B}$ by the corresponding lowercase bold letter $\mathbf{w}$; thus, $\mathbf{w} = \mathrm{PG}(W).$
\begin{prop}\label{lem_obser}
 Let $\mathbf{u}$ be a Baer subplane of $\PG(2,q^2)$ and let $\lg v\rg$ be a point. Then
     \vspace{-5pt} \begin{itemize}[itemsep=0pt]
 \item[{\rm(1)}]{\rm(\cite[Theorem 2.6]{BE})}
      {\it There is a unique Baer subline of $\PG(1,q^2)$ passing  any given  three distinct points of $\PG(1,q^2);$}
 \item[{\rm(2)}] \rm{(\cite[p.26]{BE})}
     {\it  Any line $\ell$ of $\PG(2,q^2)$ passing  two points of $\mathbf{u}$ must necessarily intersect $\mathbf{u}$ in the unique line of $\mathbf{u}$,  containing  precisely $q+1$ points of $\mathbf{u};$}
   \item[{\rm(3)}] \rm{(\cite[Theorem 3.2]{BE})} {\it For any  $\lg v\rg \not\in \mathbf{u},$ among $q^2+1$ lines through $\lg v\rg $  $($in $\PG(2,q^2)$$)$, there are $q^2$ of them intersecting $\mathbf{u}$ in one point and one of them  intersecting $\mathbf{u}$ in  $q+1$ points,  forming a Baer subline of $\mathbf{u}$. In particular, no line in $\PG(2,q^2)$ is disjoint from $\mathbf{u};$ and }
   \item[{\rm(4)}] {\it If $\mathbf{u}\in\mathcal{B}$, then $ \PG(v^\perp)\cap \mathbf{u}$  is either a single point, if $\lg v\rg \not\in \mathbf{u}$, or a Baer subline contained in $\mathbf{u}$,  if $\lg v\rg \in \mathbf{u}$.  (Drived from (1) and (2)).}       \end{itemize}
\end{prop}

\vskip -6pt
 A vector $y\in V$ is called {\it non-isotropic} and {\it isotropic}  if $(y,y)\neq 0$ and  $(y,y)=0$, respectively. Let ${\cal {N}}$ and ${\cal H}$ denote the sets of non-isotropic and isotropic points, respectively, of $\PG(2,q^2)$.  A frame $\{\lg\a\rg,\lg\b\rg,\lg\g\rg\}$ is called {\it orthogonal} if the corresponding vectors $\alpha, \beta, \gamma$ are pairwise orthogonal.
 Throughout this paper we let $\mathcal{F}_{\mathbf{w}}$ denote the set of all the orthogonal frames of the Baer subplane $\mathbf{w}.$
For any   $\lg y\rg \in {\cal {N}},$    a well-known  involution $\tau_y$ in  $\PSU(3,q)$ is given by $\t_y(\lg x\rg) =\lg \t_y(x)\rg $, where
\vspace{-5pt}\begin{equation}\label{eqn_inv}
	\tau_y(x)= -x+2\frac{(x,y)}{(y,y)}y,\ \forall\, y\in V.
\end{equation}
It is a routine checking that $\tau_{\lambda y}(x)=\tau_{y}(x)$ for any $\lambda \in\FF_{q^2}^*$ and  $\tau_y^2=1$.
The following Lemma~\ref{lem-inv-bij} might be  known in some literature, but we could not find an explicit one. Therefore,  a proof is given below for reasons of  completeness and easy readings.
 \begin{lem}\label{lem-inv-bij} Let $G:=\PSU(3,q)$,
 $ {\cal I}:=\{g^2=1\mid g\in G\},
 $  $\t({\cal N})=\{\tau_{y}\mid  \lg y\rg \in {\cal N}\}$  and  ${\cal M}=\{ \PSO(3,q)^g\mid g\in G\}$.
  Then
\begin{itemize}[itemsep=0pt]
  \item[{\rm(1)}]  $\t: \lg y \rg\mapsto\t_{y}$ gives a bijection from  ${\cal N}$ to ${\cal I};$
  \item[{\rm(2)}]  For any $K\in{\cal M}$, $\tau_y\in K$ if and only if $\lg y\rg\in \mathbf{u}\cap {\cal N}$, where $G_{\mathbf{u}}=K;$ and
 \item[{\rm(3)}] If $\lg y\rg\neq\lg z\rg$, then $\tau_y\tau_z=\tau_z\tau_y$ if and only if $(y,z)=0$. Consequently,  $\lg \tau_x,\tau_y,\tau_z\rg\cong \D_4$
if and only if $\{\lg x\rg, \lg y\rg, \lg z\rg\}$ forms an orthogonal frame.
 \end{itemize}
  \end{lem}
\demo (1) Obviously,  if $\lg y \rg\neq\lg z \rg$, then $\tau_y\neq \tau_z$, where $\lg y \rg, \lg z \rg\in{\cal{N}}$. So $|\tau({\cal {N}})|\geq q^2(q^2-q+1)$, as $|{\cal N}|=q^2(q^2-q+1)$. Note that there is only one conjugacy class of involutions in $\PSU(3,q)$ \cite[Table B.1]{BGbook} and so $|{\cal {I}}|=\frac{|\PSU(3,q)|}{|\SL(2,q)\rtimes\ZZ_{\frac{q+1}{d}}|}=q^2(q^2-q+1)$, as $\N_G(\tau)\cong \SL(2,q)\rtimes\ZZ_{\frac{q+1}{d}}$ for any involution $\tau\in\PSU(3,q)$. Now, the map $\langle y \rangle \mapsto \tau_y$ from $\mathcal{N}$ to $\mathcal{I}$ is well-defined and injective, it follows that this map is in fact a bijection, as desired.
\vskip 3mm
(2) By the transitivity of $G$ on $\cal M$ under the conjugacy action, one may assume $K=\PSO(3,q)$ and $\mathbf{u}=\mathbf{w}_0=\PG(2,q)$ as defined above.
Without loss of generality, set $(x,y)=x_1y_1^q+x_2y_2^q+x_3y_3^q$, for any  $x=(x_1,x_2,x_3)$ and $y=(y_1,y_2,y_3)$ in $V$.
Let  $\lg y\rg:=\lg(y_1,y_2,y_3)\rg\in{\cal N}$.
 Then
 $\lg\tau_y(e_\imath)\rg\in\mathbf{w}_0$ with $1\leq \imath\leq 3$, as $\tau_y\in\PSO(3,q)$.
 From this, it follows that $y_\imath y_\jmath^q\in\FF_q$ for any $1\leq \imath,\jmath\leq 3$. This implies that $y=\lambda( y_1',y_2',y_3'),$ for some $\lambda\in\FF_{q^2}^*$ and $y_\imath'\in\FF_q$ with $1\leq \imath\leq 3$. Thus, $\lg y\rg\in\mathbf{w}_0$.

 Conversely, if $\lg y\rg\in\mathbf{w}_0\cap {\cal N}$, then $y=\lambda_1( y_1',y_2',y_3'),$ for some $\lambda_1\in\FF_{q^2}^*$ and $y_\imath'\in\FF_q$ with $1\leq \imath\leq 3$. For any $\lg x\rg\in\mathbf{w}_0$, $x$ can be written as $x=\lambda_2(x_1,x_2,x_3)$ for some $\lambda_2\in\FF_{q^2}^*$ and some $x_\imath\in\FF_q$ with $1\leq \imath\leq 3$. This deduces that $\lg\tau_y(x)\rg=\lg \tau_{\lambda_1^{-1}y}(x)\rg=\lg\lambda_2(x_1',x_2',x_3')\rg$ for some $x_1',x_2',x_3'\in\FF_q$, that is
 $\tau_y(\lg x\rg)\in\mathbf{w}_0$. Therefore, $\tau_y\in G_{\mathbf{w}_0}=\PSO(3,q)$.
\vskip 3mm
(3)  Check that $\tau_y\tau_z=\tau_z\tau_y$ if and only if
$(x,z)(z,y)y=(x,y)(y,z)z$,  for any $\lg x\rg\in\PG(2,q^2)$, which  holds if and only if $(y,z)=0$.
\qed
\vskip 3mm
For any $v\in V$,  define $v^\perp=\{v'\in V\mid (v,v')=0\}$ and $\PG(v^\perp)=\{\lg v' \rg\mid v'\in v^\perp \}.$ The following lemma is the key result for the proof of Theorem~\ref{main1}.

\begin{lem}\label{lem_tauyint}
Let $\mathbf{u}$ be a Baer subplane of $\mathcal{B}$ and assume that  $\langle v\rangle$ is an isotropic point of $\mathbf{u}$.  Let $\ell$ be a Baer subline that contains two other non-isotropic points $\langle y\rangle$ and $\langle  z\rangle$ outside $\mathbf{u}$. Then the following statements hold:
\begin{enumerate}
\item[{\rm(1)}] $\tau_y(\mathbf{\mathbf{u}})\cap \mathbf{\mathbf{u}}=\{\langle v\rangle\}$;
\item[{\rm(2)}] $\tau_y(\mathbf{u})\cap \tau_z(\mathbf{u})\subseteq \PG(v^{\perp})$;
\item[{\rm(3)}] $|\tau_y(\mathbf{u})\cap \tau_z(\mathbf{u})|= \begin{cases}
   q+1,\quad |\ell\cap\mathbf{u}|=1;\\
   1, \quad |\ell\cap\mathbf{u}|=2.
\end{cases}$
\end{enumerate}
\end{lem}
\demo
Recall that  $e_1=(1,0,0),\ e_2=(0,1,0)$ and $e_3=(0,0,1)$.
Then $\langle e_1\rangle\in\mathcal{H}$.
Without loss of generality, we may assume that $\mathbf{u}=\mathbf{w}_0$ by the transitivity of $G$ on $\mathcal{B}$. Moreover, since $\PSU(3,q)_{\mathbf{u}}\cong \PSO(3,q)$ acts transitively on all $q+1$ isotropic points of $\mathbf{u}$, we can further assume that $v=e_1$. In this case, $ \PG(e_1^{\perp})=\{\langle a e_1+b e_2\rangle : a,b\in\FF_{q^2}\}.$

\medskip
 Assume that  $\langle y\rangle $  is a non-isotropic point of $\PG(v^{\perp})\setminus \mathbf{u}$. Thus, we can set $
y=e_2+\alpha e_1$ for some $\alpha\in\FF_{q^2}\setminus\FF_q$.
A direct computation shows that
\[
\tau_y(e_1)=-e_1,\qquad
\tau_y(e_2)=e_2+2\alpha e_1,\qquad
\tau_y(e_3)=-e_3+2\alpha^q e_2+2\alpha^{q+1} e_1.
\]
It follows that
\begin{equation}\label{compu}
\tau_y(\mathbf{w}_0)
=
\{
\bigl\langle (-a+2b\alpha+2c\alpha^{q+1})e_1+(b+2c\alpha^q)e_2-c e_3\bigr\rangle
\mid\ a,b,c\in\FF_q,\ (a,b,c)\neq(0,0,0)
\}.
\end{equation}
Observe that a point $\langle w\rangle\in\tau_y(\mathbf{w}_0)\cap \mathbf{w}_0$ if and only if $b=c=0$. Hence, we deduce that $
\tau_y(\mathbf{w}_0)\cap\mathbf{w}_0=\{\langle e_1\rangle\}$. Then the claim (1) holds.

\medskip
 Without loss of generality, we can assume that $y=e_2+\alpha e_1$ and $ z=e_2+\beta e_1$ for some $\alpha,\beta\in\FF_{q^2}\setminus\FF_q$.
Set $
\Delta=\tau_y(\mathbf{u})\cap\tau_z(\mathbf{u}).$
Assume that $\langle u\rangle\in\Delta$. By the expressions of $\tau_y(\mathbf{w}_0), \, \tau_z(\mathbf{w}_0)$ given in (1),  we deduce that there exist  $(a_i,b_i,c_i)\in\FF_q^3$ for $i=1,2$ such that
\[
\begin{aligned}
\langle u\rangle
&=\Bigl\langle (-a_1+2b_1\alpha+2c_1\alpha^{q+1})e_1+(b_1+2c_1\alpha^q)e_2-c_1 e_3\Bigr\rangle\\
&=\Bigl\langle (-a_2+2b_2\beta+2c_2\beta^{q+1})e_1+(b_2+2c_2\beta^q)e_2-c_2 e_3\Bigr\rangle.
\end{aligned}
\]
Comparing the above expressions up to scalar, we deduce that $a_1=a_2,\, b_1=b_2=0$ when $c_1=c_2=0$ and there are other $q$ choices for $\lg u\rg$  if $\alpha-\beta \in \FF_q$; hence
\[
|\Delta|=
\begin{cases}
q+1, & \text{if }\alpha-\beta\in\FF_q,\\[2pt]
1,   & \text{if }\alpha-\beta\notin\FF_q.
\end{cases}
\]
In particular, $\Delta\subseteq\PG(e_1^{\perp})=\PG(v^{\perp})$, which proves (2).

\medskip
 Finally, we determine the Baer subline $\ell$.
Since $y-z=(\alpha-\beta)e_1$, the Baer subline $\ell$ containing $\langle e_1\rangle$, $\langle y\rangle$ and $\langle z\rangle$ is given by
\[
\ell
=
\{\langle e_1\rangle\}
\cup
\bigl\{
\big\langle (k\alpha -k\beta+\beta)e_1+e_2\big\rangle : k\in\FF_q
\bigr\}.
\]
A straightforward verification shows that
$
|\ell\cap\mathbf{w}_0|=
\begin{cases}
1, & \text{if }\alpha-\beta\in\FF_q,\\[2pt]
2, & \text{if }\alpha-\beta\notin\FF_q.
\end{cases}$

\f Comparing this with the above description of $\Delta$, we obtain the desired conclusion in (3).
\qed

\vskip 3mm
 Two Baer subplanes $\{\mathbf{w}_1, \mathbf{w}_2\}$ in $\mathcal{B}_0$ are said to be a {\it base-size set} if  $G_{\mathbf{w}_1} \cap G_{\mathbf{w}_2} = 1$.
\begin{lem}\label{char} Consider the action of $G=\PSU(3,q)$ on the set $\mathcal{B}_0$ defined above. Then  $b(G)=2$ and
two Baer suplanes $\{\mathbf{w}_1, \mathbf{w}_2\}$ in $\mathcal{B}_0$ form a base-size set if and only if $G_{\mathbf{w}_1}\cap G_{\mathbf{w}_2}$ contains no involution; and if and only if either
  $\mathbf{w}_1\cap \mathbf{w}_2=\emptyset$ or $\mathbf{w}_1\cap \mathbf{w}_2=\lg v\rg\in{\cal H}$. \end{lem}
\demo  From Theorem~\ref{suborbit} (whose proof is independent of this section) and Table \ref{number}, the first part is true.
 Therefore, by Lemma~\ref{lem-inv-bij}, $\{\mathbf{w}_1, \mathbf{w}_2\}$ is a base-size set if and only if $\mathbf{w}_1 \cap \mathbf{w}_2 \subseteq \mathcal{H}$ (where the inclusion allows for the empty set). Assume, for a contradiction, that $\{\lg \a_1\rg,\lg \a_2\rg\}\subseteq\mathbf{w}_1\cap \mathbf{w}_2\cap {\cal H}$. Then $\PG(\a_1^\perp)\cap\PG(\a_2^\perp)=\lg\b\rg\in\mathbf{w}_1\cap\mathbf{w}_2\cap{\cal N}$, as $|\PG(\a_\i^\perp)\cap\mathbf{w}_{\i}\cap{\cal H}|=1$ with $\i=1,2$. This implies that $\{\mathbf{w}_1, \mathbf{w}_2\}$ is not a base-size set.
\qed
\vskip 3mm
Let $\FF_{q^2}^*=\lg \xi \rg$ and $\FF_q^*=\lg \th \rg $. For any three row vectors $\a_1, \a_2, \a_3\in V$, denote by $\mat (\a_1, \a_2, \a_3)$ the $3 \times 3$ matrix whose rows are $\alpha_1, \alpha_2, \alpha_3$ in  order, and by $\det(\alpha_1,\alpha_2,\alpha_3)$ its determinant.  Furthermore, let $\iso(\a_1, \a_2, \a_3)$ be the Gram matrix with respect to the Hermitian form $(\cdot,\cdot)$ on $V$, i.e., the matrix $\big( (\alpha_\i, \alpha_\j) \big)_{1 \le \i,\j \le 3}$.
A diagonal matrix with diagonal entries $a_1,a_2,a_3$ will be denoted by $[a_1,a_2,a_3];$ and an anti-diagonal matrix with anti-diagonal entries $b_1,b_2,b_3$ reading from bottom left to top right will be denoted by $]b_1,b_2,b_3[$.
Recall that for any element $\hat{g} \in \mathrm{PGU}(3,q)$, we denote by $g \in \mathrm{GU}(3,q)$ a fixed preimage of $\hat{g}$ under the natural projection. We keep the notation $\mathbf{w}_0=\PG(2,q)=\PG(\lg e_1,e_2,e_3\rg_{\FF_q})$, where $e_1,e_2, e_3$ are unit vectors.
To obtain a new characterization of Baer subplanes for use in Sections 3--5, we define the following set:
 $$\mathcal{P}(W)=\{\mathbf{w}_\imath:=\PG(W_\imath)| W_\imath:=\lg \a_1,\a_2,\a_3 \rg_{\FF_q}, \,
\,  \a_\imath\in V, \,\iso (\a_1, \a_2, \a_3)\in \GL(3,q)\}.$$
\begin{lem}\label{wsigma} Set $\mathcal{B}:=\{\PG(W_0)^{\hat{g}}:=\PG(W_0^g)\mid g\in \GU(3,q)\}$.  Then
\begin{itemize}[itemsep=0pt]
\item[\rm(1)] $\mathcal{B}=\mathcal{P}(W);$
 \item[\rm(2)] Let $\mathbf{w}=\PG(W)\in \mathcal{B}$, where $W=\lg \b_1, \b_2, \b_3\rg_{\FF_q}$.  Then a point $\lg \d\rg$ in $\PG(V)$ is contained  in $\mathbf{w}$ if and only if  $(a\d, \b_\jmath)\in \FF_q$, where $\jmath=1, 2, 3$, for some $a\in \FF_{q^2}^*;$
     \item[\rm(3)] Let $W=\lg a\a_1,b\a_2,c\a_3\rg_{\FF_q}$ and $W'=\lg \a_1,\a_2,\a_3\rg_{\FF_q}$. Suppose $W=W'^g$ for some $g\in \GU(3,q)$. Then $\hat g\in \PSU(3,q)$ if and only if   $abc\in (\FF_{q^2}^*)^3$.
\end{itemize}\end{lem}
\demo (1) Suppose $\mathbf{w}_\imath\in \mathcal{B}$.  Then $\mathbf{w}_\imath=\mathbf{w}_0^{\hat{g}}$ for some $\hat{g}\in \PGU(3,q)$. Set  $W_\imath=W_0^{g}$, where $\mathbf{w}_\imath=\PG(W_\imath)$ and $\mathbf{w}_0=\PG(W_0)$.
Without loss of any generality, let the  Gram matrix  be the identity matrix $E=[1,1,1]$, under a basis.
Then $W_\imath=\lg e_1^g,e_2^g,e_3^g\rg_{\FF_q}$, where $\iso(e_1^g,e_2^g,e_3^g)=E$,  as $g$ preserves the $U$-form. So $\mathbf{w}_\imath\in \mathcal{P}(W)$.

Conversely, let $W_\imath:=\lg \a_1, \a_2, \a_3 \rg_{\FF_q} ,$ where  $\a_1, \a_2, \a_3\in V$ such that $\iso(\a_1, \a_2, \a_3)\in \GL(3,q)$. Then the $U$-subspace $W_\imath$  forms  orthogonal space over $\FF_q$ and so
 it contains two types of bases with the orthogonal Gram (diagonal) matrix $[1,1,1]$ and $[1,1,\th]$, where $\FF_q^*=\lg \th \rg $.  Moreover, $W_0$ contains two such type of bases. Therefore, there exists $g\in \GU(3,q)$ mapping  a base of $W_0$ to a base of same type of $W_\imath$. So $\mathbf{w}_0^{\hat{g}}=\mathbf{w}_\imath$.
\vskip 3mm
(2)  Suppose  $\lg \d\rg\in \mathbf{w}$.  Then $a\d=a_1\b_1+a_2\b_2+a_3\b_3$ where $a_1, a_2, a_3\in \FF_{q}$ and  $a\in \FF_{q^2}^*$, and so $(a\d, \b_\jmath)\in \FF_q$ with $\jmath=1,2,3$.
Conversely, suppose that   $(a\d, \b_\jmath)\in \FF_q$ for some $a\in\FF_{q^2}^*$, where $\jmath=1, 2, 3$. Write $a\d=a_1'\b_1+a_2'\b_2+a_3'\b_3$. Then $a_1'(\b_1, \b_\jmath)+a_2'(\b_2, \b_\jmath)+a_3'(\b_3, \b_\jmath)=t_\jmath$ with
$t_\jmath\in \FF_q$, where $\jmath=1, 2, 3$.
Solving this system of equations, we get $a_1', a_2', a_3'\in \FF_q$ so that $a\d\in W$ and $\lg \d\rg \in \mathbf{w}$.
\vskip 3mm
(3) Suppose that $(a\a_1,b\a_2,c\a_3)^T=Xg(\a_1, \a_2, \a_3)^T$,  for some $X\in \GL(3,q)$ and $g\in \GU(3,q)$. Then $abc=|X||g|$ and $|X|\in\FF_q^*$. It follows that $\hat g\in \PSU(3,q)$ if and only if $abc\in (\FF_{q^2}^*)^3$.
\qed

\begin{lem}\label{perp-intersect} Let $M, M'\leq\PGU(3,q)$ such that $M\cong M'\cong \SO(3,q)$.  Let $t_1,t_2$ and $t_3$ be involutions of $M$ such that $A=\lg t_1,t_2,t_3 \rg\cong D_4$ and $A\cap M'=1$. Then  there exists at least one involution $t'\in M'$ such that $[t', t_\imath]=1$ for some $\imath\in\{1,2,3\}$.
\end{lem}
\demo
Let $\mathbf{w},\mathbf{w}'\in\mathcal{B}$, where
$M=\PGU(3,q)_{\mathbf{w}}$ and $M'=\PGU(3,q)_{\mathbf{w}'}.$
By Lemma \ref{lem-inv-bij}, there exists a orthogonal frame $\{ \lg y_1\rg , \lg y_2\rg, \lg y_3\rg \}$  such that $t_\imath=\tau_{y_\imath}$, with $\imath=1,2,3$.
Let the Gram matrix be the identity matrix $[1,1,1]$. No loss, set $\mathbf{w}=\PG(2,q)$. By Proposition \ref{lem_obser}.(4), there exists only one point $\lg y_\imath'\rg\in\PG(y_\imath^\perp)\cap \mathbf{w}'$, where
\vspace{-5pt}$$\begin{array}{lll}
y_1'=(0,1,\xi_1),&
y_2'=\a(1,0,\xi_2),&
y_3'=\b(1,\xi_3,0),
\end{array}$$
with $\xi_\imath,\a,\b\in\FF_{q^2}^*$ and $\imath=1,2,3$.
By Lemma \ref{lem-inv-bij}, it is now sufficient to prove that $\{\lg y_1'\rg, \lg y_2'\rg, \lg y_3'\rg\}\cap {\cal N}\neq \emptyset.$
On the contrary,  assume that $\lg y_\imath'\rg\in{\cal H}$ for all $\imath$.
Then by Lemma \ref{wsigma}(1), we get that $(y_\imath',y_\imath')=0$ and $(y_\imath',y_\jmath')\in\FF_q$ for any distinct $\imath,\jmath\in\{1,2,3\}$.
This implies that $y_2'=\a(1,0,k_2\a^q\xi_1)$ and $y_3'=k\a(1,kk'\a^q,0)$,
where $k_2,k,k'\in\FF_q^*$ and $\a\in\FF_{q^2}^*$.
Therefore $\lg y_1',y_2',y_3' \rg_{\FF_q}\subseteq W'$, where $\PG(W')=\mathbf{w}'$. Then
$z:=k_2\a^{1+q}y_1'-y_2'=(-\a,k_2\a^{1+q},0)\in W'$ so that
 either $\lg -k_2^{-1}k^2k'z+y_3'\rg=\lg y_1\rg\in \mathbf{w}'$ or $\lg k_2^{-1}k^2k'z\rg=\lg y_3'\rg$. If $\lg y_1 \rg \in \mathbf{w}' \cap \mathcal{N}$, then this contradicts to the
 condition $A \cap M' = 1$. If $\lg k_2^{-1} k^2 k' z \rg= \lg y_3' \rg$, then the three points $\lg y_1' \rg, \lg y_2' \rg, \lg y_3' \rg$ in $\mathcal{H}$ are collinear. However, a Baer subline contains at most two isotropic points,  a contradiction too.
\qed
\vskip -5pt
\subsection{The proof of Theorem~\ref{main1}}\label{geoproof}
\f {\bf Proof:} Let $G = \PSU(3,q)$ and consider the unitary projective space $\PG(2,q^2)$ with Gram
 matrix $]1,1,1[$ (the anti-diagonal matrix).
Let  $\mathbf{w}_1$ and $\mathbf{w}_2$ be  distinct elements in  $\mathcal{B}_0$. Then there exists an isotropic  point $\lg w \rg \in \mathbf{w}_1\setminus \mathbf{w}_2.$  By  the transitivity of $G$ on $\mathcal{B}_0$ and that of $\SO(3,q)$ on $\mathbf{w}_0\cap {\cal H}$, there exists an element $\hat{g}\in G$ such that $\mathbf{w}_1^{\hat{g}}=\mathbf{w}_0$ and $\lg w \rg  ^{\hat{g}}=\lg e_1\rg$ with $e_1=(1,0,0)$. Thus, $\lg e_1\rg \in \mathbf{w}_0\setminus \mathbf{w}_2^{\hat{g}} .$
By the coming Lemma~\ref{base1}, there exists a point $\lg y\rg \in (\PG(e_1^\perp)\cap {\cal N})\setminus \mathbf{w}_0$ such that both  $\t_y(\mathbf{w}_0)\cap \mathbf{w}_0$ and $\t_y(\mathbf{w}_0)\cap \mathbf{w}_2^{\hat{g}}$  are either empty or an isotropic point. Then  $\t_y(\mathbf{w}_0)^{\hat{g}^{-1}}\cap \mathbf{w}_1$ and $\t_y(\mathbf{w}_0)^{\hat{g}^{-1}}\cap \mathbf{w}_2$  are either empty or an isotropic point. By Lemma~\ref{char}, both $\{ \t_y(\mathbf{w}_0)^{\hat{g}^{-1}},\mathbf{w}_1\}$ and $\{ \t_y(\mathbf{w}_0)^{\hat{g}^{-1}}, \mathbf{w}_2\}$ are base-size sets, which completes the proof.\qed
\vskip 2mm
Following the proof of Theorem \ref{main1}, the Gram matrix will be taken as $]1,1,1[$ throughout this section.
\begin{lem}\label{base1}  For any  $\mathbf{w}'$ in $\mathcal{B}_0 \setminus \{\mathbf{w}_0\}$
there exists $\lg y\rg \in  (\PG(e_1^\perp)\cap {\cal N})\setminus \mathbf{w}_0$ such that both
$\{\t_y(\mathbf{w}_0),\mathbf{w}_0\}$
and $\{\t_y(\mathbf{w}_0 ), \mathbf{w}'\}$ are base-size sets.
\end{lem}
\demo  Let $q\ge 7$. No loss of any generality, let $\mathbf{w}_0\setminus\mathbf{w}'$ contains  the isotropic point $\lg e_1\rg$.
   For  any  $\lg y\rg \in (\PG(e_1^\perp)\cap {\cal N})\setminus \mathbf{w}_0$,
 from   Lemma~\ref{lem_tauyint} we know  $\t_{y}(\mathbf{w}_0)\cap \mathbf{w}_0=\lg e_1\rg $, meaning $\{ \mathbf{w}_0, \t_y(\mathbf{w}_0)\}$  is a base-size set following Lemma~\ref{char}.  The lemma is proved provided we may  find such a $\lg y\rg $ such that $\t_y(\mathbf{w}_0)$ intersects  no nonisotropic point with  $\mathbf{w}'$.
\vskip 3mm
(1) {\it Find a set $\Lambda$ of non-isotropic points in $\PG(e_1^\perp)$ which contains  $\lg y\rg$ such that $\tau_y(\mathbf{w}_0)\cap \mathbf{w}'\cap\PG(e_1^\perp)=\emptyset.$}
\vskip 2mm
Since $\lg e_1\rg \not\in \mathbf{w}'$, by  Proposition~\ref{lem_obser}, we get
 $\PG(e_1^\perp) \cap  \mathbf{w}'=\{ \lg  e_2+\lambda_0 e_1\rg \}$,  for some $\lambda_0 \in \FF_{q^2}.$  Suppose that $\t_{y}(\mathbf{w}_0)$ contains this point $\lg e_2+\ld_0 e_1\rg $ too, where  $\lg y\rg\in(\PG(e_1^\perp)\cap {\cal N})\setminus \mathbf{w}_0$. No loss, write $y=e_2+\mu e_1$, where $\mu\in\FF_{q^2}$. Then we have  $\t _y(\lg z\rg )=\lg e_2+\ld_0 e_1\rg$ for some $z\in \mathbf{w}_0$, equivalently, $\t_y(\lg e_2+\ld_0 e_1\rg)\in \mathbf{w}_0.$ Since
$$\begin{array}{l}
 \t_y(e_2+\ld_0 e_1)=e_2+(2\mu -\ld_0)e_1,
 \end{array}$$
 that is $\mu \in\frac 12\ld_0 +\FF_q$, so that $y\in e_2+(\frac 12\ld_0 +\FF_q)e_1.$
 Based on the above observations, we introduce three sets:
$$\Delta=\FF_{q^2}\setminus (\FF_q\cup (\frac{1}{2}\ld_0+\FF_q)),\quad \Lambda=\{\lg e_2+xe_1 \rg\mid x\in\Delta\} \quad  {\rm and} \quad  {\cal W}=\{\t_y(\mathbf{w}_0)\mid  \lg y\rg \in \Lambda\}.$$
By Lemma \ref{lem_tauyint},  $\tau_y(\mathbf{w}_0)\neq\tau_z(\mathbf{w}_0)$ for distinct points $\lg y\rg,\lg z\rg\in\Lambda$.
Then $|{\cal W}|=|\Lambda |=|\Delta|=q^2-(1+\e)q,$ where $\e=0~{\rm or}~1$, depending on $\ld_0\in\FF_q$ or $\ld_0\in\FF_q^2\setminus\FF_q$.   For any $\t_y(\mathbf{w}_0)\in {\cal W}$, we have
 $\t_y(\mathbf{w}_0)\cap \mathbf{w}'\cap \PG(e_1^\perp)=\emptyset.$
\vskip 3mm
 (2) {\it Find $\t_{y}(\mathbf{w}_0)\in{\cal W}$ intersecting no nonisotropic points with $\mathbf{w}'$.}
\vskip 2mm
For the contrary, assume  that for any $\lg y\rg \in \Lambda$, we have
$\t_{y}(\mathbf{w}_0)\cap \mathbf{w}'$ contains at least one non-isotropic point.
There are totally $q^2+1$ lines through $\lg e_1\rg $  represented by
$$\ell_{\infty}=\PG(e_1^\perp)=\PG(\lg e_1, e_2\rg)  \, {\rm and}\, \ \ell_{x}=\PG(\lg e_1, e_3+xe_2\rg), \,   x\in \FF_{q^2},$$
 noting that the union of these lines cover all points of $\PG(2, q^2)$.
 Recall  $ \ell_{\infty}\cap \mathbf{w}' =\lg  e_2+\lambda_0 e_1\rg$ and moreover, $\PG((e_2+\lambda_0 e_1)^\perp)=\PG(\lg e_1, e_3-\lambda_0^q e_2\rg)=\ell_{-\ld_0^q }.$
 Since  $\lg  e_2+\lambda_0 e_1\rg \in \mathbf{w}'$, by  Proposition~\ref{lem_obser}(4)  we get  $|\ell_{-\ld_0^q }\cap \mathbf{w}'|=q+1.$
 Clearly, every point $\lg w\rg \in \mathbf{w}'$ is contained in the  line  $\PG(\lg e_1, w\rg)$, which must be   $\ell_x$ for some $x\in \FF_{q^2}\cup \{\infty\}$.
In what follows we shall carry out the proof by the following three steps.
  \vskip 2mm  {\it Step 1:  Drive an equation: \begin{equation}\label{eqn_claimF_q}
  \bigcup_{\lg y\rg \in \Lambda} (\t_{y}(\mathbf{w}_0) \cap  \mathbf{w}')=\bigcup _{x\in 2\Delta } (\ell_{x} \cap \mathbf{w}')
\end{equation}
 \vskip -2mm
\f where both sides contain exactly $q^2-(1+\e)q$ nonisotropic points of $\mathbf{w}'$; here $\e=0$ or $1$, depending on whether $\ld_0\in\FF_q$ or not.}
\vskip 3mm
Firstly we show that the left-hand side is a subset of the right-hand side in  Eq(\ref{eqn_claimF_q}).
In fact, for any $\lg y\rg\in\Lambda$, one may write  $y=e_2+\alpha e_1$, where $\a\in \Delta=\FF_{q^2}\setminus (\FF_q\cup (\frac{1}{2}\ld_0+\FF_q))$.
 We have shown that $\tau_y(\mathbf{w}_0)\cap \ell_{\infty }\cap \mathbf{w}'=\emptyset$, meaning that there exists no point in  $\tau_y(\mathbf{w}_0)\cap \mathbf{w}'$ contained in $\ell_{\infty }.$
 Suppose a point $\lg z\rg$ of $\tau_y(\mathbf{w}_0)\cap \mathbf{w}'$  is contained in $\ell_x$. Now, we only need to prove $x\notin \FF_q\cup (\ld_0+\FF_q)$, see the definition of $\Delta$. On the contrary, assume that $x\in \FF_q\cup (\ld_0+\FF_q)$.
 Then there exists $w=(a, b, -1)\in W_0$ with $(a,b\in\FF_q)$, and $r\in \FF_{q^2}$ such that
$$\lg\tau_y(w)\rg=\lg z\rg =\lg \big(-a +2b\alpha -2\alpha^{q+1})e_1+(b-2\alpha^q)e_2+e_3\rg =
\lg r e_1+xe_2+e_3\rg\in\ell_x,$$ which follows from Eq\eqref{compu}.
This forces $\alpha \in \FF_q\cup (\frac{1}{2}\ld_0+\FF_q)$, a contradiction.

\vskip 2mm
Secondly, we show that both sides of  Eq(\ref{eqn_claimF_q}) have  the same cardinality.
On the one hand, since  $|\ell_x \cap \mathbf{w}'|=1$ if and only if $x\ne -\ld_0^q$ (by Proposition \ref{lem_obser}, noting that $|\ell_{-\lambda^q}\cap \mathbf{w}'|=q+1$), we have $|\bigcup_{x\in 2\Delta }(\ell_x \cap \mathbf{w}')|=q^2-(1+\e)q=|\Lambda|$, for $\e=0~{\rm or}~1$ depending on $\ld_0\in\FF_q$ or $\ld_0\in\FF_{q^2}\setminus\FF_{q}$. On the other hand,
for any distinct $\lg y_1\rg, \lg y_2\rg \in \Lambda,$  by Lemma~\ref{lem_tauyint} we have
$$ \Pi:=(\t_{y_1}(\mathbf{w}_0) \cap  \mathbf{w}')\cap (\t_{y_2}(\mathbf{w}_0) \cap  \mathbf{w}')\subseteq   (\t_{y_1}(\mathbf{w}_0) \cap \t_{y_2}(\mathbf{w}_0))\cap \mathbf{w}' \subseteqq \ell_{\infty }\cap \mathbf{w}' \cap \t_{y_2}(\mathbf{w}_0)=\emptyset, $$
while $\Pi=\emptyset$ implies $|\bigcup_{\lg y\rg \in \Lambda }(\t_{y}(\mathbf{w}_0) \cap  \mathbf{w}')|\ge |\Lambda|$, as desired.

Moreover, under the assumption that for any $\lg y\rg \in \Lambda$,
$\t_{y}(\mathbf{w}_0)\cap \mathbf{w}'$ contains at least one nonisotropic point, we have $|\t_{y}(\mathbf{w}_0) \cap  \mathbf{w}'\cap{\cal N}|=1$.

\vskip 3mm
{\it Step 2:  Show that  $q+1$ isotropic points of $\mathbf{w}'$ are contained in $\bigcup _{x\in (\ld_0 + \FF_q)\cup\{\infty\}}  (\ell_{x} \cap \mathbf{w}').$}
\vskip 2mm
From the proof of Step 1,  one may see  that
  $\bigcup_{x\in 2\Delta}(\ell_x \cap \mathbf{w}')$ consists of  $q^2-(1+\e)q$ nonisotropic points, which implies that
    $${\small\mathbf{w}'\cap {\cal H}\subseteq\bigcup _{x\in \FF_q\cup (\ld_0+\FF_q)\cup \{\infty\}}  (\ell_{x} \cap \mathbf{w}').}$$
The statement of Step 2 is clearly true when $\ld_0 \in \FF_q$. So assume $\ld_0 \in \FF_{q^2}\setminus \FF_q$.
For the contrary, there exists an isotropic point, say $\lg z_2\rg $   contained in  $\bigcup _{x\in\FF_q}  (\ell_{x} \cap \mathbf{w}')$.
Then we may write  $z_2=\a_1 e_1+xe_2+e_3  $ with $x\in\FF_q$, where $\a_1+\a_1^q+x^2=0,$  as $(z_2,z_2)=0$.  Recall
 $ z_1:=\ld_0 e_1+e_2\in e_1^\perp\cap W',$
 with $\mathbf{w}'=\PG(W').$ Let $z_3 =\lg z_1,z_2 \rg\cap z_1^{\perp}$ be the intersecting point of two lines in $V(3,q^2)$, where $\PG(z_1^\perp)=\ell_{-\ld_0^q}$.
 So, we have $\lg z_3\rg\in\mathbf{w}'.$
 Then, no loss,
 there exist $m, s,t\in\FF_{q^2}$ such that
 \vspace{-5pt}$$z_3=mz_1+z_2=m(\ld_0 e_1+e_2)+(\a_1 e_1+xe_2+e_3)=se_1+t(-\ld_0^q e_2+e_3).$$
 Then we have  $t=1$, $m=-x-t\ld_0^q$ and  $s=\a_1-x\ld_0 -\ld_0^{q+1}$, that is  $z_3=z_2-(x+\ld_0^q)z_1.$
Set $\ell'$ be the unique Baer subline which contains the points $\lg z_\imath\rg$ with $\imath=1,2,3.$
Therefore,
$$\begin{array}{lll}
 \ell'&=&\{ \lg z_2+k(x+\ld_0^q)z_1\rg \mid k\in \FF_q\}\cup \{\lg z_1\rg \}\\
 &=&\{ \lg (\a_1 e_1+xe_2+e_3)+k(x+\ld_0^q)(\ld_0 e_1+e_2)\rg \mid k\in \FF_q\}\cup \{\lg\ld_0 e_1 +e_2 \rg\}\\
 &=&\{(p(k,x):=\lg (\a_1+k\ld_0^{q+1}+kx\ld_0)e_1+(x+kx+k\ld_0^q)e_2+e_3\rg \mid k\in\FF_q\}\\
 &&\cup \{ \lg e_2+\ld_0 e_1 \rg\},
 \end{array}$$
where $x\in\FF_q$.  Since $$x+kx+k\ld_0^q\in x+xk+k(-\ld_0 +\FF_q)=-k\ld_0 +\FF_q\in \FF_q\cup (\ld_0+\FF_q)\cup \{\infty\}$$ \vspace{-7pt} \f if and only if $k\in \{0, -1\}$, that is
$\lg p(k,x)\rg \in  \{\lg z_2\rg , \lg z_3\rg \}$.
Therefore,
$$\ell'\cap  \bigcup_{x\in \FF_q\cup (\ld_0 +\FF_q)\cup \{\infty\}}(\ell_x \cap \mathbf{w}')=\{ \lg z_1\rg, \lg z_2\rg , \lg z_3\rg \}.$$
 Note that $\ell'$ is a secant line of $\mathbf{w}'$, as $\lg z_2\rg$ is an isotropic point and $(z_2,z_1)\neq 0$. This forces  $\lg z_2\rg $ and $\lg z_3\rg $
are only two isotropic points of $\ell'$. This implies that $(z_3,z_3)=0$. Then
$$\begin{array}{lll}
(z_3,z_3)=-(x+\ld_0^q)^q(z_1,z_2)^q-(x+\ld_0^q)(z_1,z_2)+(x+\ld_0^q)^{1+q}(z_1,z_1)=0,
\end{array}$$
which forces $x+\ld_0^q=0$, as $(z_1,z_2)=\ld_0+x^q$ and $(z_1,z_1)=1$.
Therefore, $x=-\ld_0^q$, which is a contradiction to the fact that $x\in\FF_q$ and $\ld_0\notin\FF_q$, as desired.

\vskip 3mm  {\it Step 3: Get a contradiction.}
\vskip 3mm

From Step 2, all $q+1$ isotropic points of $\mathbf{w}'$ are contained in $\bigcup_{x\in (\ld_0+\FF_q)\cup\{\infty\}}  (\ell_{x} \cap \mathbf{w}').$
     We shall find two lines $\ell '$ and $\ell '' $ so that $(\ell ' \cup \ell'')\cap \mathbf{w}'=\bigcup_{x\in (\ld_0+\FF_q)\cup\{\infty\}}(\ell_{x} \cap \mathbf{w}')$. This will  arise   a contradiction, as both  $\ell' \cap \mathbf{w}'$ and  $\ell' \cap \mathbf{w}'$
  contain  at most two isotropic points.

   Firstly, set
   $\ell'= \PG((e_2+\lambda_0 e_1)^\perp)$. Then
$$\begin{array}{ll}
   &\ell'\cap \mathbf{w}'=\PG(\lg e_1, e_3-\ld_0^q e_2\rg)\cap \mathbf{w}'=\ell_{-\ld_0^q }\cap \mathbf{w}' \, \, (\in \bigcup_{x\in (\FF_q+\ld_0)\cup\{\infty\}} (\ell_{x} \cap \mathbf{w}')),
   \end{array}$$
and $|\ell'\cap \mathbf{w}'|=q+1$, as $\lg e_2+\lambda_0 e_1\rg\in \mathbf{w}'$.

   Secondly,  let  $\ell''$ to be the line containing $\lg z_1\rg:=\lg \lambda_0 e_1+e_2 \rg $ and other two isotropic points, say $\lg w_1\rg$ and $\lg w_2\rg $ of $\mathbf{w}'$.
    To find $\lg w_1\rg $ and $\lg w_2\rg $, noting that every  isotropic  point of $\mathbf{w}'$  is  contained in $\bigcup_{x\in (\ld_0+\FF_q)\cup\{\infty\}}(\ell_{x} \cap \mathbf{w}')$ shown as above,
    we write $ w_\imath=\alpha_\imath e_1+(a_\imath+\ld_0)e_2+e_3)\in \ell_{a_\imath+\ld_0}$, where $\imath\in \{1, 2\}$, $\alpha_\imath\in \FF_{q^2}$, $a_\imath\in \FF_q$. Since $\lg z_1\rg, \lg  w_1\rg$ and $\lg w_2 \rg$ are collinear, we get
{\small\[\left|\begin{array}{ccc}
   	\lambda_0  & \alpha_1 &\alpha_2 \\
   	1 & a_1+\ld_0 & a_2+\ld_0 \\
   	 0 & 1 & 1
   \end{array}\right|=0, \]}
  giving $\ld_0(a_1-a_2)=(\alpha_1-\alpha_2)$.
 Now we have that
$$ w_1- w_2= (\alpha_2-\alpha_1) e_1+(a_1-a_2)e_2=t(\ld_0 e_1+e_2 )=t z_1,\,  {\rm where} \, t:=a_1-a_2\in\FF_q^*,$$
 $$\begin{array}{lll}  \ell'' \cap \mathbf{w}' & = & \{\lg z_1 \rg\}\cup \{\lg w_2+ kt(\lambda_0 e_1+e_2) \rg\mid  k\in \FF_q\}\\
  &=&  \{\lg z_1 \rg\}\cup \{\lg \alpha_2e_1+(a_2+\ld_0)e_2+e_3)+ kt(\lambda_0 e_1+e_2) \rg\mid  k\in \FF_q\}
  \end{array}$$
 $$\begin{array}{lll}
  &=&\{\lg z_1 \rg\}\cup \{\lg (\alpha_2+kt\ld_0)e_1+(a_2+kt+\ld_0)e_2+e_3 \rg\mid  k\in \FF_q\}\\
  &=&\{\lg z_1 \rg\}\cup \{\lg (\alpha_2+(k'-a_2)\ld_0)e_1+(k'+\ld_0)e_2+e_3 \rg\mid  k'\in \FF_q\}\\
   &\subseteq & \bigcup _{x\in (\FF_q+\ld_0)\cup\{\infty\}} (\ell_{x} \cap \mathbf{w}')),
  \end{array}
 $$
 where $k':=a_2+kt$.
 Moreover, we get that $|\ell''\cap \mathbf{w}'|=q+1$.
Till now, we get two lines $\ell'$ and $\ell''$ such that $(\ell'\cup\ell'')\cap \mathbf{w}'\subseteq  \bigcup_{x\in (\FF_q+\ld_0)\cup\{\infty\}} (\ell_{x} \cap \mathbf{w}'))$. Since
both sides of this inequality has the cardinality $2q+1$, they are the same set, as desired.  \qed

\section{$\PSL(3,q)\lvertneqq G$ and $\SO(3,q)\lvertneqq M$}\label{ell}
 In this section, we mainly prove the following theorem.
\vskip -7pt  \begin{theorem} \label{main2}  Let $G=\PSU(3,q)\ZZ_{m_1}$ and
    $\a=M=\PSO(3,q)\ZZ_{m_1}\le G$, with  $q=p^m\ge 7$ for an odd prime $p$,  $m_1\mid 2m$ and $m_1\ne 1$.
Under the primitive right multiplication action of $G$ on $\Omega=[G:M]$, we have $b(G)=2$ and the BG-Conjecture is true, provided $\sqrt{q}\geq 17$ or $45$ if $3\nmid (q+1)$ and $3\di (q+1)$, respectively.
  \end{theorem}

  From now on, unless otherwise specified, we assume $d:=(3,q+1)$,  $q\geq 7$ and   $\sqrt{q}\geq 15d$.
  Let $T=\PSU(3,q)$, $M_0=\PSO(3,q)\cong \SO(3,q)$,  $G=\PSU(3,q)\ZZ_{m_1}$,
    $M=M_0\ZZ_{m_1}\le G$ where $m_1\mid 2m$ and $m_1\ne 1$.
Consider the primitive permutation representation of $G$ on $\Omega=[G:M]$ with right multiplication action, where $|T|=q^3(q^3+1)(q^2-1)/d$ and $|\O|=q^2(q^3+1)/d$.
Set $\a=M$.
Let $\Gamma_{nr}(\a)$ and $\Gamma_{r}(\a)$ be the union of nontrivial nonregular suborbits and of regular suborbits of $T$  relative to $\a$, respectively; and let
$\Gamma_{nr}'(\a)$ and $\Gamma_{r}'(\a)$ be the union of nontrivial nonregular suborbits   and of regular suborbits of  $G$  relative to $\a$, respectively.
Let $n'(\a)$ denote the number of points $\g \in\G_r(\a)$ that are fixed by some nontrivial element of $M\setminus M_0$.

The logical  structure of the proof of Theorem~\ref{main2} is  illustrated by the  diagram:
{\small $$\begin{array}{llllll}
            &\swarrow {\rm Lem~\ref{measure}} & & & &\\
           & \swarrow{\rm Lem~\ref{n'r}}      &\swarrow {\rm Th~\ref{suborbit}}   & \leftarrow {\rm Lem~\ref{fix}}  &&\\
{\rm Th~\ref{main2}}     & \leftarrow {\rm Lem~\ref{b(G)}}   &\leftarrow {\rm Th~\ref{main21}}   &  \leftarrow {\rm Lem~\ref{AB}}  &  \leftarrow {\rm Lem~\ref{used1}}  &  \leftarrow {\rm Lem~\ref{used0} }\\
                                                     &&  &\nwarrow \underline{{\rm Lem~\ref{C}}} &  \leftarrow \cdots\cdots &\\
                                                    &&&   \nwarrow {\rm Lem~\ref{D}}&   & \\
                                                      &&& \nwarrow {\rm Lem~\ref{DD}} \\
                    &\swarrow  {\rm Lem~\ref{Dqpm1}} &                             &                               &\swarrow {\rm Lem~\ref{ld12}}& \swarrow {\rm Lem~\ref{3-1}}            \\
                    &\swarrow  {\rm Lem~\ref{q+1}}&  \swarrow {\rm Lem~\ref{hlD}}    &  \swarrow {\rm Lem~\ref{H1_3}}   &  \leftarrow {\rm Lem~\ref{ld3}}    & \leftarrow {\rm Lem~\ref{absolute}}\\
\underline{{\rm Lem~\ref{C}}} &  \leftarrow {\rm Lem~\ref{B1B2}}    &  \leftarrow {\rm Lem~\ref{l-h}}  &   \leftarrow {\rm Lem~\ref{L_1(g)3}}  & \leftarrow {\rm Lem~\ref{n2d}} &\leftarrow {\rm Lem~\ref{D4Z2}}\\
                   &  \nwarrow {\rm Lem~\ref{C3}}    &  \nwarrow {\rm Lem~\ref{f2}}     &                             &  \nwarrow {\rm Lem~\ref{n2}}            &\\
\end{array}
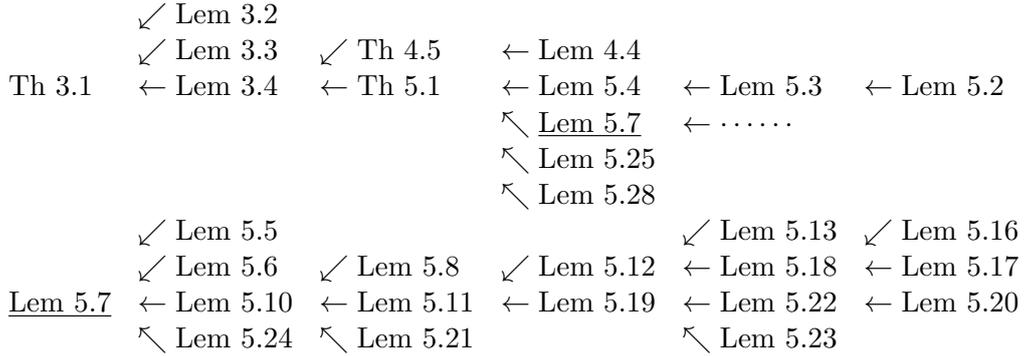
$$}
  As the first leave of the diagram,  Theorem~\ref{main2} is proved by   the following three lemmas.
\begin{lem} \label{measure}
The BG-Conjecture is true for $(G, M)$, provided that for any $\a, \a'\in \O$, at least one of the following two conditions is satisfied:
 \vspace{-5pt}\begin{equation}\label{lGM}
\ell(G, M):=|\Gamma_{nr}'(\a)\cap \Gamma_{nr}'(\a' )|-|\Omega|+2|\Gamma_{r}'(\a)|\gneqq 0;
\end{equation}
\begin{equation}\label{2GM}
2n'(\a) \lneqq \ell(T, M_0):=|\Gamma_{nr}(\a)\cap \Gamma_{nr}(\a')|-|\Omega|+2|\Gamma_{r}(\a)|.
\end{equation}

\end{lem}
\demo (1)  The BG-Conjecture  holds if and only if for any $\a,\a'\in\Omega$, \vspace{-5pt}
$$\begin{array}{ll}
&\Gamma_r'(\alpha)\cap \Gamma_r'(\a')\neq \emptyset\\
&\Leftrightarrow\Gamma_{nr}'(\alpha)\cup\Gamma_{nr}'(\a')\varsubsetneq \Omega\setminus\{\alpha,\a'\}
\Leftrightarrow|\Gamma_{nr}'(\alpha)|+|\Gamma_{nr}'(\a' )|-|\Gamma_{nr}'(\alpha)\cap\Gamma_{nr}'(\a' )|
\lvertneqq |\Omega|-2\\
&\Leftrightarrow |\Gamma_{nr}'(\alpha)\cap\Gamma_{nr}'(\a')|\gvertneqq  (|\Omega|-|\Gamma_r'(\alpha)|-1)+(|\Omega|-|\Gamma_r'(\alpha')|-1)-|\Omega|+2\\
&\Leftrightarrow\ell(G, M):=|\Gamma_{nr}'(\a)\cap \Gamma_{nr}'(\a' )|-|\Omega|+2|\Gamma_{r}'(\a)|\gneqq 0.

\end{array}
$$

(2) Some regular suborbits of $T$ might be fixed by some nontivial elements in $M\setminus M_0$ with $G_\a=M$,  meaning that they are  contained in some  nonregular suborbit of $G$.
 Therefore,  we have $|\Gamma_{r}'(\a)|\le |\Gamma_{r}(\a)|$ and $|\Gamma_{nr}'(\a)\cap \Gamma_{nr}'(\a')|\ge |\Gamma_{nr}(\a)\cap \Gamma_{nr}(\a')|$.
 Hence,
 $$\begin{array}{ll}
   &{\rm Eq}(\ref{lGM})\, {\rm holds} \\
    &\Leftrightarrow (|\Gamma_{nr}'(\a)\cap \Gamma_{nr}'(\a')|-|\Omega|+2|\Gamma_{r}(\a)|)
    +[(|\Omega|-2|\Gamma_{r}(\a)|)-(|\Omega|-2|\Gamma_{r}'(\a)|)]\gvertneqq  0\\
  &\Leftarrow      (|\Gamma_{nr}(\a)\cap \Gamma_{nr}(\a')|-|\Omega|+2|\Gamma_{r}(\a)|)-2(|\Gamma_{r}(\a)|-|\Gamma_{r}'(\a)|)\gvertneqq 0\\
  &\Leftrightarrow  \ell(T, M_0)\gvertneqq 2(|\Gamma_{r}(\a)|-|\Gamma_{r}'(\a)|)
  \Leftarrow     2n'(\a)\lneqq\ell(T, M_0). \hskip 112mm \Box
    \end{array} $$

  \begin{lem}\label{n'r}  Let $n'(r)$ be the number of regular suborbits of $T$, which are fixed by a subgroup $K\leq M\setminus M_0$ of prime order $r$.
  Then  $n'(r) \leq\frac{q_1^2}{d}+\frac{3}{2d},$ where $q_1=p^{\frac{m}{r}}$ with $r\ne 2$; and $n'(2)\leq q+\frac{5}{2}$ for  $d=1$ and  $\frac q2+\frac{1}{3}$  for $d=3$.
    \end{lem}
 \demo
 Suppose that $\a^{sM_0}$ is a regular suborbit of $T$ contained in a nonregular suborbit $\a^{sM}$ of $G$. Then $\a^{sM}$ is
 setwise fixed by a subgroup $K$, say $K=\lg uf'\rg $ of prime order $r$, where $u\in M_0$ and $f'=f^{\frac{m_1}{r}},$ where $m_1\mid 2m$ and $m_1\ne 1$. Clearly,  $q(q^2-1)\le |\a^{sM}|\lvertneqq m_1q(q^2-1)$.
In what follows we shall  derive  an upper bound of $n'(r)$, which depends on whether or not $r=2$.
\vskip 3mm
(1) Suppose that $r=2$. Then  $f': a\mapsto a^q$ and $f'\le \Z(M)$. In this case, $K\le M_1:=M_0\times \lg f'\rg$ and $K\cong \ZZ_2$.
Since $M_1$ (and so $M$) contains  five  conjugacy classes of subgroups $\ZZ_2$, say  $t_1, t_2\in M_0$, $f'$, $t_1f'$ and  $t_2f'$, we may set $K=\lg f'\rg $ or   $\lg t_\i f'\rg$  where $\i=1, 2$.

First, assume $d=1$. Now we compute
$\Fix(K)$, where $\Fix(K)$ denotes the number of fixed points of $K$ on $\O$.
By Proposition \ref{fieldconjugate}.(1),   $\lg f'\rg $, $\lg t_1f'\rg $ and $\lg t_2f'\rg $ are conjugate to each other in $\PSU(3,q)=T\le G$, as $d=1$.
By Proposition \ref{man}, for $\i=1,2$, we get
$${\small\begin{array}{ll}
&\Fix(\lg f' \rg)=\Fix(\lg t_\imath f' \rg)=\frac{|\N_G(\lg f' \rg)|}{|\N_{M}(\lg f' \rg)|}+\frac{|\N_G(\lg t_1 f' \rg)|}{|\N_{M}(\lg t_1 f' \rg)|}+\frac{|\N_G(\lg t_2 f' \rg)|}{|\N_{M}(\lg t_2 f' \rg)|}\\
&=\frac{|\PGO(3,q)|}{|\PSO(3,q)|}+\frac{|\PGO(3,q)|}{|\D_{2(q\pm 1)}|}+\frac{|\PGO(3,q)|}{|\D_{2(q\pm 1)}|}\leq1+\frac{q(q+1)}{2}+ \frac{q(q+ 1)}{2}=1+q^2+q.\end{array}}$$
Suppose $\lg f'\rg$ fixes a point in  $\a^{sM}$. Since $f'\le \Z(M)$, it would then fixes $\a^{sM}$ pointwise. But this leads to a contradiction, because $|\a^{sM}|>\Fix(\lg f'\rg).$
So $K=\lg t_\i f'\rg $, where $\i=1,2$.
Moreover, since  $K$ fixes the number $|\N_{M}(\lg t_\i f'\rg): \lg t_\i f'\rg|=\frac{2|\D_{2(q\pm 1)}|}{2}=2(q\pm1)$ of points in  $\a^{sM}$, where $\i=1,2$,
the number $n'(2)$ of such suborbits  $\a^{sM}$ is at most   $2\times\frac{q^2+q+1}{2(q-1)}=\frac{q^2+q+1}{q-1}=q+2+\frac{3}{q-1}\leq q+\frac{5}{2},$ as $q\geq 7$.

Second, assume $d=3$.
Since $\N_{\PGammaU(3,q)}(\lg f' \rg)=\PGO(3,q)\rtimes\lg f \rg=\PSO(3,q)\rtimes\lg f\rg=\N_{\PSigmaU(3,q)}(\lg f' \rg)$ and $|\PGammaU(3,q):\PSigmaU(3,q)|=d$,  we have that $\lg f'\rg$, $\lg t_1f' \rg$, $\lg t_2f'\rg$ are not conjugate to each other in $\PSigmaU(3,q)$ and all of them have normalizer $\PGO(3,q)\rtimes\lg f \rg$ in $\PSigmaU(3,q)$.
By Proposition~\ref{man} again,  we have
$${\small\begin{array}{ll}
\Fix(\lg f' \rg)=\frac{|\N_G(\lg f' \rg)|}{|\N_{M}(\lg f' \rg)|}=\frac{|\PGO(3,q)|}{|\PSO(3,q)|}=1,&
\Fix(\lg t_\imath f' \rg)=\frac{|\N_G(\lg t_\imath f' \rg)|}{|\N_{M}(\lg t_\imath f' \rg)|}=\frac{|\PGO(3,q)|}{|\D_{2(q\pm 1)}|}=\frac{q(q\mp 1)}{2} \end{array}},$$
with $\imath=1,2.$ Thus, $\lg f'\rg $ just fixes one point, that is $\a$, and so $K=\lg t_\imath f'\rg$. With the same arguments as in the last paragraph, $K$
fixes the number $2(q\pm1)$ of points in each   $\a^{sM}$.
Therefore, $n'(2)\leq \frac{q(q-1)}2/2(q+1)+ \frac{q(q+1)}2/2(q-1)\le \frac q2+\frac{1}{3},$ as $q\ge 7$.
 \vskip 3mm

(2) Suppose that $r\neq 2$. By \cite[Lemma 3.1.17]{BGbook}, $uf'$ is conjugate to $f'$ by an element of $M_0$. Therefore, $K=\lg uf'\rg$ is conjugate to $\lg f' \rg$ in $M$.  Set $q_1:=p^{\frac{m}{r}}$. Then by Proposition~\ref{man}, we have

$${\small\Fix(K)=\Fix(\lg f' \rg)=\frac{|\N_G(\lg f' \rg)|}{|\N_M(\lg f' \rg)|}=\frac{|\PSU(3,q_1)|}{|\PGL(2,q_1)|}
=\frac{q_1^3(q_1^3+1)(q_1^2-1)}{dq_1(q_1^2-1)}=\frac{q_1^2(q_1^3+1)}{d}.
}$$
 Since $\lg f' \rg$ fixes the number {\small$|\N_{M_0:\lg f' \rg}(\lg f' \rg):\lg f'\rg|=q_1(q_1^2-1)$} of points in each $\a^{sM}$, we get
$n'(r)\leq \frac{\Fix(\lg f' \rg)}{q_1(q_1^2-1)}=\frac{q_1(q_1^3+1)}{d(q_1^2-1)}= \frac{1}{d}q_1^2+\frac{q_1}{d(q_1-1)}\leq\frac{q_1^2}{d}+\frac{3}{2d},$ as $q_1\ge 3$.
\qed

 \begin{lem} \label{b(G)}
  Let $\PSU(3,q)\lneqq G\leq \PSigmaU(3,q)$. Then $b(G)=2$  for any $q\ge 7$ and $2n'(\a)< \ell(T, M_0)$ holds provided $\sqrt{q}\geq 17$ or $45$ if $3\nmid (q+1)$ and $3\di (q+1)$, respectively.\end{lem}
   \demo  Set $m_1=2^{e_0}m'$, where $m'=r_1^{e_1} r_2^{e_2}\cdots r_h^{e_h}$, where $r_1<r_2 \cdots <r_h$ are odd primes, $e_0\geq0$ and $e_\imath\ne 0$ for $\imath\geq 1$. Then $m'\ge r_15^{h-1}$ for $h\geq 1$, that is $h\le 1+ \log_5 \frac{m'}{r_1}.$
Let $n'$ be the number of regular suborbits of $T$ which are fixed by a subgroup of prime order in $M\setminus M_0.$
If $m'=1$, then $n'=n'(2)$ and $n'(\a)=n'(2)|M_0|\leq q^3n'(2)$, where $n'(2)$ is given in Lemma~\ref{n'r}; and if $m'\neq 1$,
then by  Lemma~\ref{n'r} we get  \vspace{-5pt}$$\begin{array}{l} n'\leq  n'(2)+\sum\limits_{\tiny{\begin{array}{l} r\mid  m', r\geq3,  \\ r\, {\rm is\, an\, odd\, prime}\end{array} }} n(r)\le  n'(2)+\sum\limits_{\tiny{\begin{array}{l} r\mid  m', r\geq3,  \\ r\,{\rm is\, an\, odd\, prime}\end{array} }}(\frac{1}{d} p^{\frac{2m}{r}}+\frac{3}{2d})
 \leq n'(2)+(\frac{q}{d}+\frac{3}{2d})h,\end{array} $$
\vspace{-8pt}
so that   $n'(\a)=n'|M_0|\le   q^4[n'(2)/q+(\frac{1}{d}+\frac{3}{2dq})(\log_5\frac{m'}{r_1}+1)],$
where $1\leq h\leq\log_5\frac{m'}{r_1}+1.$
\vskip 3mm
 Theorem~\ref{suborbit} will establish the lower bound  $|\Gamma_r(\a)|\ge  \frac{1}{3d}q^5-\frac{1}{3d}q^4$. Combining this lower bound for $|\Gamma_r(\a)|$ with the upper bounds for $n'(2)$ and $n'(\a)$ obtained above,
  one can verify that $\Gamma'_r(\a)\geq\Gamma_r(\a)-n'(\a)\gneq  0,$ implying  $b(G)=2$.
\vskip 2mm
 It will be shown in Theorem~\ref{main21} that for $\sqrt{q}\geq 15d$, we get
$$ {\small |\Gamma_{nr}(\a)\cap \Gamma_{nr}(\a' )|\ge \left\{
                         \begin{array}{ll}
                           \frac{4}{9}q^5-\frac{1}{3}q^4\sqrt{q}-22.1555q^4, & \hbox{$d=1$;} \\
                           \frac{4}{27}q^5-\frac{31}{27}q^4\sqrt{q}-16.8372q^4, & \hbox{$d=3$.}
                         \end{array}
                       \right.}
$$
With the  above upper bound of $n'(\a)$ and lower bounds of $|\Gamma_{r}(\a)|$ and $|\Gamma_{nr}(\a)\cap \Gamma_{nr}(\a')|$,
one may check that $2n'(\a)\lneq  \ell(T, M_0):=|\Gamma_{nr}(\a)\cap \Gamma_{nr}(\a')|-|\Omega|+2|\Gamma_{r}(\a)|$ for all cases with either $\sqrt{q}\geq 17$ and $d=1$, or $\sqrt{q}\geq 15d$ and $d=3$.  This proves the lemma.
\qed
\vskip 3mm  In Sections 4 and 5,  we shall   find lower bounds for $|\Gamma_{r}(\a)|$ and   $|\Gamma_{nr}(\a)\cap \Gamma_{nr}(\a' )|.$
\vspace{-3pt} \section{A lower bound of $|\Gamma_{r}(\a)|$}
Throughout this section,  set $G=T=\PSU(3,q)$ and $M=M_0=\PSO(3,q)$, where $q=p^m\ge 7$. Let $K=M_{\b}$ for some $\b\in \O\setminus\{\a\}$, where $\a=M$.  Then we shall determine
the suborbits of $G$ relative to $\a$ so that  $|\Gamma_{r}(\a)|$ is measured.

Throughout this paper, set $\FF_{q^2}^*=\lg \xi\rg$, $\FF_{q}^*=\lg \th\rg $ and $i^2=\th$ for some $i\in\FF_{q^2}$.
   Let   $f :\xi\mapsto \xi^p$ be the  automorphism of $\FF_{q^2}$ so that $|f|=2m$.
   Let  $\delta=[\mu , 1, 1]Z$ and $T=\PSU(3,q)$, where $|\d|=d=1, 3$.  Then  $\PGU(3,q)=T\rtimes \lg \d\rg $. The field automorphism of $T$ induced by $f$ is redenoted by $f$ again so that
 $\Aut(T)=T\rtimes \lg \d, f\rg .$ The first three propositions will be  used later.
 \begin{prop}{\rm(Tabel 8.5 and Tabel 8.6 in \cite{Bray})} \label{maximal}
  Every maximal subgroup of $\PSU(3,q),$  where $q=p^m$ for an odd  prime $p$, is isomorphic to one of the following groups:
{\rm(1)}  $[q^{1+2}]\rtimes\frac{q^2-1}{d};$
{\rm(2)} $\SL(2,q)\rtimes\ZZ_{\frac{q+1}{d}};$
{\rm (3)}  $(\ZZ_{\frac{q+1}{d}}\times \ZZ_{q+1})\rtimes\SS_3$ where $q\neq 5;$
{\rm (4)} $\ZZ_{\frac{q^2-q+1}{d}}\rtimes\ZZ_3$ where $q\neq 3,5;$
{\rm (5)} $\PSU(3,q').(\frac{q+1}{q'+1},3)$ where $q=q'^r$, $r$ odd prime$;$
{\rm(6)} $\SO(3,q)$ where $q$ odd and $q\geq7;$
 {\rm(7)} $\ZZ_3^2\rtimes Q_8$ if $3\mid\mid (q+1)$, $q=p\equiv 2(\mod 3)$ and $q\geq 11;$
  {\rm(8)} $\ZZ_3^2\rtimes\SL(2,3)$ if $9\mid (q+1)$, $q=p\equiv 2(\mod 3)$ and $q\geq 11;$
  {\rm (9)} $\PSL(2,7)$ where $q=p\equiv 3,5,6(\mod 7)$ and $q\neq 5;$
  {\rm (10)}
 $\A_6$ where $q=p\equiv 11,14 (\mod 15);$
  and {\rm (11)}
 $A_7$ and  $M_{10}$, where $q=5$.
\end{prop}
 \begin{prop}\label{fieldconjugate} {\rm\cite[p.77,79]{BGbook}}
Let $T=\PSU(3, q)$ and $f$ the field automorphism. Then
\vspace{-5pt}\begin{itemize}[itemsep=0pt]
\item[\rm(1)]  Let  $x=f^\i\in\lg f\rg$  such that  $x$ is of order a prime $r$.  Then for any $g\in \PGU(3,q)$,
   $|gx|=r$ if and only if $gx=x^{g_1}$  for some $g_1\in \PGU(3,q)$.
\item[\rm(2)] Let $z\in\PGU(3,q)$ be an element of odd prime order $r\ne p$.  Then $z^T=z^{\PGU(3,q)}$.
\end{itemize}
\end{prop}
\begin{prop}\label{man}
{\rm \cite{Man}} \
Let $G$ be a transitive group on $\O$ and let $H=G_\a$
for some $\a\in \O$.
Suppose that $K\le G$ and at least one $G$-conjugate of
$K$ is contained in $H$. Suppose further that
$G$-conjugates of $K$ which are contained in $H$ form $t$
conjugacy classes of $H$ with representatives $K_1, \cdots , K_t$.
Then $K$ fixes $\sum_{\i=1}^{t}|\N_G(K_\i):\N_H(K_\i)|$ points of $\O$.
\end{prop}

Take the Gram matrix  $J=]1,1,1[$ (anti-diagonal matrix).
In $\SU(3,q)$, set
{\small$$\begin{array}{lll}
q(a,b)={{\left(
\begin{array}{ccc}
1&a&b \\
0&1&-\overline{a} \\
0&0&1\\
\end{array}
\right)  }},&  h(k)={{ \left(
\begin{array}{ccc}
k^{-q}&0&0 \\
0&k^{q-1}&0\\
0&0&k\\
\end{array}
\right)  }}, &\, \tau={{ \left(
\begin{array}{ccc}
0&0&1 \\
0&-1&0\\
1&0&0\\
\end{array}
\right)  }},
\end{array}$$}
\vspace{-10pt}$${\small Q:=\{ \hat{q}(a,b)\mid a,b\in \FF_{q^2}, b+\overline{b}+a\overline{a}=0 \}, \quad H:=\{\hat{h}(k)\mid k\in\FF_{q^2}^*\}.}$$
Then $Q$ is a Sylow $p$-subgroup of $G=\PSU(3,q)$ of order $q^3$ and $H$ is a cyclic subgroup of $G$ of order $\frac{q^2-1}{d}$,
with relations
 $${\small\begin{array}{lll}
 &\hat{q}(a,b)\hat{q}(a',b')=\hat{q}(a+a',b+b'-a\overline{a'}), \quad
&\hat{\tau} \hat{h}(k)\hat{\tau}=\hat{h}(k^{-q}),\\
&\hat{h}(k)^{-1}\hat{q}(a,b)\hat{h}(k)=\hat{q}(k^{2q-1}a,k^{q+1}b),
&\hat{q}(a,b)^{-1}=\hat{q}(-a,\overline{b}).\\
\end{array}}
$$
Recall $M=\PSO(3,q)\cong \PGL(2,q)$. It is well-known that for odd $q > 3$, the maximal subgroups of $\PGL(2,q)$ are isomorphic to $\PSL(2,q)$; $\PGL(2,q')$ (where $q$ is a prime power of $q'$); $\ZZ_p^m \rtimes \ZZ_{q-1}$; $\D_{2(q \pm 1)}$; or $\SS_4$ (where $q \equiv \pm 3 \pmod{8}$).
Consider the action of $G$ by right multiplication on the set of right cosets $\O:=[G:M]$.
 The following key  lemma concerning their structure is used  to find the possible subgroups
$K=M_\b\leq M=G_{\a}\in \O$ with $\b\in \O \setminus \{\a\}$.
\begin{lem}\label{fix}
Let $G=\PSU(3,q)$ and $G_{\a}=M=\SO(3,q)$, where $q=p^m\geq 7$.
Then
\vspace{-5pt}\begin{itemize}[itemsep=0pt]
\item[\rm(1)] Let $P\in \Syl_p(M)$  and $P_1\le P$ where $|P_1|=p$. Then {\small$\Fix(P_1)=\Fix(P)=q$.}
\item[\rm(2)] Let $P'\cong \ZZ_p^{m'}$ where $m'\leq m$ and $A'\le \N_G(P')$.
 Then {\small$\Fix(P'\rtimes A')=1$} and $q$, if $|A'|\ge 3$ and $|A'|=2$, respectively.
\item[\rm(3)] Let $L\le D\le M$, where $L\cong \ZZ_l$ and $D\cong \D_{2l}$, where  $l>2$. Then {\small$\Fix(L)=\Fix(D)=\frac{q+1}{d}$.}
\item[\rm(4)] $M$ contains totally  two conjugacy classes of subgroups $\D_4$ and {\small$\Fix(D)=\frac{(q+1)^2}{d}$} for any $\D_4\cong D\leq M$.
\item[\rm(5)] $M$ contains totally two conjugacy classes of subgroups $\ZZ_2$ and {\small$\Fix(I)=\frac{q^2(q+1)}{d}$} for any $\ZZ_2\cong I\leq M$.
\item[\rm(6)]  {\small$\Fix(\SS_4)=\Fix(\A_5)=1$}; and {\small$\Fix(\A_4)=1, 3$}  for  $d=1, 3$, respectively.
\end{itemize}
\end{lem}
\demo
(1) Without loss of any generality, let $P:=\lg \hat{q}(a,b)\mid a,b\in \FF_q , 2b+a^2=0\rg\cong \ZZ_p^m\in \Syl_p(M)$. Fix an element $\hat{q}(a_1,b_1)\in P,$ and define $P_1:=\lg \hat{q}(a_1,b_1)\rg \cong \ZZ_p$.
Since $\Z(Q)=\{\hat{q}(0,b)\mid b\in\FF_{q^2}, b+\overline{b}=0\}$ where $Q$ was defined as above, it follows that
$P\setminus\{\hat{q}(0,0)\}\subseteq Q\setminus \Z(Q)$.
For any $\hat{q}(x_1,y_1),\hat{q}(x_2,y_2)\in Q$ and $\hat{h}(k)\in H$, we have
\vspace{-5pt}$${\small\begin{array}{ll} &\hat{q}(x_2,y_2)^{-1}\hat{q}(x_1,y_1)\hat{q}(x_2,y_2)=\hat{q}(x_1,x_1^qx_2-x_1x_2^q+y_1);\\
&\hat{h}(k)^{-1}\hat{q}(x_1,y_1)\hat{h}(k)=\hat{q}(k^{2q-1}x_1,k^{q+1}y_1).\end{array}}$$
Therefore, $\N_G(P)=\lg \hat{q}(x,y),\hat{h}(k) \mid y+\overline{y}+x^2=0, x\in \FF_q, y\in\FF_{q^2}, k\in\FF_q^*\rg$
is of order
$q^2(q-1)$ and $Q_1:=\lg \hat{q}(x,y)\mid y+\overline{y}+x^2=0, x\in \FF_q, y\in\FF_{q^2}\rg\leq\C_G(P)$.
By Proposition~\ref{man},
{\small$\Fix(P)=|\N_G(P):\N_M(P)|=\frac{q^2(q-1)}{q(q-1)}=q$,}
as $\N_M(P)\cong\ZZ_p^m:\ZZ_{q-1}$.
Suppose $\hat{h}(k) \in \N_G(P_1)$. Then $k\in\FF_{q}^*$ by the given condition. For such $k$
we have $\hat{h}(k) \in \N_M(P_1)$.
Consequently, {\small$\Fix(P_1)=|\N_G(P_1):\N_M(P_1)|=\frac{q^2}{q}=q$,}
and so $\Fix(P_1)=\Fix(P)=q$.
\vskip 3mm
(2)  Let $H_{q}:=\lg \hat{h}(k')\mid k'\in \FF_q^* \rg$ and $P'$ a subgroup of $P$ of order $p^{m'}$, where $0\lvertneqq m'\leq m$. Then $\N_M(P'):=P'\rtimes H'\cong\ZZ_p^{m'}\rtimes\ZZ_{p^{m'}-1}$.
By the proof of (1), we have $\N_G(P'\rtimes A')\leq \N_G(P')\leq Q_1\rtimes H_q$ for any subgroup $1\lneqq A'\leq H'$, where $Q_1$ is defined in the proof of (1). Let $\hat{q}(x_1,y_1)\in Q$ and $\hat{h}(k_1)\in H$. Then
 $$\begin{array}{l}
\hat{q}(x_1,y_1)^{-1}\hat{h}(k_1)\hat{q}(x_1,y_1)=\hat{g},~ \text{with}~{\small g={{\left(
\begin{array}{ccc}
k_1^{-1}&x_1k_1^{-1}-x_1&k_1^{-1}y_1+x_1x_1^q+k_1y_1^q \\
0&1&-x_1^q+k_1x_1^q \\
0&0&k_1\\
\end{array}
\right)  }}.}
\end{array}$$
If $y_1\in\FF_{q^2}\setminus\FF_q$ and $\hat{q}(x_1,y_1)\in \N_G(P)$, then $\hat{q}(x_1,y_1)^{-1}\hat{h}(k_1)\hat{q}(x_1,y_1)\in M$ if and only if $k_1^2=1$.
 Now, following the proof of (1), we conclude that $\N_G(P'\rtimes\lg \hat{h}(k_0)\rg)=Q_1\rtimes Y$ where $k_0^2=1$ for some  $Y\leq M$; and $\N_G(P'\rtimes A')=P\rtimes X\leq M$ where $A'\le H'$ with $|A'|\geq 3$ for some  $X\leq M$.
So
{\small$\Fix(P'\rtimes A')=1$} and $q$, if $|A'|\ge 3$ and $|A'|=2$, respectively.
\vskip 3mm
(3)  First suppose    $\D_{2(q\pm1)}\cong D\leq M.$  Since $M$ contains  only one  class of subgroups $\D_{2(q\pm1)}$, where $\N_M(D)=D$ and $\N_{G}(D)\cong\ZZ_{\frac{q+1}{d}}\times\D_{2(q\pm1)}$, we get
{\small$\Fix(D)=\frac{|\N_G(D)|}{|\N_M(D)|}=\frac{q+1}{d}.$}

Secondly suppose $\ZZ_l\cong L\leq D\cong\D_{2(q\pm1)}$ with $l\mid q\pm1$ and $l\ge 3$. Let $r$ be a prime divisor of $l$ and $R\in \Syl_r(L)$  so that $L\le C_G(R)$.   Checking Proposition~\ref{maximal}, we know that $G$ contains a maximal subgroup
$W=\SL(2, q)\rtimes \ZZ_{\frac{q+1}d}$ of order $q(q-1)(q+1)^2/d$. Therefore, $W$ contains a  Sylow $r-$subgroup of $G$ and then one may assume $R\le L\le W$ up to conjugacy.
Then $\N_W(L)\cong \ZZ_{\frac{q+1}{d}}\times (\ZZ_{q\pm1}\rtimes\ZZ_2)$, which is maximal in $W$. Check that every maximal subgroup of $G$ containing a subgroup $\ZZ_{\frac{q+1}{d}}\times (\ZZ_{q\pm1}\rtimes\ZZ_2)$ is conjugate to  $W$.
Therefore $\N_G(L)=\N_W(L).$ Finally, we have {\small$\Fix(L)=\frac{|\N_{G}(L)|}{|\N_M(L)|}=\frac{{2(q+1)(q\pm1)}/{d}}{2(q\pm1)}=\frac{q+1}{d}.$}

From the above arguments, we get {\small$\Fix(L)=\Fix(D)$}, which implies {\small$\Fix(D_{2l})=\Fix(L)$} for any $L\le D_{2l}\le D$, where $l\ge 3$.
\vskip 3mm
 (4)   Note that $M$ contains totally  two conjugacy classes of subgroups $\D_4$, with representatives  $\D_4^A$ and $\D_4^B$,  where $\D_4^A\leq \PSL(2,q)$.
Recall that $\N_M(\D_4^A)\cong \SS_4$ and $\N_M(\D_4^B)\cong \D_8$. By checking Proposition~\ref{maximal},  we easily get that $D_4^A$ and $\D_4^B$ are conjugate  in $G$ and $\N_G(\D_4^A)\cong\N_G(\D_4^B)\cong \ZZ_{\frac{q+1}{d}}\times \ZZ_{q+1}\rtimes\SS_3$. This implies that
{\small$\Fix(\D_4^A)=\Fix(\D_4^B)=\frac{|\N_G(\D_4^A)|}{|\N_M(\D_4^A)|}
+\frac{|\N_G(\D_4^B)|}{|\N_M(\D_4^B)|}=\frac{(q+1)^2}{d}.$}

\vskip 3mm
(5) Note that $M$ contains totally  two conjugacy classes of subgroups $\ZZ_2$, with representatives  $\ZZ_2^A$ and $\ZZ_2^B$,  where $\ZZ_2^A\leq \PSL(2,q)$.
 Recall that $\N_M(\ZZ_2^A)\cong\D_{2(q-\epsilon)}$ and $\N_M(\ZZ_2^B)\cong\D_{2(q+\epsilon)}$, where $q\equiv \epsilon\pmod 4$. Moreover, $G$ contains only one class of involutions such that
 $\N_G(\ZZ_2^A)\cong \N_G(\ZZ_2^B)\cong \SL(2,q)\rtimes\ZZ_{\frac{q+1}{d}}$. Therefore
{\small$\Fix(\ZZ_2^A)=\Fix(\ZZ_2^B)=\frac{|\N_G(\ZZ_2^A)|}{|\N_M(\ZZ_2^A)|}
+\frac{|\N_G(\ZZ_2^B)|}{|\N_M(\ZZ_2^B)|}=\frac{q(q^2-1)(q+1)}{2d(q+1)}+\frac{q(q^2-1)(q+1)}{2d(q-1)}=\frac{q^2(q+1)}{d}.$}

\vskip 3mm
(6) Suppose $L\cong \A_4$. Note that $\N_M(L)\cong \SS_4$ and there is only conjugacy class subgroup of $\A_4$ in $M$. Checking Proposition~\ref{maximal}, we get $\N_G(L)\lesssim \ZZ_{\frac{q+1}{d}}\times\ZZ_{q+1}\rtimes\SS_3$.
 Clearly, $L\cong \ZZ_2\times \ZZ_2\rtimes\ZZ_3$ and $\SS_4\cong\ZZ_2\times \ZZ_2\rtimes\SS_3\lesssim \N_G(L)$.
 Now, we construct the subgroups of $\PSO(3,q)$ which is isomorphic to $\SS_4$, under the Gram matrix $[1,1,1].$
Set $D:=\{\hat{z}(b,c)\mid z(b,c):=[b,c,(bc)^{-1}],b^{q+1}=c^{q+1}=1, b,c\in\FF_{q^2}\}$, $x:=[1,\omega,\omega^2]$ ($\omega:=\xi^{\frac{q^2-1}{3}}$), $y:={\tiny\left(
\begin{array}{ccc}
0 & 0 & 1 \\
1& 0 & 0 \\
0 & 1 & 0 \\
\end{array}
\right)}
,
\text{and}\,\, s:={\tiny\left(
\begin{array}{ccc}
1 & 0 & 0 \\
0& 0 & 1\\
0 & 1 & 0 \\
\end{array}
\right)}.$
One may check that $\hat{x}\in D$ only if $d=3$ with $d:=(3,q+1)$, $D\cong\ZZ_{\frac{q+1}{d}}\times\ZZ_{q+1}$,
$R=\C_G(\lg\hat{x}\rg)=D\rtimes\lg \hat{y}\rg\cong D\rtimes\ZZ_3$
and $D\rtimes\lg \hat{y},\hat{s}\rg\cong\ZZ_{\frac{q+1}{d}}\times\ZZ_{q+1}\rtimes \SS_3.$
Further, one may get $I\rtimes \lg\hat{y}\rg\cong \A_4$ and $I\rtimes \lg\hat{y},\hat{s}\rg\cong \SS_4$ where $I:=\{[1,-1,-1],[-1,-1,1],[-1,1,-1],[1,1,1]\}.$
By direct checking, we have $\N_G(L)=I\times\lg\hat{x}\rg\rtimes\lg\hat{y},\hat{s}\rg\cong\ZZ_d\rtimes \SS_4.$ This implies that {\small$\Fix(L)=\frac{|\N_G(L)|}{|\N_M(L)|}=d.$}
Similary, we have {\small$\Fix(\SS_4)=1.$}
Since $\SO(3,q)\cong\PGL(2,q)$ contains exactly $\frac{q(q^2-1)}{60}$ subgroups isomorphic to $\A_5$ when $q\equiv \pm1\pmod{10}$. Moreover, all such subgroups are conjugate in $M$. Now, we have  $\N_M(\A_5)=\A_5$. Since a direct inspection of Proposition~\ref{maximal} reveals that $\N_G(\A_5)=\A_5$, it follows that
{\small$\Fix(\A_5)=\frac{|\N_G(\A_5)|}{|\N_M(\A_5)|}=1.$}
\qed

\begin{theorem} \label{suborbit} All nonregular suborbits of $T$ on $[T: M_0]$ are listed in Table 1 and
\begin{equation} \label{l0} |\Gamma_{r}(\a)|\ge (q^5-q^4)/3d. \end{equation}
\end{theorem}
\demo Set $q\equiv \epsilon(\mod 4)$. Let   $\cal K:=\{K_\i\di 1\le \i\le 8\}$  be
  the respective representatives of conjugacy classes of  subgroups $D_{2(q+1)}$, $D_{2(q-1)}$, $\D_4^A$, $\D_4^B$, $\ZZ_p^m\rtimes\ZZ_2$, $\ZZ_2^A$, $\ZZ_2^B$  and $\A_4$ of $M$, where $\ZZ_2^A, \D_4^A\leq \PSL(2,q)$. And $K_\i$ ($\i=1, 2, 5, 8$) are maximal  in $\cal K$;
  $K_3\le K_\i\cong \D_{2(q-\e)}, K_8$;
  $K_4\le K_1, K_2$;
  $K_6\leq K_\i, \i=1,2,3,4,8$; $K_7\le K_\i,\i=1,2,4,5$. By Lemma~\ref{fix}, every point stabilizer $M_\b$  of $M$ is conjugate to $K_\i\in \cal K$, where $\b\ne \a$.
 Let {\small$k_\i:=\Fix(K_\i)$} and
$x_\i$ the number of the suborbits of $G$ with the point stabilizer $K_\i$.

For every maximal one   $K_\i$ where $\i=1, 2, 5, 8$, $k_\i=|\N_G(K_\i)|/|\N_N(K_\i)|$ so that $x_\i=(k_i-1)/|\N_{M}(K_\i):K_\i|$, see Table 1.

 Since  $K_3$ and $K_4$ are conjugate in $G$, we get  $k_3=k_4=|\N_G(K_3):\N_M(K_3)|+|\N_G(K_4):\N_M(K_4)|$. Since $K_3\le K_\i\cong \D_{2(q-\e)}, K_8$, it  fixes $|\N_{M}(K_3):\N_{K_\i}(K_3)|=3$ (resp. $|\N_{M}(K_3):\N_{K_8}(K_3)|=2$ and $|\N_M(K_3):K_3|=3$) points in every $M$-orbit with point stabilizer $K_\i$ (resp. $K_8$ and $K_3$). Hence $x_3=(k_3-1-3x_\i -2x_{8})/3$. With the same reason, we may get $x_4=[k_4-1-x_1(1+2)]/2$, as $x_1=x_2$ and see Table 1.

 Since $K_6$ and $K_7$ are conjugate in $G$, we get $k_6=k_7=|\N_G(K_6):\N_M(K_6)|+|\N_G(K_7):\N_M(K_7)|$.  Since, in each subgroup $\D_4^B$, there are two subgroups $\ZZ_2$ which are conjugate to $\ZZ_2^A$(resp. $\ZZ_2^B$) in $M$ for $q\equiv 1\pmod 4$ (resp. $q\equiv -1\pmod 4$), we have {\small$\Fix_{K_4}(K_6)=2|\N_M(K_7):\N_{K_4}(K_7)|=q-1$}
and
{\small$\Fix_{K_4}(K_7)=|\N_M(K_7):\N_{K_4}(K_7)|=\frac{q+1}{2}$}
for $q\equiv 1\pmod4$;
and {\small$\Fix_{K_4}(K_6)=|\N_M(K_7):\N_{K_4}(K_7)|=\frac{q+1}{2}$} and $\Fix_{K_4}(K_7)=2|\N_M(K_7):\N_{K_4}(K_7)|=q-1$
for $q\equiv -1\pmod4$, where {\small$\Fix_{K_4}(K_\j)$} denote the number of fixed points of $K_\j$ in each suborbit with stabilizer $K_4$ with $\j=6,7$.
With the same reason, we get
$${\small\begin{array}{lll} x_6&=&[k_6-1-(\frac{q-\epsilon}{2}+\frac{q-\epsilon}{2}+1)x_2-\frac{3(q-\epsilon)}{2}x_3-\Fix_{K_4}(K_6)x_4-\frac{q-\epsilon}{2}x_8]/(q-\epsilon);\\
 x_7&=&[k_7-1-(\frac{q+\epsilon}{2}+\frac{q+\epsilon}{2}+1)x_2-\Fix_{K_4}(K_7)x_4-(q+\epsilon)x_5]/(q+\epsilon),\\
 \end{array}}$$
 as $x_1=x_2.$
 Finally, using the values $x_\i$ where $1\le \i \le 8$ in Table 1, one  gets
$|\G_r(\a)|=|\O|-1-\sum_{\i=1}^{8}x_\i |M:K_\i|=\frac{1}{3d}q^5-\frac{1}{3d}q^4+\frac{2}{3}q^3+\frac{1}{3d}q^2-\frac{10-d}{9}q\ge \frac{1}{3d}q^5-\frac{1}{3d}q^4.$
\qed

\begin{table}[htp]
\caption{Suborbits of $(T, M_0)$}
\label{number}
{\footnotesize
\begin{center}
\begin{tabular}{|c |c|c|c|c|c|}
 \hline
 $d=1, 4|q-1$  \\
\hline
 i&$K_i$ &$\N_{M_0}(K_i)$ &  $\N_T(K_i)$& $k_i$& $x_i$  \\
\hline
  $1$& $\D_{2(q+1)}$ & $\D_{2(q+1)}$ & $\ZZ_{\frac{q+1}{d}}\times \ZZ_{q+1}\rtimes\ZZ_2$ & $\frac{q+1}{d}$& $\frac{q+1}{d}$-1  \\
 \hline
$2$& $\D_{2(q-1)}$& $\D_{2(q-1)}$ & $\ZZ_{\frac{q^2-1}{d}}\rtimes\ZZ_2$ & $\frac{q+1}{d}$& $\frac{q+1}{d}$-1  \\
 \hline

$3$&$\D_4^A$& $\SS_4$ & $\ZZ_{\frac{q+1}{d}}\times\ZZ_{q+1}\rtimes\SS_3$ & $\frac{(q+1)^2}{d}$ & $\frac{q^2-q}{6}$  \\
  \hline

$4$&$\D_4^B$    & $\D_8$ &$\ZZ_{\frac{q+1}{d}}\times\ZZ_{q+1}\rtimes\SS_3$ &$\frac{(q+1)^2}{d}$ &$\frac{q^2-q}{2}$ \\
  \hline
$5$&$\ZZ_p^m\rtimes\ZZ_2$& $[q]\rtimes \ZZ_{q-1}$ & $[q^2]\rtimes \ZZ_{q-1}$ & $q$ & $2$  \\
  \hline

$6$&$\ZZ_2^A$    & $\D_{2(q-1)}$ &$\SL(2,q)\rtimes\ZZ_{\frac{q+1}{d}}$ &$\frac{q^2(q+1)}{d}$ &$\frac{q^3+6q^2-3q-4}{4(q-1)}$ \\
  \hline
$7$&$\ZZ_2^B$   & $\D_{2(q+1)}$ &$\SL(2,q)\rtimes\ZZ_{\frac{q+1}{d}}$&$\frac{q^2(q+1)}{d}$&$\frac{3q^3-15q-12}{4(q+1)}$ \\
  \hline
$8$&$\A_4$    & $\SS_4$ &$\SS_4$ &$1$ &$0$ \\
  \hline

  $d=3, 4\di q-1$  \\
\hline
  $1$& $\D_{2(q+1)}$ & $\D_{2(q+1)}$ & $\ZZ_{\frac{q+1}{d}}\times \ZZ_{q+1}\rtimes\ZZ_2$ & $\frac{q+1}{d}$& $\frac{q+1}{d}$-1  \\
 \hline
$2$& $\D_{2(q-1)}$& $\D_{2(q-1)}$ & $\ZZ_{\frac{q^2-1}{d}}\rtimes\ZZ_2$ & $\frac{q+1}{d}$& $\frac{q+1}{d}$-1  \\
 \hline

$3$&$\D_4^A$& $\SS_4$ & $\ZZ_{\frac{q+1}{d}}\times\ZZ_{q+1}\rtimes\SS_3$ & $\frac{(q+1)^2}{d}$ & $\frac{(q+1)(q-2)}{18}$  \\
  \hline

$4$&$\D_4^B$    & $\D_8$ &$\ZZ_{\frac{q+1}{d}}\times\ZZ_{q+1}\rtimes\SS_3$ &$\frac{(q+1)^2}{d}$ &$\frac{q^2-q+4}{6}$ \\
  \hline
$5$&$\ZZ_p^m\rtimes\ZZ_2$& $[q]\rtimes \ZZ_{q-1}$ & $[q^2]\rtimes \ZZ_{q-1}$ & $q$ & $2$  \\
  \hline

$6$&$\ZZ_2^A$    & $\D_{2(q-1)}$ &$\SL(2,q)\rtimes\ZZ_{\frac{q+1}{d}}$ &$\frac{q^2(q+1)}{d}$ &$\frac{q^3+6q^2-7q}{12(q-1)}$ \\
  \hline
$7$&$\ZZ_2^B$   & $\D_{2(q+1)}$ &$\SL(2,q)\rtimes\ZZ_{\frac{q+1}{d}}$&$\frac{q^2(q+1)}{d}$&$\frac{q^3-9q-8}{4(q+1)}$ \\
  \hline
$8$&$\A_4$    & $\SS_4$ &$\ZZ_3\rtimes\SS_4$ &$3$ &$1$ \\
  \hline
  $d=1, 4\di q+1$  \\
\hline $1$& $\D_{2(q+1)}$ & $\D_{2(q+1)}$ & $\ZZ_{\frac{q+1}{d}}\times \ZZ_{q+1}\rtimes\ZZ_2$ & $\frac{q+1}{d}$& $\frac{q+1}{d}$-1  \\
 \hline
$2$& $\D_{2(q-1)}$& $\D_{2(q-1)}$ & $\ZZ_{\frac{q^2-1}{d}}\rtimes\ZZ_2$ & $\frac{q+1}{d}$& $\frac{q+1}{d}$-1  \\
 \hline

$3$&$\D_4^A$& $\SS_4$ & $\ZZ_{\frac{q+1}{d}}\times\ZZ_{q+1}\rtimes\SS_3$ & $\frac{(q+1)^2}{d}$ & $\frac{q^2-q}{6}$  \\
  \hline

$4$&$\D_4^B$    & $\D_8$ &$\ZZ_{\frac{q+1}{d}}\times\ZZ_{q+1}\rtimes\SS_3$ &$\frac{(q+1)^2}{d}$ &$\frac{q^2-q}{2}$ \\
  \hline
$5$&$\ZZ_p^m\rtimes\ZZ_2$& $[q]\rtimes \ZZ_{q-1}$ & $[q^2]\rtimes \ZZ_{q-1}$ & $q$ & $2$  \\
  \hline

$6$&$\ZZ_2^A$    & $\D_{2(q+1)}$ &$\SL(2,q)\rtimes\ZZ_{\frac{q+1}{d}}$ &$\frac{q^2(q+1)}{d}$ &$\frac{q^3-3q-2}{2(q+1)}$ \\
  \hline
$7$&$\ZZ_2^B$   & $\D_{2(q-1)}$ &$\SL(2,q)\rtimes\ZZ_{\frac{q+1}{d}}$&$\frac{q^2(q+1)}{d}$&$\frac{q^3+2q^2-5q+2}{2(q-1)}$ \\
  \hline
$8$&$\A_4$    & $\SS_4$ &$\SS_4$ &$1$ &$0$ \\
  \hline
 $d=3, 4\di q+1$    \\
\hline
  $1$& $\D_{2(q+1)}$ & $\D_{2(q+1)}$ & $\ZZ_{\frac{q+1}{d}}\times \ZZ_{q+1}:\ZZ_2$ & $\frac{q+1}{d}$& $\frac{q+1}{d}$-1  \\
 \hline
$2$& $\D_{2(q-1)}$& $\D_{2(q-1)}$ & $\ZZ_{\frac{q^2-1}{d}}\rtimes\ZZ_2$ & $\frac{q+1}{d}$& $\frac{q+1}{d}$-1  \\
 \hline

$3$&$\D_4^A$& $\SS_4$ & $\ZZ_{\frac{q+1}{d}}\times\ZZ_{q+1}\rtimes\SS_3$ & $\frac{(q+1)^2}{d}$ & $\frac{(q+1)(q-2)}{18}$  \\
  \hline

$4$&$\D_4^B$    & $\D_8$ &$\ZZ_{\frac{q+1}{d}}\times\ZZ_{q+1}\rtimes\SS_3$ &$\frac{(q+1)^2}{d}$ &$\frac{q^2-q+4}{6}$ \\
  \hline
$5$&$\ZZ_p^m\rtimes\ZZ_2$& $[q]\rtimes \ZZ_{q-1}$ & $[q^2]\rtimes \ZZ_{q-1}$ & $q$ & $2$  \\
  \hline

$6$&$\ZZ_2^A$    & $\D_{2(q+1)}$ &$\SL(2,q)\rtimes\ZZ_{\frac{q+1}{d}}$ &$\frac{q^2(q+1)}{d}$ &$\frac{q^3-3q-2}{6(q+1)}$ \\
  \hline
$7$&$\ZZ_2^B$   & $\D_{2(q-1)}$ &$\SL(2,q)\rtimes\ZZ_{\frac{q+1}{d}}$&$\frac{q^2(q+1)}{d}$&$\frac{q^3+2q^2-13q+10}{6(q-1)}$ \\
  \hline
$8$&$\A_4$    & $\SS_4$ &$\ZZ_3\rtimes\SS_4$ &$3$ &$1$ \\
  \hline
\end{tabular}
\par\vspace{0.5em} 
\footnotesize
Notation: $[n]$ denotes an unspecified soluble group of order $n$.
\end{center}}
\end{table}

 \section{A lower bound of $|\Gamma_{nr}(\a)\cap \Gamma_{nr}(\a')|$}
 Theorem~\ref{main21} is the main result of this section.
\begin{theorem}\label{main21}
Using the notation from Section~\ref{ell},
consider the right multiplication action of $T$ on the cosets $[T:M_0]$, where $d=(3, q+1)$ and $q\geq 7$ is odd. Suppose that $\sqrt{q}\geq 15d$. Then
$|\Gamma_{nr}(\a)\cap \Gamma_{nr}(\a')|\ge \frac{4}{9}q^5-\frac{1}{3}q^4\sqrt{q}-22.1555q^4$ for $d=1$; and
$|\Gamma_{nr}(\a)\cap \Gamma_{nr}(\a')|\ge \frac{4}{27}q^5-\frac{31}{27}q^4\sqrt{q}-16.8372q^4$ for $d=3$.
\end{theorem}
\demo From now on, set
$G=T=\PSU(3,q)$ and $M=M_0=\PSO(3,q)$.  Since $M$ is maximal in $G$, we identify each coset $Mg$ with its point stabilizer $M^g$ so that $\O=\{ M^g\mid g\in G\}$, with the conjugacy action. Set $\a=M$.
  Let $\a':=M'=M^{g'}$ for any $g'\in G$. Then we shall get a lower bound of  $|\Gamma_{nr}(\a)\cap \Gamma_{nr}(\a')|$.
 \vskip 3mm
The possible isomorphic  types for $M\cap M'$ are:
 \vspace{-5pt}\begin{equation}\label{A_i} A_q=\ZZ_p^m.\ZZ_2,\,  A_{q\pm 1}=\D_{2(q\pm 1)}, \,  A_3=\D_4,  A_4=\Alt(4)\, ({\rm if}\,  3\mid q+1), A_2=\ZZ_2 \, {\rm or}\, A_1=1.\end{equation}
All of these groups $A_\imath$ contain at least one involution, except for $A_1$. Therefore, set
\vspace{-5pt}\begin{equation}\label{I} I=\{ g\in M\mid |g|=2\}\quad {\rm and}\quad I'=\{ g\in M'\mid |g|=2\}.\end{equation}
Clearly,  $0 \le |I \cap I'| \le q+2$, giving
     \vspace{-5pt}\begin{equation}\label{lw}\Gamma_{nr}(\a)\cap \Gamma_{nr}(\a')=\{M_\imath\in \O  \setminus \{M, M'\}\mid t, t'\in  M_\imath,\,\, {\rm for \, some}\,\,  t\in I, t'\in I' \}.\end{equation}
 For the reasen of a simplicity, we turn to consider a subset ${\cal W}$ of $\Gamma_{nr}(\a)\cap \Gamma_{nr}(\a')$, where
  \vspace{-5pt}\begin{equation}\label{W_1} \cal W=\{ M_\imath\in \O \setminus \{M, M'\}\mid t, t'\in  M_\imath\,\, {\rm for \, some}\, \, t\in I, t'\in I' \,\, {\rm but}\,\,  t, t'\notin I\cap I'\}.\end{equation}
To determine $|{\cal W}|$, we introduce a set of incident triples (under the inclusion relation)
 \vspace{-5pt}\begin{equation}\label{E} \cal E=\{(t, M_\imath, t')\mid  M_\imath\in {\cal W}, t\in M\cap M_\imath, t'\in M'\cap M_\imath, t, t'\notin I\cap I'\},\end{equation}
\f while for any given $ t\in I\setminus I'$ and $M_\imath\in {\cal W}$, we introduce two respective sets of pairs
\vspace{-5pt}\begin{equation}\label{E(t)} \cal E(t)=\{(M_\imath,t')\mid (t, M_\imath, t')\in \cal E\}\,\, {\rm and}\,\,  \cal E(M_\imath)=\{(t,t')\mid (t, M_\imath, t')\in \cal E\},\end{equation}
so that  we may count  $|\cal E|$ by two different ways, that is,
\vspace{-5pt}\begin{equation} \label{2way} \sum_{t\in I\setminus I'}|\cal E(t)|=|\cal E|=\sum_{M_\imath\in {\cal W}}|\cal E(M_\imath)|.\end{equation}
Further, according to $M_\imath\cap M\cong A_j$, $M_\imath\cap M'\cong A_k$,  $M_\imath\cap M\cap M'\cong A_l$, where $j, k,l\in \{1,2,3,4,q,q\pm1\}$, as mentioned in Eq(\ref{A_i}),
 the  set ${\cal W}$ is the union of the following subsets ${\cal W}_{j, k, l}$, where
\vspace{-5pt}\begin{equation}\label{jkl}\begin{array}{lll}
    {\cal W}_{j,k, l}&=&\{M_\imath\in  {\cal W}\mid M_\imath\cap M\cong A_j, M_\imath\cap M'\cong A_k, M_\imath\cap M\cap M'=A_l\}\end{array}\end{equation}
and $j,k,l\in\{1,2,3,4,q,q\pm1\}$,
reminding that $(M_\imath\cap M) \cap (M_\imath\cap M')=M\cap M_\imath\cap M'$. Moreover, set
\begin{equation}\label{n}
 n_{j, k, l}=|{\cal W}_{j,k,l}|.\end{equation}
 Then by Eq(\ref{2way}) and Eq(\ref{jkl}),  we have
 \vspace{-5pt}$$\begin{array}{ll}
 &\sum_{t\in I\setminus I'}|\cal E(t)|=\sum_{j,k, l} n_{j,k, l} |\{\cal E(M_\imath): M_\imath\in {\cal W}_{j,k, l}\}|,
 \end{array}$$
\begin{equation}\label{w1}
\begin{array}{lll}|{\cal W}|&=&\sum\limits_{j,k, l}n_{j,k, l}=\sum\limits_{t\in I\setminus I'}|\cal E(t)|
-(\sum\limits_{t\in I\setminus I'}|\cal E(t)|-\sum\limits_{j,k, l}n_{j,k,l})\\
&=&\sum\limits_{t\in I\setminus I'}|\cal E(t)|
-\sum\limits_{j,k, l} n_{j,k,l}(|\{\cal E(M_\imath): M_\imath\in {\cal W}_{j,k, l}\}|-1)\\
&=&{\hm A}-{\hm B}-{\hm C}-{\hm D},\quad {\rm where} \end{array}
\end{equation}
$$\begin{array}{ll}
&{\hm A}=\sum\limits_{t\in I\setminus I'}|\cal E(t)|,\hskip 2cm
{\hm B}=\sum\limits_{4, q\not\in \{j, k\}}n_{j, k, 1}(|\{\cal E(M_\imath): M_\imath\in{\cal W}_{j,k,1}\}|-1),\\
&{\hm C}=\sum\limits_{4, q\not\in \{j, k\},l\ne 1}n_{j, k, l}(|\{\cal E(M_\imath): M_\imath\in{\cal W}_{j,k,l}\}|-1),\\
&{\hm D}=\sum\limits_{\{4, q\}\cap  \{j, k\}\ne \emptyset} n_{j, k, l}(|\{\cal E(M_\imath): M_\imath\in{\cal W}_{j,k, l}\}|-1).\\
\end{array}$$
\f In the following Lemmas~\ref{AB},~\ref{C} and~\ref{D} and~\ref{DD},   it will be  respectively shown that
$$\begin{array}{lll}
{\hm A}&\geq &q^5-1.0355 q^4,\,    \hbox{$d=1;$}\,\,  \frac 13q^5-0.3342 q^4, \,  \hbox{$d=3;$}\,\\
{\hm B} &\leq & \frac{5}{9}q^5+\frac{1}{3}q^4\sqrt{q}+15.25 q^4,  \, \hbox{$d=1$;} \,
                  \frac{5}{27}q^5+\frac{31}{27}q^4\sqrt{q}+10.74 q^4,  \, \hbox{$d=3$;}\\
 {\hm C}&\leq&   3.85q^4, \hbox{$d=1$;} \,  3.26 q^4,  \hbox{$d=3$};\\
               {\hm D}&\leq&   2.02 q^4, \, \hbox{$d=1$;}  \,  2.503 q^4,\,  \hbox{$d=3$}.  \\
\end{array}$$
\f By using    Eq(\ref{lw}), Eq(\ref{W_1}) and Eq(\ref{w1}),  we get $|\Gamma_{nr}(\a)\cap \Gamma_{nr}(\a')|\ge|{\cal W}|= {\hm A}-{\hm B}-{\hm C}-{\hm D}$.  Inserting these bounds  of ${\hm A}$, ${\hm B}$, ${\hm C}$, ${\hm D}$, we get  a lower bound for $|\Gamma_{nr}(\a)\cap \Gamma_{nr}(\a')|$   as shown in Theorem~\ref{main21}, omitting the details.
\qed

\subsection{Measure ${\hm A}$}
Note that every dihedral subgroup of $\PSO(3,q)$ is isomorphic to $D_4$, $D_{2p}$ or $D_{2l}$ with $l\mid q\pm1$ and $l>2$. To measure ${\hm A}$,  for any $t\in I$ and $s=2, p$ or $s\mid (q\pm 1)$ but $s\ne 2$,   set
  \vspace{-5pt}$$ I'(t,s)=|\{t'\in I'\setminus{I}\mid o(tt')=s\}|.$$
Denote by $\Fix_{\mathcal{H}}(t)$ and $\Fix_{\mathcal{H}}(D)$ the number of isotropic points fixed by $t$ and $D\leq \PSU(3,q)$, respectively. To measure $I'(t,p)$, we need the following lemma.
\begin{lem}\label{used0}
Let $t\in G=\PSU(3,q)$ be an involution and $D\cong \D_{2p}$ a subgroup of $G$.
 Then we have
(1) $\Fix_{\mathcal{H}}(t)=q+1$; (2) $\Fix_{\mathcal{H}}(D)=1$ and for any $\lg y_1\rg, \lg y_2\rg\in \mathcal{N}$,  $\lg \tau_{y_1},\tau_{y_2}\rg\cong \D_{2p}$ if and only if $\lg y_1\rg, \lg y_2\rg\in \PG(v^\perp)$ for some $\lg v\rg\in{\cal H}$.\end{lem}
\demo (1)  Let $t$ be an involution. By Lemma \ref{lem-inv-bij}, $t=\tau_{y}$ for a point $\lg y\rg\in\mathcal{N}$. Then for any $\lg x\rg\neq\lg y\rg$,  $\tau_y(\lg x\rg)=\lg x\rg$ if and only if  $\lg x\rg\in \PG(y^\perp)$.
Since each line contains $q+1$ isotropic points in $\PG(2,q^2)$, we get
 $\Fix_{\mathcal{H}}(t)=\Fix_{\mathcal{H}}(\tau_y)=|\PG(y^\perp)\cap {\cal H}|=q+1$.
\vskip 3mm
(2) From the proof of Lemma~\ref{fix}.(2),  we  know  that $D\leq G_{\lg v\rg}$ for some $\lg v\rg\in{\cal H}$, so that   $\Fix_{{\cal H}}(D)=\frac{|\N_G(D)|}{|\N_{G_{\lg v\rg}}(D)|}=\frac{|Q_1\rtimes Y|}{|Q_1\rtimes Y|}=1$, meaning   $\Fix_{\mathcal{H}}(D)=1$, where $Q_1$ and $Y$ were defined in the proof of Lemma~\ref{fix}.(2).

 Suppose that $D:=\lg \tau_{y_1},\tau_{y_2}\rg\cong \D_{2p}$.
 The above arguments follow that  $\lg y_1\rg,\lg y_2\rg\in\PG(v^\perp)\cap {\cal N}$, where $\lg v\rg$ is the unique isotropic point fixed by $D$. Conversely, suppose that  $\lg y_1\rg, \lg y_2\rg\in \PG(v^\perp)\cap {\cal N}$ for some $\lg v\rg\in{\cal H}$.  Then $\PG(y_1^\perp)\cap \PG(y_2^\perp)=\lg v\rg$. Note that every   dihedral subgroup of $G$ is isomorphic to a subgroup of $\D_{2p}$, $\D_{2(q\pm 1)}$ or $\D_{\frac{2(q^2-1)}{d}}$.
 If $D\lesssim\D_{2(q\pm 1)}$ or $\D_{\frac{2(q^2-1)}{d}}$, then $D$ contains a unique central involution, say $\lg z\rg\in{\cal N}$ such that $\lg y_1\rg,\lg y_2\rg\in\PG(z^\perp)$, contradicting with $\lg y_1\rg,\lg y_2\rg\in\PG(v^\perp)$ for $\lg v\rg \in {\cal H}$.
 Therefore, $D\cong D_{2p}$.
\qed

\begin{lem}\label{used1}
 Let $t\in I$. Then $I'(t,2)=0, 1$   if $t\in I\setminus I';$ and $I'(t,p)\leq 3q-1$.
\end{lem}
   \demo
The first conclusion follows from Proposition \ref{lem_obser} and Lemma \ref{lem-inv-bij} directly.
To prove the second one,  let $\mathbf{w},\mathbf{w}'\in\mathcal{B}_0$ with $M=G_{\mathbf{w}}$ and $M'=G_{\mathbf{w}'}$. Set $t=\tau_y$ with $\lg y\rg\in\mathbf{w}\cap {\cal N}$. Then by Lemma~\ref{used0}, we get
$I'(t,p)=|\{t'\in I'\setminus{I}\mid o(tt')=p\}|=|\Upsilon|,$ where  \vspace{-5pt}$$
\Upsilon=\{ \lg y'\rg\in({\cal N}\cap \mathbf{w'})\setminus\mathbf{w}\mid \{\lg y\rg,\lg y'\rg\}\subseteq \PG(v^\perp),~\text{for some}~\lg v\rg\in{\cal H}\}.$$
\f Let $\cal H_1=\PG(y^\perp)\cap \cal H$. Then $|\cal H_1|=q+1$ and  $\mathbf{w'}$ contains $k$ points of  $\cal H_1$, where $k\le 2.$
(1) First suppose $k=0$. Then for any $\lg v'\rg\in \cal H_1$,
by Proposition~\ref{lem_obser}.(4),  $|\PG(v'^\perp)\cap \mathbf{w'}|=1$ which implies there is at most one $\lg y'\rg$ contained in $\Upsilon$. Therefore,  $|\Upsilon|\le q+1$ in this case;
(2) Secondly, suppose  $k=1$, saying  $\cal H_1\cap \mathbf{w'}=\{ \lg v_1\rg \}.$  Then $\mathbf{w'}$ contains $q$ nonisotropic points of $\PG(v_1^\perp)$ so that we get  at most $q$ points $\lg y'\rg \in \Upsilon$ related to $\lg v_1^\perp\rg $. Together with other (at most $q$) points $\lg y'\rg $ related  to    points in $\HH_1\setminus \{ \lg v_1\rg \}$, we get $|\Upsilon|\le 2q$; (3) Finally, assume $k=2$. Similarly we get  $|\Upsilon|\le q+q+(q-1)=3q-1,$ as desired.
\qed
\begin{lem}\label{AB}
$ {\hm A}\geq q^5-1.0355 q^4$ for $d=1$, and $\frac 13q^5-0.3342 q^4$ for $d=3$.
\end{lem}
\demo From Table \ref{number}, we get that every subgroup $D_4$ (resp. $D_{2p}$, $D_2({q+1)}$ and $D_2({q-1)}$) of $G$ is contained in $\frac{(q+1)^2}d$ (resp. $q$, $\frac{q+1}d$ and  $\frac{q+1}d$) subgroups $M_\i\cong \SO(3,q)$.
By  Lemma \ref{used1} and  facts  $|I'|=q^2$ and $|I \cap I'| \le q+2$, we get that for any $t\in I\setminus I'$,
$$
\begin{array}{lll}
|\cal E(t)|&=&|\{(M_\imath,t')\mid (t, M_\imath, t')\in \cal E\}|\\
&=&\frac{(q+1)^2}{d}I'(t,2)+qI'(t,p)+\frac{q+1}{d}I'(t,q+1)+\frac{q+1}{d}I'(t,q-1)\\
&=&\frac{(q+1)^2}{d}I'(t,2)+qI'(t,p)+\frac{q+1}{d}[q^2-|I\cap I'|-I'(t,2)-I'(t,p)]\\
&=&\frac{q^2(q+1)}{d}+\frac{q^2+q}{d}I'(t,2)+\frac{(d-1)q-1}{d}I'(t,p)-\frac{q+1}{d}|I\cap I'|\\
&\ge &\frac{q^2(q+1)}{d}-\frac{q+1}{d}(q+2)+\frac{(d-1)q-1}{d}I'(t,p)
\ge   q^3-6q-1,  \hbox{$d=1;$} \,
         \frac{q^3-3q-2}{3},\, \hbox{$d=3.$}
     \end{array}$$
Since $|I\setminus I'|\geq q^2-q-2$, we have
${\hm A}=\sum_{t\in I\setminus I'}|{\cal E}(t)|\geq (q^2-q-2)\min_{t\in I\setminus I'}|\cal E(t)|$. Inserting the above lower bound for $|\cal E(t)|,$ we obtain the bound stated in the lemma.
\qed

\subsection{Measure ${\hm B}$}
To obtain an upper bound for ${\hm B}$, we first need the following two lemmas.
\begin{lem}\label{Dqpm1}
Let $M_1$ and $M_2$ be any two distinct subgroups in $\Omega:=\{M^g\mid g\in G\}$ with corresponding Baer subplanes $\mathbf{w}_1$ and $\mathbf{w}_2$ (i.e., $G_{\mathbf{w}_\imath} = M_\imath$ for $\imath=1,2$).
Suppose  $M_1 \cap M_2 \cong \D_{2(q+1)}$ (resp. $\D_{2(q-1)}$), then there is a point $\lg y_0\rg\in{\cal N}$ such that $\PG(y_0^\perp)\cap \mathbf{w}_1=\PG(y_0^\perp)\cap \mathbf{w}_2\subseteq \mathbf{w}_1\cap \mathbf{w}_2$, which is a Baer subline containing $0$ (resp. $2$)  isotropic points.
 \end{lem}
\demo Suppose $M_1\cap M_2\cong \D_{2(q\pm1)}$. Let $\mathcal{I}_0$ be the set of involutions of $M_1\cap M_2$. Write  $\mathcal{I}_0=\{  \tau_{y_0},\tau_{y_\imath}\mid \lg y_0\rg, \lg y_\i\rg\in\mathcal{N}, [\tau_{y_0},\tau_{y_\i}]=1,  1\leq \i\leq q\pm1\}$.
  Then by Lemma \ref{lem-inv-bij}, $\lg y_\i\rg\in \PG(y_0^\perp)$ and $\{\lg y_0\rg,\lg y_\i\rg\mid 1\leq \i\leq q+1\}\subseteq \mathbf{w}_1\cap \mathbf{w}_2$.
  Following the fact that $\PG(y_0^\perp)\cap\mathbf{w}_1\cap\mathbf{w}_2$ is a Baer subline, we know that it  containing $0$ (resp. $2$)  isotropic points provided  $M_1\cap M_2\cong \D_{2(q+1)}$ (resp. $\D_{2(q-1)}$).
\qed

\begin{lem}\label{q+1}
Let $n_{j,k,l}$ be given in {\rm Eq(\ref{n})}. Then
$n_{q+1,q-1,1}=n_{q-1,q+1,1}=n_{q+1,q+1,1}=0.$
\end{lem}
\demo
We only prove $n_{q+1,q-1,1}=0$ in detail and have similar arguments for others.
On the contrary, we assume that there exists $M_\imath\in \mathcal{W}$ such that $M_\imath\cap M\cong \D_{2(q+1)}$, $M_\imath\cap M'\cong \D_{2(q-1)}$ and $M\cap M_\i\cap M'=1$. Set $M:=G_{\mathbf{w}}$, $M_\i:=G_{\mathbf{w}_\i}$ and $M':=G_{\mathbf{w}'}$, where $\mathbf{w},\mathbf{w}_\i,\mathbf{w}'\in\mathcal{B}_0$.

Note that each Baer subline has  $q+1$ points and at most two of them are isotropic.
In viewing of  the proof of Lemma \ref{Dqpm1}, since $M_\imath\cap M\cong \D_{2(q+1)}$ (resp. $M_\imath\cap M'\cong \D_{2(q-1)}$), the Baer subline $\ell_1\in \mathbf{w}_\i\cap \mathbf{w}$ (resp. $\ell_2\in \mathbf{w}_\i\cap \mathbf{w'}$) contains no isotropic points (resp. $2$ isotropic points).
However, since $M\cap M_\imath\cap M'=1$ (containing no involutions),  $\ell_1\cap \ell_2$ ($\subseteq \mathbf{w}_\i$)  must be  isotropic, getting a contradiction.  Therefore, $n_{q+1,q-1,1}=0$.
\qed

\begin{lem}\label{C} ${\hm B}\leq  \frac{5}{9}q^5+\frac{1}{3}q^4\sqrt{q}+15.25 q^4,  \, \hbox{$d=1$;} \,
                  \frac{5}{27}q^5+\frac{31}{27}q^4\sqrt{q}+10.74 q^4,  \, \hbox{$d=3$.}$
\end{lem}
\demo Recall ${\hm B}=\sum\limits_{4, q\not\in \{j, k\}}n_{j,k,1}(|\{\cal E(M_\imath): M_\imath\in{\cal W}_{j,k,1}\}|-1)$.
By Lemma~\ref{q+1}, we get  $n_{q+1,q+1,1}=n_{q+1,q-1,1}=n_{q-1,q+1,1}=0.$
 Inserting the values of $|\{\cal E(M_\imath): M_\imath\in{\cal W}_{j,k,1}\}|-1$ (see Eq(\ref{E(t)})) for each $M_\imath\in {\cal W}_{j,k, 1}$, we get
$$\begin{array}{ll}
{\hm B}&=2(n_{2,3,1}+n_{3,2,1})+(q+1)(n_{2,q+1,1}+n_{q+1,2,1})+(q-1)(n_{2,q-1,1}+n_{q-1,2,1})+8n_{3,3,1}\\
&+[3(q+2)-1](n_{3,q+1,1}+n_{q+1,3,1})+(3q-1)(n_{3,q-1,1}+n_{q-1,3,1})
+(q^2-1)n_{q-1,q-1,1}.\\
\end{array}$$
For a technical reason, we write ${\hm B}$ to be a sum of ${\hm B_\i},$ where $1\le \i\le 4$ such that every number $(|\{\cal E(M_\imath): M_\imath\in {\cal W}_{j,k,1}\}|-1)n_{j, k, 1}$ is cut into four parts, while each part is one term of ${\hm B_\i}$,
$$\begin{array}{lll}
{\hm B_1}&=&2(n_{2,3,1}+\frac{q+1}{2}n_{2,q+1,1}+\frac{q-1}{2}n_{2,q-1,1}+n_{3,3,1}+\frac{q+1}{2}n_{3,q+1,1}+\frac{q-1}{2}n_{3,q-1,1}),\\
{\hm B_2}&=&2(n_{3,2,1}+\frac{q+1}{2}n_{q+1,2,1}+\frac{q-1}{2}n_{q-1,2,1}
+3n_{3,3,1}+\frac{3(q+1)}{2}n_{q+1,3,1}+\frac{3(q-1)}{2}n_{q-1,3,1}),\\
{\hm B_3}&=&(q^2-1)n_{q-1,q-1,1},\\
{\hm B_4}&=&2[(q+2)n_{3,q+1,1}+qn_{3,q-1,1}
  +n_{q+1,3,1}+n_{q-1,3,1}].\\
\end{array}$$
\f For example,  for each $M_\imath\in {\cal W}_{3,q+1,1}$, we have $3(q+2)$ triples $(t, M_\imath, t')\in {\cal E}$. Therefore,   the number $(|\{\cal E(M_\imath): M_\imath\in {\cal W}_{3,q+1,1}\}|-1)n_{3, q+1, 1}$ is equal to $(3q+5)n_{3, q+1, 1}$, while
we put $2\cdot \frac{q+1}2n_{3, q+1, 1}$ in ${\hm B_1}$, $0\cdot n_{3, q+1, 1}$ in ${\hm B_2}$, $0 \cdot n_{3, q+1, 1}$ in ${\hm B_3}$ and  $2(q+2)n_{3, q+1, 1}$ in ${\hm B_4}$.

It will be  shown  in Lemmas~\ref{B1B2} and ~\ref{C3}, that
$$ \begin{array}{l}
 {\bf B}_1 \leq \left\{
    \begin{array}{ll}
   \frac{q(q^2-1)}{3}( \frac{2}{3}q^2+q\sqrt{q}+13.34 q)+5.1q^4, & \hbox{$d=1$;} \\
    \frac{q(q^2-1)}{3}(\frac{2}{9}q^2+\frac{16}{9}q\sqrt{q}+12.4367q)+1.7q^4, & \hbox{$d=3$,}
    \end{array}
  \right.\\
   {\bf B}_2  \leq \left\{
    \begin{array}{ll}
    \frac{q(q^2-1)}{3}(q^2+8q), & \hbox{$d=1$;} \\
     \frac{q(q^2-1)}{3}(\frac{1}{3}q^2+\frac{5}{3}q\sqrt{q}+9.67q), & \hbox{$d=3$,}
    \end{array}
  \right. \\
 {\hm B_3}\leq q(q^2-1)(q+1), \, {\hm B_4}\le \frac{2}{d}q(q+1)(q+2)(q+1-d),\end{array}$$
  so that  an upper bound for ${\hm B}$ as shown in the lemma may be obtained by summing these above upper bounds for  ${\hm B_\i}$, where $1\le \i\le 4$.\qed

\subsubsection{Upper bounds of ${\bf B}_1$ and ${\bf B}_2$}\label{secC12}
Upper bounds of  ${\hm B_1}$ and ${\hm B_2}$ will be given in Lemma~\ref{B1B2}. First we need to give a characterization for them.
 For any given $\D_4\cong D'\le M'$ and  $\D_4\cong D\le M$, set
\begin{equation}\label{lh}
 \begin{array}{ll}
&{\bf R} (M, M', D')=\{ M_\imath\mid  M_\imath\in {\cal W}_{j,k, 1}, j\in \{2, 3\},k\in \{3, q+1, q-1\},   D'\le M_\imath\},\\
&{\bf L} (M, M', D)=\{ M_\imath \mid  M_\imath\in {\cal W}_{j,k, 1}, j\in \{3, q+1, q-1\}, k\in \{2, 3\}, D\le M_\imath\},\\
&{\bf H} (M, M', D)=\{(M_\imath, t)\mid  M_\imath\in {\cal W}_{j,k, 1},  j\in \{3, q+1, q-1\}, k\in \{2, 3\},\\
&  \hskip 3cm
D\le M_\imath, t\in M_\imath\cap M'\}.\\
\end{array} \end{equation}
The number of subgroups $M_\imath$ contained in the union $\bigcup_{\D_4\cong D'\le M'} {\bf R}(M, M', D')$ is determined by the count of subgroups $\D_4$ within $\D_{2(q\pm 1)}$. Specifically, since $\D_{2(q+1)}$ and $\D_{2(q-1)}$ contain $(q+1)/2$ and $(q-1)/2$ subgroups $\D_4$ respectively, the totals are $\frac{q+1}{2} n_{j, q+1, 1}$, $\frac{q-1}{2} n_{j, q-1, 1}$ and $n_{j,3,1}$ with $j=2,3$. Concurrently, the union $\bigcup_{\D_4 \cong D\le M} {\bf H}(M, M', D)$ contains $\frac{3(q+1)}{2} n_{q+1, 3, 1}$ pairs $(M_\imath, t)$ for $M_\imath\in{\cal W}_{q+1,3,1}$, with other cases following similar arguments.
 Note that $M\cong \PGL(2,q)$ contains  $n=\frac{q(q^2-1)}{6}$ subgroups isomorphic to $\D_4$. From the above arguments and the definitions of ${\hm B_1}$ and ${\hm B_2}$, we have
\vspace{-5pt}\begin{equation}\label{D=D'} \begin{array}{lll} {\hm B_1}&=&2\sum\limits_{\D_4\cong D'\le M'}|{\bf R} (M, M', D')|\, \, {\rm and}\,\,
  {\hm B_2}=2\sum\limits_{\D_4\cong D\le M} |{\bf H} (M, M', D)|.\end{array}\end{equation}

To derive upper bounds for ${\bf B}_1$ and ${\bf B}_2$, we need to  make a preparatory translation of their group-theoretic description into a geometric one.

Set $\mathbf{w}_0:=\PG(2,q)$.
We shall identify $\Omega$ with the set $\mathcal{B}_0=\{\mathbf{w}_0^g\mid g\in\PSU(3,q)\}$.
So every couple $M$ and $M'$ are identified with two Baer subplanes $\mathbf{w}$ and $\mathbf{w}'$ in $\mathcal{B}_0$, while by Lemma \ref{lem-inv-bij} any $D_4\cong D\leq M$ is identified with a orthogonal frame $\mathcal{{\hm D}}$ of $\mathbf{w}$,
noting that $M_\imath\in {\cal W}_{j,k, 1}$ gives  $ \mathcal{{\hm D}}\cap \mathbf{w}'=\emptyset.$
Under this identification,  the sets ${\cal W}$, ${\cal W}_{j,k, 1}$,   ${\bf L} (M, M', D)$ and ${\bf H} (M, M', D)$ are respectively identified with
\begin{equation}\label{'}
\begin{array}{lll}{\cal W}'&=&\{ \mathbf{w}_\i\in  \mathcal{B}_0\setminus \{\mathbf{w}, \mathbf{w}'\}\mid  \lg \g\rg , \lg \g'\rg \in \mathbf{w}_\i,  {\rm for \, some}\,  \\&&\lg \g \rg \in {\cal N}\cap (\mathbf{w}\setminus \mathbf{w}'),
\lg \g'\rg \in {\cal N}\cap (\mathbf{w}'\setminus \mathbf{w})\},\\
{\cal W}_{j,k, 1}'&=&\{\mathbf{w}_\i\in  {\cal W'}\mid G_{\mathbf{w}}\cap G_{\mathbf{w}_\i}\cong A_j,
G_{\mathbf{w}_\i}\cap G_{\mathbf{w}'}\cong A_k, G_{\mathbf{w}_\i}\cap\\ && G_{\mathbf{w}}\cap G_{\mathbf{w}'}=1\}; \\
{\bf L}'(\mathbf{w}, \mathbf{w}', \cal D)&=&\{\mathbf{w}_\i\mid \mathbf{w}_\i\in {\cal W'}_{j,k, 1},   j\in \{3, q+1, q-1\},  k\in \{2, 3\}, \\&&{\cal D}\subseteq \mathbf{w}_\i\cap \mathbf{w}\},\\
{\bf H}'(\mathbf{w}, \mathbf{w}', \cal D)&=&\{(\mathbf{w}_\i,\lg \g'\rg)\mid \mathbf{w}_\i\in {\cal W}_{j,k, 1}', j\in \{3, q+1, q-1\}, k\in \{2, 3\},\\
&&  {\cal D}\subseteq \mathbf{w}\cap \mathbf{w}_i, \lg \g'\rg \in \mathbf{w}_\i\cap \mathbf{w}'\cap {\cal N}\}, \end{array} \end{equation}
where $A_j$, $A_k$ are defined in Eq(\ref{A_i}) and $\mathcal{D}$ is a given orthogonal frame of $\mathbf{w}$.
By the transitivity of the action of $G=\PSU(3,q)$ on $\O=\{M^g\mid g\in G\}$ and the symmetry of notation for $M$ and $M'$, by Eq(\ref{D=D'}), we have
$\mathbf{B}_1\leq 2 \max\limits_{M'\in\mathcal{M}\setminus\{M\}}(\sum\limits_{\D_4\cong D\le M}|{\bf L} (M, M', D)|).$
As defined, $\mathcal{F}_{\mathbf{w}}$ is the set of all orthogonal frames of the Baer subplane $\mathbf{w}$. With the geometric description in Eq(\ref{'}),  two inequalities in Eq(\ref{D=D'}) may be redescribed as \begin{equation}\label{B12}
{\bf B}_1\leq
2\sum\limits_{\mathcal{D}\in\mathcal{F}_{\mathbf{w}}}\max\limits_{\mathbf{w}'\in\mathcal{B}_0\setminus\mathbf{w}}
\{|{\bf L'} (\mathbf{w}, \mathbf{w}', \mathcal{D})|\} \, \, {\rm and}\, \,
{\bf B}_2\le 2 \sum\limits_{\mathcal{D}\in\mathcal{F}_{\mathbf{w}}}\max\limits_{\mathbf{w}'\in\mathcal{B}_0\setminus\mathbf{w}}
\{|{\bf H'} (\mathbf{w}, \mathbf{w}', \mathcal{D})|\}.
\end{equation}

The remaining is  to measure ${\bf L}'(\mathbf{w}, \mathbf{w}', \cal D)$ and ${\bf H}'(\mathbf{w}, \mathbf{w}', \cal D)$ for any two fixed Baer subplanes $\mathbf{w},\mathbf{w}'\in\mathcal{B}_0$ and a fixed orthogonal frame $\mathcal{D}\in \mathbf{w}$ which is disjoint from $\mathbf{w}'$. To do this, we need some notation. Recall that $[x,y,z]$ denotes  a diagonal  matrix.
We take the identity matrix as our Gram matrix $J$.
Set $g:=g_{A}+g_{B}i\in\{{\hm g_1}, {\hm g_2}\},$ where
\begin{eqnarray} \label{g1g2}
{\hm g_1}={\small\left(
  \begin{array}{ccc}
    0 & 1 & m_1+n_1i \\
    e+i& t & 0\\
    t& i & m_2+n_2i \\
  \end{array}
\right)}\,\,\text{and}\,\,\,\, {\hm g_2}=
{\small\left(
  \begin{array}{ccc}
    0 & 1 & m_1+n_1i \\
   s& -1 & 0\\
    i& si & m_2+n_2i \\
  \end{array}
\right),}\end{eqnarray}
where $\FF_q^*=\lg \th \rg $, $i^2=\th $  and $e,s,t,m_1,m_2,n_1,n_2\in\FF_q$,  with $m_2n_1-m_1n_2=1, t\ne 0$ for $g={\hm g_1}$; and $m_2n_1-m_1n_2=s, s\ne 0$ for $g={\hm g_2}$.
Further, set
\begin{equation}\label{bc}
{\small\begin {array}{ll} &(b^*, c^*)\in \{ (-b+i, 1), (1,-c+i), (-b+i,-c+i)\},   a'\in \{1, a+i\}, a, b, c\in \FF_q; \\
  &[1,  b^*, c^*]=E_A+E_Bi; \\
& A_{(a', b^*, c^*)}=\left\{
               \begin{array}{ll}
                 g_AE_B+g_BE_A, & \hbox{$a'=1$;} \\
                 a(g_AE_B+g_BE_A)+(g_AE_A+g_BE_B\th), & \hbox{$a'=a+i;$}
               \end{array}
             \right.
\end{array}}\end{equation}
\begin{equation}\label{HL}
 {\small\begin{array}{lll}
 {\bf H}(\mathcal{D},g)&=&\{((b^*,c^*),\lg Y \rg)\mid YA_{(a',b^*,c^*)}=0, Y\in\FF_q^3, YY^{T}\neq 0, \\&&  \rank(A_{(a',b^*,c^*)})=2, b^*c^*\in(\FF_{q^2}^*)^3\};\\
  {\bf L}(\mathcal{D},g)&=&\{(b^*,c^*)\mid  ((b^*,c^*),\lg Y \rg)\in {\bf H}(\mathcal{D},g)\,\, \text{for some}\,\, Y\in\FF_q^3\,\, \text{with}\,\,Y {Y}^{T}\neq 0\}.\\
 \end{array}}
 \end{equation}
With the notation given in Eq(\ref{g1g2}), Eq(\ref{bc}) and Eq(\ref{HL}),  we are ready to measure ${\bf L}'(\mathbf{w}, \mathbf{w}', \cal D)$ and ${\bf H}'(\mathbf{w}, \mathbf{w}', \cal D)$.
\begin{lem}\label{hlD}
Let $\mathbf{w},\mathbf{w}'\in\mathcal{B}_0$ be any two given Baer subplanes and let $\mathcal{D}\in\mathcal{F}_{\mathbf{w}}$ be a given orthogonal frame of $\mathbf{w}$ which is disjoint from $\mathbf{w}'$.
Then
\vspace{-5pt}$$\quad |{\bf L'}(\mathbf{w}, \mathbf{w}', \cal D)|\leq|{\bf L}(\mathcal{D},g)|\quad  {\rm and}\quad  |{\bf H'}(\mathbf{w}, \mathbf{w}', \cal D)|\leq|{\bf H}(\mathcal{D},g)|.$$
 \end{lem}
\demo Let $\mathbf{w}=\PG(W)$, $\mathbf{w}_\i=\PG(W_\i)$, $\mathbf{w}'=\PG(W')$ and let $\mathcal{\hm D}=\{\lg\a\rg,\lg\b\rg,\lg\g\rg\}$ be a fixed  orthogonal frame of $\mathbf{w}$ satisfying   $\mathcal{\hm D}\cap \mathbf{w'}=\emptyset$. Set $\De=(\a, \b, \g)^T$, where $A^T$ denotes the transpose of the matrix $A$, and set $W'=\lg \a', \b', \g'\rg_{\FF_q} $, where $(\a', \b', \g')^T=g\De$ for some $g\in \SU(3,q)$.
Then $\mathbf{w}'=\{\lg Yg\Delta\rg\mid Y\in\FF_q^3\}.$

For any $(\mathbf{w}_\i,\lg \g'\rg)\in {\bf H'}(\mathbf{w}, \mathbf{w}', \cal D)$, since ${\cal D}\subseteq \mathbf{w}_\i\in {\cal W'}_{j,k, 1},$ where $j\in \{3, q+1, q-1\}$ and $k\in \{2, 3\},$ we may set
$W_\i=\lg \a, b' \b, c'\g\rg_{\FF_q} $ for some $b', c'\in \FF_{q^2}^*$.
Set $E=[1, b', c']$.  Then for any  point $\lg \d\rg \in \mathbf{w}_\i\cap \mathbf{w}'\cap{\cal N}$, we have
$XE\De =\d=a'Yg\De$ for some row vectors $X, Y\in \FF_q^3$, $YY^T\ne 0$, $a'\in \FF_{q^2}^*$,
 which is equivalent to  $XE=a'Y g$, that is
 \vspace{-5pt}$$X=Yg(a'E^{-1}).$$
 Replacing    $Y$ with $mY$ for some $m\in \FF_q^*$, we may let $a'\in \{ 1, a+i\}$ and $E^{-1}=[1, b^*, c^*]$,  where $b^*\in \{1, -b+i\}, c^*\in \{1, -c+i\}$ for some $a, b, c\in \FF_q.$
  Set $g:=g_A+g_Bi$ and $a'E^{-1}=F_A+F_Bi$, where $g_A,g_B,F_A,F_B\in\mathbf{M}_{3\times 3}(\FF_q)$. Then  we have
 \vspace{-5pt}$$X=Y(g_AF_A+\th g_BF_B)\quad {\rm and}\quad Y(g_AF_B+g_BF_A)=0.$$
 Consequently, $\d$ has a solution if and only if the following equation admits a nonzero solution $Y$ over $\FF_q$,
 \vspace{-5pt} \begin{equation}\label{Y} Y(g_AF_B+g_BF_A)=0, \quad YY^T\ne 0.\end{equation}
 We carry out the proof by the following two steps.
 \vskip 3mm
{\it Step 1: Show  $g\in \{ {\hm g_1}, {\hm g_2}\}$.}
\vskip 3mm

Since $\mathbf{w}'=\{\lg Yg\De\rg\mid Y\in\FF_q^3\}$ and for any $Q\in \GL(3,q)$, by a  linear transformation of the basis $g\Del$ with respect to $Q$, we have  $\mathbf{w}'=\{\lg YQg\De\rg\mid Y\in\FF_q^3\}$.
By Proposition \ref{lem_obser}, we have $|\PG(\a^{\perp})\cap \mathbf{w}'|=1$, which implies that $\PG(\a^\perp)\cap \{\lg Yg\De\rg\mid Y\in\FF_q^3\}=\lg (0,1,m_1+n_1i)\De\rg$, where $(m_1,n_1)\in\FF_{q}^2\setminus\{(0,0)\}$.
 Using Lemma \ref{wsigma} and  by rechoosing $Q$,   we may let $g\in\{{\hm g_1},{\hm g_2}\}$, as defined in Eq(\ref{g1g2}).
Remind
 that each row of ${\hm g_\i}\De$ $(\i=1,2)$ corresponds to a point in $\mathbf{w}'$, and these points are unnecessarily nonisotropic.
\vskip 3mm
{\it Step 2:  Show $|{\bf L'}(\mathbf{w}, \mathbf{w}', \cal D)|\leq|{\bf L}(\mathcal{D},g)|\quad  {\rm and}\quad  |{\bf H'}(\mathbf{w}, \mathbf{w}', \cal D)|\leq|{\bf H}(\mathcal{D},g)|.$}
\vskip 3mm

Let $g:=g_{A}+g_{B}i\in\{{\hm g_1}, {\hm g_2}\}$. Set $E^{-1}=[1, b^*, c^*]$ where $b^*\in \{1, -b+i\}, c^*\in \{1, -c+i\}$ for some $a, b, c\in \FF_q;$ and  $a'E^{-1}=F_A+F_Bi$;
 and $A=A_{(a',b^*,c^*)}:=g_AF_B+g_BF_A$, that is   Eq(\ref{bc}). Now, $\mathbf{w}_\i=\{\lg XE\De \rg\mid X\in \FF_q^3\}=\{\lg X[1,(b^*)^{-1},(c^*)^{-1}]\De \rg\mid X\in \FF_q^3\}$ and $\mathbf{w}'=\{\lg Yg\De \rg\mid Y\in \FF_q^3\}$. From the above arguments,
$\mathbf{w}_\i\cap \mathbf{w}'$ contains a point  $\lg \d\rg $  if and only if Eq(\ref{Y}) has a nonzero solution $Y$ over $\FF_q$,
that is $YA_{(a',b^*,c^*)}=0$ and $YY^T\neq 0$.
Followed by the computations for the suborbits in Table 1, we know that $|\mathbf{w}_\i\cap \mathbf{w}'\cap {\cal N}|\in\{0, 1, 3, q, q+2\}$, where each value in this set gives the number of involutions in $M_\i \cap M'$.
Moreover concretely, for our $g$ and $\mathcal{D}$, we remark that
\vskip 3mm
(1)  $|\mathbf{w}_\i\cap \mathbf{w}'\cap {\cal N}|=0$: in this case, for any  $a'\in \{1, a+i\di a\in \FF_q\}$,  either $\rank(A)=3$  or
 $\rank(A)=2$ and $YA=0$ has a solution $\lg Y\rg $ with $YY^T=0$;
\vskip 3mm
(2)   $|\mathbf{w}_\i\cap \mathbf{w}'\cap {\cal N}|=1  \, (resp.\, 3):$  in this case,      $|A|=0$  has only one  (resp. three distinct) solutions $a'$
such that $\rank(A)=2$ and  $YA=0$  has a solution $\lg Y\rg $ with   $YY^T\ne 0$;
\vskip 3mm
(3) If $\rank(A)=1$ for some $a'\in\{1,a+i\}$, then $|\mathbf{w}_\i\cap \mathbf{w}'\cap{\cal N}|\geq q$ and $G_{\mathbf{w}_\i}\cap G_{\mathbf{w}'}\cong \D_{2(q\pm1)}$ or $\ZZ_p^m.\ZZ_2$;
\vskip 3mm
(4) If $G_{\mathbf{w}_\i}\cap G_{\mathbf{w}'}\cong \D_{2(q\pm 1)}$, then $\rank(A)=1$ for some $a'\in\{1,a+i\}$.
\vskip 3mm
The first two items are obviously. Now, we are showing the last two ones. Followed from Table 1 and Lemmas \ref{used0} and \ref{Dqpm1}, we have that $|\mathbf{w}_\i\cap \mathbf{w}'\cap {\cal N}|\geq q$
if and only if there is a Baer subline $\ell\in\mathbf{w}_\i\cap \mathbf{w}'$.
 Firstly, assume $\rank(A)=1$ for some $a'$. Then the solutions $Y$ of $YA=0$ forms a Baer subline $\ell\subseteq\{\lg Yg\Delta\rg\mid YA=0\}$ which contains $q\pm 1$ or $q$ nonisotropic points, according $G_{\mathbf{w}_\i}\cap G_{\mathbf{w}'}\cong D_{2(q\pm 1)}$ or $\ZZ_p^m.\ZZ_2$.
 Secondly, assume $G_{\mathbf{w}_\i}\cap G_{\mathbf{w}'}\cong \D_{2(q\pm 1)}$. This implies that $\mathbf{w}_\i\cap \mathbf{w}'=\{\ell,\ell^\perp\}$ consists of $q+2$ points, where $\ell$ is a Baer subline and $\ell^\perp$ is a nonisotropic point. Then there exists two solutions $\lg Y_1\rg$ and $\lg Y_2\rg$ satisfying $Y_\i A=0$ for some $a'\in\{1,a+i\}$, where $\i=1,2$. It follows that for any $k_1,k_2\in\FF_q$, we have that $(k_1Y_1+k_2Y_2)A=0$ for such $a'$, as $|\{1,a+i\mid a\in\FF_q\}|=q+1<|\mathbf{w}\cap\mathbf{w}'|=q+2$.
Since all these solutions obtained from the same $a'$, we have $\rank(A)=1$.

\vskip 3mm
Continue our proof of the lemma.   Suppose that $\mathbf{w}_\i\in \mathcal{B}_0=\{\mathbf{w}_0^g\mid g\in G\}$. Then by Lemma~\ref{wsigma},
\vspace{-5pt}\begin{equation}\label{3}  b^*c^*\in(\FF_{q^2}^*)^3.\end{equation}
All triples  $((b^*,c^*),\lg Y \rg)$ and pairs $(b^*,c^*)$ satisfying  Eq(\ref{Y}), Eq(\ref{3}) and $\rank(A_{(a',b^*,c^*)})=2$ form ${\bf H}(\mathcal{D},g)$ and  ${\bf L}(\mathcal{D},g)$, respectively.
 Moreover, following the above arguments, we know that every pair $(\mathbf{w}_\i, \lg \g'\rg )$ in  ${\bf H'}(\mathbf{w}, \mathbf{w}', \cal D)$  is determined by one such triple $((b^*,c^*)$, $\lg Y \rg)$, while
 every  $\mathbf{w}_\i$ in  ${\bf L'}(\mathbf{w}, \mathbf{w}', \cal D)$  is determined by one such pair $(b^*,c^*)$. Remark that  some solutions  of the system of Eq(\ref{Y}), Eq(\ref{3}) and $\rank(A_{(a',b^*,c^*)})=2$ may give
 some extra   $\mathbf{w}_\i$, which are contained  in ${\cal W'}_{j,k, l}$ where $l\ne 1$ or $k=4$. This finishes our proof. \qed
\vskip 3mm

Followed from the proof of Lemma~\ref{hlD}, one may  have  the following remark.
\begin{rem}\label{type} Taking into account of notation from Lemma \ref{hlD}: $b^*, c^*,  b'\in \{1, b+i\} $, $c'\in \{1, c+i\} $ for $b, c\in \FF_q$ and $\De=(\a,\b,\g)^T$,
  we set  $\mathbf{w}_{(\Delta,b',c')}=\{\lg X[1,(b^*)^{-1},(c^*)^{-1}]\De \rg\mid X\in \FF_q^3\}=\{\lg X[1,b',c']\De \rg\mid X\in \FF_q^3\}$. Then
 every Baer subplane containing $\cal D$  is of the form  $\mathbf{w}_{(\Delta,b',c')}$,  while $\mathbf{w}_{(\Delta,b',c')}$ is said to be of {\it type} $[b', c']$ related to  $\cal D$. Let $\cal W_\cal D$ be the set of all such $\mathbf{w}_{(\Delta,b', c')}$, where $(b',c')\ne (1, 1)$.  Clearly, $|\cal W_{\cal D}|=q^2+2q$ and $\mathbf{w}_{(\Delta,b',c')}\in\mathcal{B}_0$ if and only if $b'c'\in(\FF_{q^2}^*)^3$. Thus ${\cal W_\cal D}$  is a union of three subsets with the respective length $q$, $q$ and $q^2:$
\begin{equation}\begin{array}{ll}
 \cal W_{\cal D,1}=\{ \mathbf{w}_{(\Delta,b+i,1)}\mid b\in \FF_q \},&
    \cal W_{\cal D,2}=\{ \mathbf{w}_{(\Delta,1,c+i)}\mid  c\in \FF_q \}, \\
     \cal W_{\cal D,3}=\{ \mathbf{w}_{(\Delta,b+i,c+i)}\mid b, c\in \FF_q \}.&\end{array}
      \end{equation}
  \end{rem}

Let   $\cal D _0=\{ \lg \a_0\rg ,  \lg \b_0\rg  ,  \lg \g_0\rg \}$ be any orthogonal frame of $\mathbf{w}$ and $E=[1,b',c']$, where $\Del_0=(\a_0, \b_0, \g_0)^T$.
  Following the above  notation as in Lemma 5.8, if $XE\Delta_0=a'Y\Del' $,   for some $X, Y\in \FF_q^3$,  then   $\lg Y\Del'\rg\in\mathbf{w}_{(\Delta_0, b',c')}$, that is to say the point $\lg Y\Delta'\rg$ is in the Bare subplane of type $[b',c']$ corresponding to the frame $\cal D _0$.
 Suppose that $\cal D _0\cap \mathbf{w}'=\emptyset .$    Followed from Lemma 2.7,   for any one, say   $\d \in \{ \a_0 ,  \b_0 ,  \g_0 \}$,  we have $|\PG(\d^\perp)\cap\mathbf{w'}|=1$. Therefore, $\mathbf{w'}$ contains the following set
 $$\Phi(\cal D _0,\mathbf{w'})=\{ \lg (0, 1, \rho )\Delta_0\rg , \lg (\tau, 1, 0)\Delta_0\rg , \lg (\tau, 0, -\rho )\Delta_0\rg \}, $$
  for some $\rho, \tau \in\FF_{q^2}^*$.  Totally,  we have the following six cases:
$$\begin{array}{lll}
&{\rm (i)}\,  \rho^{q+1}\ne -1;  &{\rm (ii)}\,  \rho^{q+1}=-1, \tau^{q+1}\ne -1\,    \tau \not\in \FF_q;\\
&{\rm (iii)}\,   \rho^{q+1}=-1, \tau^2\ne -1,  \tau \in \FF_q; &{\rm (iv)}\,    \rho^{q+1}=\tau^{q+1}=-1, \tau, \t^{-1}\rho \not\in \FF_q;\\
&{\rm (v)}\,  \rho^{q+1}=\tau^{2}=-1,  \rho \not\in \FF_q^*, \tau \in \FF_q; &
{\rm (vi)}\,  \rho^{q+1}=\tau^{q+1}=-1,  \t^{-1}\rho \in \FF_q.\end{array} $$
Let $F=F_1\cup F_2$ be the set of  orthogonal frames of $\mathbf{w}$ which is disjoint from $\mathbf{w}'$, where
 \begin{equation}\label{F1F2}
\begin{array}{lll}
F_1=\{\cal D _0\in F\di \Phi(\cal D _0,\mathbf{w'}) \, {\rm is\,  one\,  of\,  cases\,}  (i), (ii), (iv)\},\\
F_2=\{\cal D _0\in F\di \Phi(\cal D _0,\mathbf{w'}) \, {\rm is\,  one\,  of\,  cases\,}  (iii), (v), (vi)\}.
\end{array}
\end{equation}
Moreover, set ${\bf f}_2=\max\limits_{\mathbf{w}'\in\mathcal{B}_0\setminus\{\mathbf{w}\}}|F_2|$ and
\begin{equation}\label{lhf}
\begin{array}{lll}
 {\bf l}_{F_1}=\max\limits_{\mathcal{D}\in F_1,\mathbf{w}'\in\mathcal{B}_0\setminus\{\mathbf{w}\}}\{|{\bf L}'(\mathbf{w},\mathbf{w}',\mathcal{D})|\}
,&&
{\bf h}=\max\limits_{\mathcal{D}\in F_1\cup F_2,\mathbf{w}'\in\mathcal{B}_0\setminus\{\mathbf{w}\}}\{|{\bf L}'(\mathbf{w},\mathbf{w}',\mathcal{D})|\}.\\

\end{array}
\end{equation}

\begin{lem}\label{B1B2} Two numbers ${\hm B_1}$ and ${\hm B_2}$ have the following  upper bounds:
$$\begin{array}{l}
 {\bf B}_1 \leq \left\{
    \begin{array}{ll}
   \frac{q(q^2-1)}{3}( \frac{2}{3}q^2+q\sqrt{q}+13.34 q)+5.1q^4, & \hbox{$d=1$;} \\
    \frac{q(q^2-1)}{3}(\frac{2}{9}q^2+\frac{16}{9}q\sqrt{q}+12.4367q)+1.7q^4, & \hbox{$d=3$,}
    \end{array}
  \right.\\
   {\bf B}_2  \leq \left\{
    \begin{array}{ll}
    \frac{q(q^2-1)}{3}(q^2+8q), & \hbox{$d=1$;} \\
     \frac{q(q^2-1)}{3}(\frac{1}{3}q^2+\frac{5}{3}q\sqrt{q}+9.67q), & \hbox{$d=3$.}
    \end{array}
  \right. \\
 \end{array}$$
\end{lem}
\demo
Recall that $\mathcal{F}_{\mathbf{w}}$ denotes the set of all orthogonal frames of the Baer subplane $\mathbf{w}$, and that $|\mathcal{F}_{\mathbf{w}}|=\frac{q(q^2-1)}{6}.$ Moreover, from Table 1, each orthogonal frame lies in exactly $\frac{(q+1)^2}{d}$ distinct Baer subplanes. Combining Eq(\ref{B12}) with the preceding arguments, we obtain the following bounds:
$$\begin{array}{l}
\begin{array}{lll}
{\bf B}_1&\leq&
\sum\limits_{\mathcal{D}\in F}
\max\limits_{\mathbf{w}'\in\mathcal{B}_0\setminus\mathbf{w}}
2\{|{\bf L'} (\mathbf{w}, \mathbf{w}', \mathcal{D})|\}
\leq 2(|\mathcal{F}_{\mathbf{w}}|-F_2){\bf l}_{F_1}+(\frac{(q+1)^2}{d}-1)\times2{\bf f}_2\\
&\leq&\frac{q(q^2-1)}{3}{\bf l}_{F_1}+\frac{q^2+2q+1-d}{d}\times2{\bf f}_2,
\end{array}\\
\begin{array}{lll}
{\bf B}_2&\le& \sum\limits_{\mathcal{D}\in F}\max\limits_{\mathbf{w}'\in\mathcal{B}_0\setminus\mathbf{w}}
2\{|{\bf H'} (\mathbf{w}, \mathbf{w}', \mathcal{D})|\}
\leq \frac{q(q^2-1)}{3}{\bf h}.
\end{array}\end{array}$$
Moreover, it will be shown in Lemmas \ref{l-h} and \ref{f2}, that ${\bf f}_2\leq\frac 52q^2+4q-\frac 12$ and
$$
\begin{array}{lll}
{\bf l}_{F_1}\leq \left\{
    \begin{array}{ll}
    \frac{2}{3}q^2+q\sqrt{q}+13.34 q, & \hbox{$d=1$;} \\
    \frac{2}{9}q^2+\frac{16}{9}q\sqrt{q}+12.4367q, & \hbox{$d=3$,}
    \end{array}
  \right.&&
{\bf h}\leq \left\{
    \begin{array}{ll}
     q^2+8q, & \hbox{$d=1$;} \\
     \frac{1}{3}q^2+\frac{5}{3}q\sqrt{q}+9.67q, & \hbox{$d=3$.}
    \end{array}
  \right. \\
\end{array}
$$
Then we get upper bounds in the  lemma  for ${\hm B_1}$ and ${\hm B_2}$ by  substituting these bounds for  ${\bf l}_{F_1}$,  ${\bf h}$ and ${\bf f}_2.$
\qed
\vskip 3mm
We now establish upper bounds for the parameters, beginning with ${\bf l}_{F_1}$ and ${\bf h}$.

\begin{lem} \label{l-h}
  Let ${\bf l}_{F_1}$ and ${\bf h}$ be defined in {\rm Eq(\ref{lhf})}. Then
\vspace{-5pt}$$
 {\bf l}_{F_1} \leq \left\{
    \begin{array}{ll}
    \frac{2}{3}q^2+q\sqrt{q}+13.34 q, & \hbox{$d=1$;} \\
    \frac{2}{9}q^2+\frac{16}{9}q\sqrt{q}+12.4367q, & \hbox{$d=3$,}
    \end{array}
  \right.\quad
   {\bf h}  \leq \left\{
    \begin{array}{ll}
     q^2+8q, & \hbox{$d=1$;} \\
     \frac{1}{3}q^2+\frac{5}{3}q\sqrt{q}+9.67q, & \hbox{$d=3$.}
    \end{array}
  \right. \\
 $$
  \end{lem}
\demo
 Given a orthogonal frame $\mathcal{D}=\{\lg\a\rg, \lg\b\rg, \lg\g\rg\}\in \mathbf{w}$ with $\mathbf{w}=\PG(W)$,  let ${\bf H}(\mathcal{D},g)$ and  ${\bf L}(\mathcal{D},g)$ be  defined in  Eq(\ref{HL}).
For convenience, we introduce the following two subsets of them, respectively:
\begin{equation} \label{H1L1} \begin{array}{lll}
{\bf H}_1(\mathcal{D},g)&=&\{((b^*, c^*),\lg Y \rg)\in  {\bf H}(\mathcal{D},g)\mid a', b^*, c^*\ne 1\};\\
{\bf L}_1(\mathcal{D},g)&=&\{(b^*, c^*)\in  {\bf L}(\mathcal{D},g)\mid a', b^*, c^*\ne 1 \}. \end{array}
\end{equation}
 We now discuss the following three cases, separately. \vskip 3mm
(1) Case $b^*, c^*\ne 1$ and $a'=1$.
\vskip 3mm
 In our case,  $E_A+E_Bi=[1, -b + i, -c + i]$, that is  $E_A=[1, -b, -c]$ and $E_B=[0, 1, 1]$. Then
 $A:=A_{(1, -b+i, -c+i)}=g_AE_B+g_BE_A$,  where  $g:=g_A+g_Bi\in\{{\hm g_1},{\hm g_2}\}$.

  First,  let
$g={\hm g_1}=g_A+g_Bi$ where
$$   g_A={\small \left(
                          \begin{array}{ccc}
                            0 & 1 & m_1 \\
                            e & t & 0 \\
                            t & 0 & m_2 \\
                          \end{array}
                        \right)},\,
 g_B={\small\left(
                          \begin{array}{ccc}
                            0 & 0 & n_1 \\
                            1 & 0 & 0 \\
                            0 & 1 & n_2 \\
                          \end{array}
                        \right)},\,
                        A=g_AE_B+g_BE_A={\small\left(
                          \begin{array}{ccc}
                            0 & 1& m_1-n_1c \\
                            1 & t & 0 \\
                            0 & -b & m_2-n_2c\\
                          \end{array}
                        \right)}, $$
                         \f
where  $m_1n_2-m_2n_1=-1, m_1+n_1i, t\neq 0.$ Then $|A|=n_1bc-m_1b+n_2c-m_2=0$  if and only if    $b=\frac{m_2-n_2c}{n_1c-m_1}$ and $n_1c-m_1\ne 0$.
    Since  $\rank(A)=2$,  for giving $b$ and $c$ the equation $YA=0$ has one solution $\lg Y\rg $, so that the number of triples
    $((b^*,c^*),\lg Y \rg)$ is no more than $q-1$.

Second, let  $g={\hm g_2}=g_A+g_Bi,$  where $E_A=[1,-b,-c],$ $E_B=[0,1,1],$
$$ g_A={\small\left(
                          \begin{array}{ccc}
                            0 & 1 & m_1 \\
                            s & -1 & 0 \\
                            0& 0 & m_2 \\
                          \end{array}
                        \right)},\,   g_B={\small\left(
                          \begin{array}{ccc}
                            0 & 0 & n_1 \\
                            0& 0 & 0 \\
                            1 & s & n_2 \\
                          \end{array}
                        \right)},\, A=g_AE_B+g_BE_A={\small\left(
                          \begin{array}{ccc}
                            0 & 1& m_1-n_1c \\
                            0 & -1 & 0 \\
                            1 & -bs & m_2-n_2c\\
                          \end{array}
                        \right)}, $$
where $m_1n_2-m_2n_1=-s\neq 0,$ $m_1+n_1i\neq 0$. Similarly,  we get $|A|=-n_1c+m_1=0$ if and only if $c=\frac{m_1}{n_1}$ and $n_1\neq 0$.
Note that in this case, $YA=0$ if and only if $\lg Y\rg=\lg(1,1,0)\rg$, so that
the number of triples
    $((b^*,c^*),\lg Y \rg)$ is no more than $q$.

\vskip 3mm
(2) Case either $b^*=1$ or $c^*=1$.
\vskip 3mm
 Let $\De=(\a, \b, \g)^T$ and $\mathbf{w}=\{\lg X\De\rg\mid X\in\FF_q^3\}$ for $G_{\mathbf{w}}=M$.
 If either $A_{(a',b^*,c^*)}=A_{(a',1,-c+i)}$ or $A_{(a',-b+i,1)}$, then by Lemma \ref{lem-inv-bij}, we have that $G_{\mathbf{w}_\i}\cap M\cong \D_{2(q\pm1)}$, for the respective $\mathbf{ w}_\i
 =\{\lg X[1,b^*, c^*]\De\rg\mid X\in\FF_q^3\}$.
By Remark \ref{type} and Table \ref{number}, there are at most the number  $\frac{2q}{d}$  (resp. $\frac{6q}{d}$) contributes to ${\bf L}(\mathcal{D},g)$ (resp. ${\bf H}(\mathcal{D},g)$).
\vskip 3mm
(3) Case  $a', b^*, c^*\ne 1$:  Then $((b^*,c^*),\lg Y \rg)$ (resp. $(b^*,c^*)$) in ${\bf H}(\mathcal{D},g)$ (resp. ${\bf L}(\mathcal{D},g)$) is contained in  ${\bf H}_1(\mathcal{D},g)$ (resp. ${\bf L}_1(\mathcal{D},g)$). It will be shown in Lemmas~\ref{H1_3} and ~\ref{L_1(g)3} that for any $\mathbf{w}'\in\mathcal{B}_0\setminus\mathbf{w}$ with $\mathbf{w}'=\{\lg Yg\Delta\mid Y\in\FF_q^3\rg\},$ we have
 $$
\begin{array}{l}
 \max\limits_{\mathcal{D}\in F}|{\bf H}_1(\mathcal{D},g)|\leq \left\{
    \begin{array}{ll}
    q^2+q, & \hbox{$d=1$;} \\
   \frac{1}{3}q^2+\frac 53q\sqrt{q}+6.67q,& \hbox{$d=3$,}
    \end{array}
  \right.\\
   \max\limits_{\mathcal{D}\in F_1}|{\bf L}_1(\mathcal{D},g)| \leq \left\{
    \begin{array}{ll}
     \frac{2}{3}q^2+q\sqrt{q}+10.34 q, & \hbox{$d=1$;} \\
     \frac{2}{9}q^2+\frac{16}{9}q\sqrt{q}+10.77 q, & \hbox{$d=3$.}
    \end{array}
  \right.
 \end{array}$$

Summing upper bounds of the above three cases,  we get upper bounds of ${\bf l}_{F_1}$ and ${\bf h}$ as shown in the lemma.\qed

\subsubsection{Upper bound of $|{\bf H}_1(\mathcal{D},g)|$ with $\mathcal{D}\in F$ } \label{w3}
Taking account of Eq(\ref{bc}), Eq(\ref{HL}) and Eq(\ref{H1L1}), we have
 \vspace{-5pt}$$\begin{array}{lll}
 {\bf H}_1(\mathcal{D},g)&=&\{((-b+i,-c+i),\lg Y \rg)\mid YA=0, Y\in\FF_q^3, Y {Y}^{T}\neq 0,
  \rank(A)=2, \\&&(-b+i)(-c+i)\in(\FF_{q^2}^*)^3\},\\
  {\bf L}_1(\mathcal{D},g)&=&\{(-b+i,-c+i)\mid  ((-b+i,-c+i),\lg Y \rg)\in {\bf H}_1(\mathcal{D},g)\,\, \text{for some}\,
  Y\in\FF_q^3\},
 \end{array}
 $$
where \vspace{-5pt}$$A:=A_{(a+i, -b+i, -c+i)}=a(g_AE_B+g_BE_A)+(g_AE_A+g_BE_B\th),$$
 $g_A+g_Bi\in \{ {\hm g_1}, {\hm g_2} \}$ and $E_A+E_Bi=[1, -b+i, -c+i]$.
Write \vspace{-5pt}\begin{equation}\label{abc} h(a,b,c):=|A|=f(a,c)b-g(a,c).\end{equation}
By Eq(\ref{Y}) and Eq(\ref{3}), we need to solve the system
\vspace{-5pt}\begin{equation} \label{YZ} h(a,b,c)=0\quad {\rm and}\quad  (-b+i)(-c+i)\in (\FF_{q^2}^*)^3.\end{equation}

Obverse that  given  $b$ and $c$,  $h(a, b, c)$ is a polynomial with indeterminate  $a$, which is of degree at most 3.
Suppose that  $\rank(A)=2.$  Then given $a, b, c$, $YA=0$ has one nonzero solution for $\lg Y\rg$.
If  $M_\imath\cap M'$ contains at most three involutions, then $M_\imath\cap M'\cong \ZZ_2$, $D_4$ or $\Alt(4)$, that is $\mathbf{w}_\i\cap \mathbf{w}'$ contains either one or three nonisotropic points. Therefore,
given $b$ and $c$, each  solution $\lg Y\rg$ from ${\hm H}_1(\mathcal{D},g)$ corresponds to one $a$ and thus the number of $\lg Y\rg $ is either one or three with $YY^T\neq 0$. Moreover,
from the above arguments, every couple $b, c$ determines the unique $\mathbf{w}_\i$ (i.e. $M_\imath$).
 To measure ${\bf H}_1(\mathcal{D},g)$, we introduce two sets:
   \vspace{-5pt}\begin{equation}\label{Lambda12}
   \begin{array}{lll}
\Lambda_1:&=&\{(a,b,c)\in \FF_q^3\mid h(a,b,c)=0, f(a,c)\neq 0,
(-b+i)(-c+i)\in(\FF_{q^2}^*)^3\},\\
\Lambda_2:&=&\{(a,b,c)\in \FF_q^3\mid f(a,c)=g(a,c)=0, (-b+i)(-c+i)\in(\FF_{q^2}^*)^3\}. \\
\end{array}
\end{equation}

\begin{lem}\label{H1_3}
Let $\mathcal{D}\in F$ (see Eq(\ref{F1F2})). Then
$ |{\bf H}_1(\mathcal{D},g)|\le q^2+q$ for $d=1;$ and
$ \frac{1}{3}q^2+\frac 53q\sqrt{q}+6.67 q,$ for $d=3$.
 \end{lem}
\demo
(1) Suppose that $d=1$. Then $\FF_{q^2}^*=(\FF_{q^2}^*)^3$ and so $(-b+i)(-c+i)\in(\FF_{q^2}^*)^3$ is automatically true. Since $b=\frac{g(a,c)}{f(a,c)}$ $(f(a,c)\neq 0)$, $(a, b, c)$ has at most $q^2$ choices, that is  $|\Lambda_1|\leq q^2$. Following the coming  Lemma~\ref{ld12},  $|\Lambda_2| \leq q$ and  consequently, $|\mathbf{H}_1(g)| \leq q^2 + q$.
\vskip 3mm
(2) Suppose that $d=3$. Then $\Lambda_{1}=\Lambda_{11}\cup \Lambda_{12}$, where
\vspace{-5pt}\begin{equation}
\begin{array}{lll}
\Lambda_{11}:=\{(a,b,c)\in\Lambda_1\mid b=-c\}&\text{and}&
\Lambda_{12}:=\{(a,b,c)\in\Lambda_1\mid b\neq-c\}.
\end{array}
\end{equation}
Suppose that   $b=-c$ and $f(a, c)\ne 0$.  Then we have at most $q$ pairs $(b,c)$, where every  pair corresponds at most three $\lg Y\rg$ such that $((b^*,c^*),\lg Y\rg)\in{\hm H}_{1}(\mathcal{D},g)$, as $\mathbf{w}_\i\in {\cal W'}_{j,k, 1}$ where $k=2, 3$. So that $|\Lambda_{11}|\leq 3q.$

Following Lemmas~\ref{ld12} and ~\ref{ld3}, we have
$|{\bf H}_1(g)|\le  3q+|\Lambda_2|+|\Lambda_{12}|\le 3q+q+(\frac{1}{3}q^2+\frac 53q\sqrt{q}+\frac 83q)\le \frac{1}{3}q^2+\frac 53q\sqrt{q}+\frac {20}3q.$
\qed

\begin{lem}\label{ld12} With the above notation, we have  $|\Lambda_2|\le q$;    and for any $a_0 \in \FF_q\setminus \{-e\}$,  we have  $(f(a_0, c), g(a_0, c))=1$, where
$e$ is an  entry  of ${\hm g_1}$.
 \end{lem}
\demo   (1) Set $g={\hm g_1}$. Then
$$\begin{array}{lll} h(a,b,c)&=&{\small\left|
                     \begin{array}{ccc}
                       0 & a-b & (am_1+\th n_1)-(m_1+an_1)c \\
                       e+a & at-tb & 0 \\
                       t & \th-ab & (am_2+\th n_2)-(m_2+an_2)c \\
                     \end{array}
                   \right|}\\
 &=&\left|A_1+aB_1, (aA_2+\th B_2)-(A_2+aB_2)b, (aA_3+\th B_3)-(A_3+aB_3)c\right|\\
 &=&f(a,c)b-g(a,c),  \end{array}$$
\vskip -10pt \f where   $e,s,t,m_1,m_2,n_1,n_2\in\FF_q$,  with $m_2n_1-m_1n_2=1, t\ne 0$; $f(a,c):=A+tB-aC$ and $g(a,c):=aA+atB-\th C$ with
\vspace{-3pt}\begin{equation} \label{ABC1}\begin{array}{ll} &A:=(e+a)[(am_2+\th n_2)-(m_2+an_2)c], \\ &B:=t[(am_1+\th n_1)-(m_1+an_1)c],\\ &C:=(e+a)[(am_1+\th n_1)-(m_1+an_1)c].\end{array} \end{equation}

(1.1) To measure  $|\Lambda_2|$,   suppose $f(a,c)=g(a,c)=0$. It follows from $af(a, c)-g(a, c)=0$ that $(a^2-\th)C=0$  and so $C=0$. Therefore, $a=-e$ or $c=\frac{am_1+\th n_1}{m_1+an_1}$. If $a=-e$, then $C=A=0$ and so $B=0$, and $c=\frac{am_1+\th n_1}{m_1+an_1}=\frac{-em_1+\th n_1}{m_1-en_1}$. If
$c=\frac{am_1+\th n_1}{m_1+an_1}$, then  $B=C=0$, so that $A=0$, that is, $a=-e$, as $(am_2+\th n_2)-(m_2+an_2)c\neq 0$ for $c=\frac{am_1+\th n_1}{m_1+an_1}$. Thus every triple of $\Lambda_2$ is of the form
$(-e, b, \frac{-em_1+\th n_1}{m_1-en_1})$ where $b\in \FF_q$. So $|\Lambda_2|\le q.$

(1.2) Since $\deg(f(a_0, c), c), \deg(g(a_0, c), c)\le 1$ for any $a_0\in\FF_q$, they have a common divisor if and only if they have a common root. According to the proof of (1.1), this occurs if and only if $a_0 = -e$. Therefore,
 $f(a_0, c)$ and $g(a_0, c)$ are coprime provided $a_0\ne -e$.

\vskip 3mm (2)  Set  $g={\hm g_2}.$ Then
$$\begin{array}{lll} h(a,b,c)&=&{\small\left|
                     \begin{array}{ccc}
                       0 & a-b & (am_1+\th n_1)-(m_1+an_1)c \\
                       s & -a+b & 0 \\
                       a & \th s-abs & (am_2+\th n_2)-(m_2+an_2)c \\
                     \end{array}
                   \right|} \\
&=&\left|A_1+aB_1, (aA_2+\th B_2)-(A_2+aB_2)b, (aA_3+\th B_3)-(A_3+aB_3)c\right|\\
 &=&f(a,c)b-g(a,c),  \end{array}$$
\vskip -10pt \f  where  $m_2n_1-m_1n_2=s, s\ne 0$;  $f(a,c):=A-B-asC,$  $g(a,c):=aA-aB-\th sC$, where
 \vspace{-5pt}\begin{equation} \label{ABC2}\begin{array}{lll}
 &A:=s[(am_2+\th n_2)-(m_2+an_2)c],&
 B:=a[(am_1+\th n_1)-(m_1+an_1)c],\\
 &C:=s[(am_1+\th n_1)-(m_1+an_1)c].&\end{array} \end{equation}

Again suppose that  $f(a,c)=g(a,c)=0$. Then we have $C=0$, that is,  $c:=\frac{am_1+\th n_1}{m_1+an_1}$, as $s\neq 0$. Thus,  $A\neq 0$ but $B=C=0$, which is impossible, that is  $\Lambda_2=\emptyset.$
 Similar to  (1.2),   $f(a_0, c)$ and $g(a_0, c)$ are co-prime for any $a_0$, as $\Lambda_2=\emptyset.$\qed

\vskip 3mm
To measure  $|\Lambda_{12}|$,  assume $b\ne -c$ and $f(a, c)\ne 0$.  Set  $(-b+i)(-c+i)=(x+yi)^3$ for some $x,y\in\FF_q$ and $z=\frac{x}{y}$. Then our  system  Eq(\ref{YZ}) can be
 described as
\vspace{-5pt}\begin{equation}  b=\frac{g(a, c)}{f(a, c)}\,\,  {\rm and}\,  \, \label{z}(bc+\th)(3z^2+\th)=-(z^3+3z\th)(b+c).\end{equation}
  \vspace{-5pt}  Inserting $b=\frac{g(a, c)}{f(a, c)}$ to the second equation, we get
\begin{equation} \label{h'} h'(a,c,z):=F(a, c)(z^3+3z\th)+G(a,c)(3z^2+\th)=0,\end{equation}
where $F(a,c):=g(a,c)+f(a,c)c$ and $G(a,c):=g(a,c)c+f(a,c)\th$, $f(a,c)\neq 0$ and $b+c=\frac{g(a,c)}{f(a,c)}+c\neq 0$. Now, we set
\vspace{-5pt}\begin{equation} \label{ld12'} \Lambda_{12}':=\{(a,c,z)\in \FF_q^3\mid h'(a,c,z)=0\}.\end{equation}

Before going on we insert two propositions. The first one  will play a key role in this paper. It relies on Weil's deep work on algebraic geometry and   number theory. To use it,
 $h'(a_0,c,z)$ is necessarily  absolutely irreducible, for certain $a_0\in \FF_q$. This will be proved in Lemma~\ref{absolute}.
\vspace{-5pt}  \begin{prop} {\rm\cite[Theorem 6.57] {LN}} \label{weil}
Let $g(x)\in\FF_q[x]$ with $\deg(g)=d$ be such that $f(x,y)=y^2-g(x)$ is absolutely  irreducible polynomial. Then the number $N$ of solutions of $f(x,y)=0$ in $\FF_q^2$ satisfies
$|N-q|\leq (d-1)\sqrt{q}.$
\end{prop}
\begin{prop}   {\rm\cite[Corollary 1.15]{Hirs}} \label{disc}
  Let $f(x)=\sum_{\i=0}^3a_ix^i$ be a polynomial over $\FF_q$ with three roots $x_1, x_2, x_3$ in some field $K\ge \FF_q$. Let $D=\Pi_{\i\le \j}(x_\i-x_\j)^2$ be its discriminant. Then
    $f(x)$ has multiple roots, exactly one root or three distinct roots in $\FF_q$, only if  $D$ is zero, a nonzero nonsquare or a nonzero square element of $\FF_q$, respectively.
    \end{prop}

\begin{lem}\label{3-1}  $|\Lambda_{12}|\leq\frac 13 |\Lambda_{12}'|$.\end{lem}
\demo  Set $\mu=\frac{z^3+3z\th}{3z^2+\th}$ and  $m(z,\mu)=z^3-3\mu z^2+3\th z-\mu\th$. For any given $\mu_0\in \FF_q$, one may  check that  the discriminant
of $m(z, \mu_0)$  is $(-3)\th 2^23^2(\mu_0^2-\th)^2$, a quare, as $q\equiv 2 \pmod 3$. By Proposition~\ref{disc}, three distinct $z$ give a unique $\mu$. This implies that each triple $(a, b, c)\in\Lambda_{12}$, gives three solutions $(a,c,z_\i)$ where $\i=1, 2, 3$ in Eq(\ref{ld12'}). Since $z$ is free whenever $f(a,c) = g(a,c) = 0$, each solution $(a, c)$ corresponds to multiple triplets $(a, c, z)$ in $\Lambda_{12}'$. However, these do not yield valid solutions $(a, b, c)$ in $\Lambda_{12}$. Consequently, $|\Lambda_{12}|\leq\frac 13 |\Lambda_{12}'|$.\qed

\begin{lem} \label{absolute}
Suppose $d=3$. Set $\daleth=\{a_0\mid (F(a_0,c),G(a_0,c))= c^2-\th\}$. Then
\vspace{-5pt}\begin{enumerate}[itemsep=5pt]
\item[\rm(1)]  $|\daleth|\leq 3$, and $F(a_0,c)$ and $G(a_0,c)$ are coprime, for any  $a_0\in\FF_q\setminus{\{-e,\Delta\}}$;
\item[\rm(2)]   for any $a_0\in\FF_q\setminus{\{-e,\daleth \}}$,  $h'(a_0,c,z)$ is absolutely irreducible.
\end{enumerate} \end{lem}
\demo  (1)  Set $m(a_0,c):=(F(a_0,c),G(a_0,c))$. If $a_0=-e$, then from  Eq(\ref{ABC1}),  we have $f(a,c)=tB,$ and $g(a,c)=atB$ so that $m(a_0,c)\ne 1$. Let  $a_0\in\FF_q\setminus{\{-e\}}$ and   assume that $m(a_0,c)\neq 1$.
 By Eq(\ref{h'}), $F(a,c):=g(a,c)+f(a,c)c$ and $G(a,c):=g(a,c)c+f(a,c)\th$.  Then it follows that $cF(a_0,c)-G(a_0,c)=(c^2-\th)f(a_0,c).$ Since    $(f(a_0,c),g(a_0,c))=1$ by Lemma~\ref{ld12},
 we get  $m(a_0,c)\mid c^2-\th$. This implies that $m(a_0,c)=c^2-\th$, as $\theta\in\FF_q$ is non-square and $c^2-\theta\in\FF_q[c]$ is irreducible.

 \vskip 3mm
 Suppose  $g={\hm g_1}$.
Following the proof of Lemma \ref{ld12}, we obtain
$$
\begin{array}{lll} F(a_0,c)&=&(a_0+c)A+(a_0+c)tB-(a_0c+\th)C,\\
&=&c^2[(m_1+a_0n_1)(-t^2+a_0e+a_0^2)-(e+a_0)(m_2+a_0n_2)]\\
&&+c[(e+a_0)(-a_0^2n_2+\th n_2)+(m_1+a_0n_1)(\th e+a_0\th-a_0t^2)\\
&&+(a_0m_1+\th n_1)(-a_0e-a_0^2+t^2)]\\
&&+a_0(e+a_0)(a_0m_2+\th n_2)+(a_0m_1+\th n_1)(a_0t^2-\th e-\th a_0)
\end{array}\\$$
$$\begin{array}{lll}
G(a_0,c)&=&(a_0c+\th)A+(a_0c+\th)tB-(a_0+c)\th C\\ &=&c^2[-a_0(e+a_0)(m_2+a_0n_2)-(m_1+a_0n_1)(a_0t^2-\th e-\th a_0)]\\
&&+c[(e+a_0)(a_0^2m_2-\th m_2)+(a_0t^2-\th e-\th a_0)(a_0m_1+\th n_1)\\
&&+(m_1+a_0n_1)\th(a_0 e+a_0^2-t^2)]\\
&&+(-a_0 e-a_0^2+ t^2)\th(a_0m_1+\th n_1)+\th(e+a_0)(a_0m_2+\th n_2).
\end{array}$$ \vspace{-5pt}
Since $m(a_0,c)= c^2-\th$, the coefficient of $c$ for both $F(a_0,c)$ and $G(a_0,c)$ is zero,
that is
$$\begin{array}{lll}
&&-a_0^3(m_1+n_2)-a_0^2(en_2+m_1 e+n_1t^2)+a_0\th(n_2+m_1)
+(\th n_2e+m_1\th e+\th n_1t^2)=0,\\
&&a_0^3(m_2+n_1\th)+a_0^2(em_2+m_1t^2+\th n_1 e)-a_0\th(m_2+\th n_1)
-\th(em_2+\th n_1e+m_1t^2)=0.
\end{array}$$
If at least one of $m_1-n_2$ and  $m_2-n_1\th$ is not zero, then  $a_0$ has  at most $3$ choices; otherwise,  we get
 $n_1a_0^2+\th n_1=0$ and $m_1a_0^2+m_1\th=0,$
forcing  that $a_0$ has  at most $2$ choices, noting  $m_2n_1-m_1n_2=1$.
Therefore,
 $|\daleth|\leq 3$.
 \vskip 3mm
 Suppose  $g={\hm g_2}$.
Following the proof of Lemma \ref{ld12}, we obtain
$$\begin{array}{lll} F(a_0,c)&=&(a_0+c)A-(a_0+c)B-(a_0c+\th)sC\\
&=&c^2[(m_1+a_0n_1)(a_0+a_0s^2)-s(m_2+a_0n_2)]\\
&&+c[-(a_0m_1+\th n_1)(a_0s^2+a_0)+(m_1+a_0n_1)(\th s^2+a_0^2)+sn_2(\th-a_0^2)]\\
&&+a_0s(a_0m_2+\th n_2)-(a_0^2+\th s^2)(a_0m_1+\th n_1)
\end{array}$$
and
$$\begin{array}{lll}
G(a_0,c)&=&(a_0c+\th)A-(a_0c+\th)B-(a_0+c)s\th C\\ &=&c^2[-a_0s(m_2+a_0n_2)+(a_0^2+s^2\th)(m_1+a_0n_1)]\\
&&+c[sm_2(a_0^2-\th)-(a_0^2+s^2\th)(a_0m_1+\th n_1)
+a_0\th(1+ s^2)(m_1+a_0n_1)]\\
&&-a_0\th(1+s^2)(a_0m_1+\th n_1)+\th s(a_0m_2+\th n_2).
\end{array}$$
Since $m(a_0,c)= c^2-\th$, the coefficient of $c$ for both $F(a_0,c)$ and $G(a_0,c)$ is zero,
that is
$$\begin{array}{lll}
&&n_1a_0^3-a_0^2(s^2m_1+sn_2)-\th n_1a_0+\th s^2m_1+\th sn_2 =0,\\
&&-m_1a_0^3+a_0^2(sm_2-\th s^2n_1)+a_0\th m_1-sm_2\th-s^2n_1\th^2=0.
\end{array}$$
Since $m_2n_1-m_1n_2=1$, we get that $a_0$ has at most $3$ choices.
Therefore,
 $|\daleth|\leq 3$.

\vskip 3mm
(2) Suppose  that  $a_0\in\FF_q\setminus{\{-e,\daleth\}}$.  Assume on the contrary that $h'(a_0,c,z)$ is  reducible. Since
$(F(a_0,c),G(a_0,c))=1$, write  $h'(a_0, c, z)=h_1(c, z)h_2(c,z)$, where
$$h_1(c,z)=g_1(c)z+f_1(c)\quad {\rm and}\quad h_2(c,z)=g_2(c)z^2+u(c)z+f_2(c),$$
 \f where, no loss, assume both $g_1(c)$ and  $f_1(c)$  are monic.
Thus, $$\begin{array}{lll}&&F(a_0,c)z^3+3G(a_0,c)z^2+3\theta F(a_0,c) z+\theta G(a_0,c)=h'(a_0,c,z)=h_1(c,z)h_2(c,z)\\ &=&g_1(c)g_2(c)z^3+(f_1(c)g_2(c)+g_1(c)u(c))z^2+(g_1(c)f_2(c)+f_1(c)u(c))z+f_1(c)f_2(c),\end{array} $$
where
$F(a_0,c):=g(a_0,c)+f(a_0,c)c$ and $G(a_0,c):=g(a_0,c)c+f(a_0,c)\th$,
 which implies
$$\begin{array}{lllll}  F(a_0,c)&=&g_1(c)g_2(c),
  \hskip 2.5cm \theta G(a_0,c)&=&f_1(c)f_2(c),\\
  3G(a_0,c)&=&f_1(c)g_2(c)+g_1(c)u(c),\quad
   3\theta F(a_0,c)&=&g_1(c)f_2(c)+f_1(c)u(c).\end{array} $$
    Since $(F(a_0,c),G(a_0,c))=1$, we get that $(f_\i(c), g_\j(c))=1$ for $\i, \j=1, 2$. Then $f_1(c), g_1(c)\mid u(c)$ so that $f_1(c)g_1(c)\mid u(c)$. Set $u(c)=f_1(c)g_1(c)u_1(c)$. Then
    we have \vspace{-5pt}$$3f_2(c)=\th g_2(c)+
    \th g_1(c)^2u_1(c) \quad {\rm and}\quad
   3\th g_2(c)=f_2(c)+f_1(c)^2u_1(c),$$
 which implies that $\deg(f_1(c)), \deg(g_1(c))\le 1$  and
$$8f_2(c)=u_1(c)(f_1^2(c)+3\th g_1^2(c))\quad {\rm and}\quad
 8\th g_2(c)=u_1(c)(3f_1^2(c)+\th g_1^2(c)).$$
\f Check that this is impossible and so $h'(a_0,c,z)$ is absolute irreducible.
\qed

\begin{lem}\label{ld3}
Suppose $d=3$. Then $|\Lambda_{12}|\leq \frac{1}{3}q^2+\frac 53q\sqrt{q}+\frac 83q.$
\end{lem}
\vspace{-5pt} \demo Given $a_0\in \FF_q$,  rewrite $h'(a_0,c,z)=\sum_{\i=0}^2 h_\i(a_0, z)c^\i$, as  $f(a, c)\ne 0$ and $\deg(h'(a_0,$ $c,z),c)=2$.
 Let $l(a_0, z)=h_1(a_0, z)^2-4h_2(a_0, z)h_0(a_0, z).$
Let $a_0\in \FF_q\setminus (\daleth \cup  \{-e\})$. By Lemma~\ref{absolute},  $h'(a_0, c, z)$ is absolutely irreducible, which implies   that  $x^2-l(a_0, z)$  is absolutely irreducible,
Noting $\deg(l(a_0, z),z)=6$,  by Proposition~\ref{weil}, we have at most   $q+5\sqrt{q}$  solutions $(z, x)$ for $x^2=l(a_0, z)$,   which implies
  there are at most   $q+5\sqrt{q}$  solutions $(c, z)$ for $h'(a_0, c, z)=0$.

Suppose $a_0\in\daleth \cup  \{-e\}$. Then $h'(a_0, c, z)=0$ has  at most the number $3q$ of  solutions  for  $(c,z)$.  Set $k=|\daleth \cup  \{-e\}|$. Then $k\le 4$.
Finally by Lemma~\ref{3-1} we have
 \vspace{-5pt}$$\hskip 25mm|\Lambda_{12}|\leq  \frac 1{3}|\Lambda_{12}'|\leq \frac{1}{3}((q-k)(q+5\sqrt{q})+3kq)\le \frac{1}{3}q^2+\frac 53q\sqrt{q}+\frac 83q.\hskip 20mm \Box $$

\subsubsection{Upper bounds of $|{\bf L}_1(\mathcal{D},g)|$ \, ($\mathcal{D}\in F_1$)\, and  ${\bf f}_2$}\label{w4}
\begin{lem} \label{L_1(g)3}
Suppose  $\mathcal{D}\in F_1$. Then
 $ |{\bf L}_1(\mathcal{D},g)|\le \frac{2}{3}q^2+q\sqrt{q}+10.34q$, for $d=1$;\,      $\frac{2}{9}q^2+\frac{16}{9}q\sqrt{q}+10.77 q$ for $d=3$.
   \end{lem}
\demo  Let $h(a, b, c)$ and $h'(a, c, z)$ be defined in Eq(\ref{abc}) and Eq(\ref{h'}), respectively.
To measure $|{\bf L}_1(\mathcal{D},g)|$, set $n_d=\frac{1}{d}|\Pi_d|,$ where  $d=1,3,$  and
\vspace{-5pt}\begin{equation}\label{nd}
\begin{array}{lll}
\Pi_1&=&\{(b_0,c_0)\in\FF_q^2\mid  h(a,b_0,c_0)=0\,\, {\rm has\, no\, three \, distinct\, solutions \, for\, } a\in\FF_q,  \\
&& {\rm where}\, f(a, c)\ne 0 \};\\
\Pi_3&=&\{(c_0,z_0)\in\FF_q^2\mid  h'(a,c_0,z_0)=0\,\, {\rm has\, no\, three\, distinct\,  solutions \, for\, } a\in\FF_q,\\
&& {\rm where}\, f(a, c)\ne 0, \,  g(a,c)+f(a,c)c \neq0\}.
\end{array}
\end{equation}
 Form the arguments before Eq(\ref{Lambda12}),  given $b$ and $c$, each  solution $\lg Y\rg$ corresponds to one $a$ and thus the number of $\lg Y\rg $ is   three, provided there are three distinct solutions $a$
 of system  of Eq(\ref{3}),   Concretely, given $(b,c)$ (resp. $(c,z)$),    if $d=1$ (resp. $d=3$),  the equation  $h(a, b, c)=0$ (resp. $h'(a, c, z)=0$)  has three distinct solutions $a$, recalling that
 there distinct $(c, z_\i)$ give same  solutions $a$ (one of three) of $h'(a, c, z)=0$.

 Recalling ${\bf H}_1(\mathcal{D},g)\subseteq\Lambda_1\cup \Lambda_2$ (resp. $\Lambda_{11}\cup \Lambda_{12}\cup \Lambda_{2}$) for $d=1$  (resp. $d=3$), we get
   $$|{\bf L}_1(\mathcal{D},g)|\leq \frac{|\Lambda_1|-n_1}{3}+n_1+|\Lambda_2|, d=1  \, {\rm and}\, |{\bf L}_1(\mathcal{D},g)|\leq \frac{|\Lambda_{12}|-n_3}{3}+n_3+|\Lambda_{11}|+|\Lambda_{2}|, d=3.$$
 Inserting upper bounds of $\Lambda_{1}$ and  $\Lambda_{2}$ for $d=1$;  $\Lambda_{11}$,  $\Lambda_{12}$ and $\Lambda_{2}$  for $d=3$  in the proof of Lemma~\ref{H1_3};  $n_1$ in Lemma~\ref{n2d}; and  $n_3$ in  Lemma \ref{n2}, we get that
\vspace{-5pt}
$$\hskip 1cm\begin{array}{lll} |{\bf L}_1(\mathcal{D},g)|&\le &\frac 13q^2+q+\frac{2}{3}(\frac 12 q^2+\frac 32q\sqrt{q}+14q-36\sqrt{q})\leq\frac{2}{3}q^2+q\sqrt{q}+10.34q, \, d=1; \,  \\
|{\bf L}_1(\mathcal{D},g)|&\leq &\frac 13(\frac 13q^2+\frac 53q\sqrt{q}+\frac{8}{3}q)+3q+q+\frac 23(\frac 16q^2+\frac {11}6q\sqrt{q}+\frac{26}{3}q+\frac{20}{3}\sqrt{q}), \\
&\le& \frac{2}{9}q^2+\frac{16}{9}q\sqrt{q}+10.77q,\,  d=3.\hskip 6.5cm\Box
\end{array} $$

First we derive  the existence of three distinct solutions $a$ of $h(a, b, c)=0$ for some paris $(b, c) $ and $h'(a, c, z)=0$ for some pairs $(c, z)$,
 which  will be a  start point to  measure  $n_3$ in Lemma~\ref{n2d} and $n_1$ in  Lemmas~\ref{n2}.

Let $\mathcal{D}_{01}=\{\lg\a_1\rg, \lg\b_1\rg , \lg\g_1\rg\}$ and  $\mathcal{D}_{02}=\{\lg\a_2\rg, \lg\b_2\rg, \lg\g_2\rg\}$  be two representatives of two classes of  orthogonal frame of $\mathbf{w}$,  which can transpired by
  an element  $g'\in \PSU(3,q)\setminus M=\SO(3,q)$, where set  $\Delta_{0,\i}=(\a_\i,\b_\i,\g_\i)^T$  for $\i=1, 2$. Then $\lg\tau_{\a_\i},\tau_{\b_\i},\tau_{\g_\i}\rg\cong\D_4$ for $\i=1,2$, and these subgroups represent the two conjugacy classes of $\D_4$ in $M$.
Recall that $F=F_1\cup F_2$ be the set of  orthogonal frames  $\mathcal{D}_0$  of $\mathbf{w}$ with  ${\mathcal{D}}_0\cap \mathbf{w}'=\emptyset $, see Eq(\ref{F1F2}).

\begin{lem}\label{D4Z2} Let $\mathcal{D}_0\in F_1$. Then  there exist at least  $\frac 1{2d}(q-5\sqrt{q})-\frac 13$   pairs $(b, c)\in\FF_q^2$ relative to  $\mathcal{D}_0$  such that $h(a,b,c)=0$ has three distinct solutions $a$, where $d=1, 3$;  and  there exist at least
 $\frac 12(q-5\sqrt{q})-1$   pairs $(c, z)\in\FF_q^2$ such that $h'(a,c,z)=0$ has three distinct solutions $a$ where $d=3$.
 \end{lem}
\demo  Let $\mathcal{D}_0\in F_1$.  Then by Eq(\ref{F1F2}), $\Phi(\cal D _0,\mathbf{w'})$ is one of  cases  $(i)$, $(ii)$ and $(iv)$.
 \vskip 3mm
  {\it Case (i):   $\rho^{q+1}\ne -1$.}
 \vskip 2mm
 In this case,  $\lg (0,1, \rho )\Delta_0\rg\in \cal N$ and by Propsosition\ref{lem_obser}.(4),   $W'$ has a basis as follows:
 $$v_1=(0,1,\xi_1)\Delta_0,\,v_2=\d_2(1,\xi_2,t\xi_1\xi_2)\Delta_0,\,v_3=\d_3(1,s\xi_2,st\xi_1\xi_2)\Delta_0,$$
where $s,t\in\FF_q$, $s\ne 1$,   $\xi_1=\rho \ne 0$, satisfying  $(v_\i,v_\j)=0$ for $\i\ne \j$; and  $\ne 0$  for $\i=\j$.
Since $\mathcal{D}_0\cap \mathbf{w}'=\emptyset,$ we get   $\{\d_2,\d_3\}\cap (\FF_{q^2}\setminus\FF_q)\neq \emptyset$.
 This allows us rechoose a basis in $\lg v_2, v_3\rg_{\FF_q}$, denoted by $v_2, v_3$ again such that $\d_2 =m'+i$ and $\d_3=n'+i$, where $m',n'\in\FF_q$.
 Then every point in $\PG(v_1^\perp) \setminus \lg v_3\rg $ is of the form
 \vspace{-5pt}$$\begin{array} {lll} \lg v_2+kv_3\rg &=&\lg (\d_2+k\d_3, \d_2\xi_2+ks\d_3\xi_2, t\d_2\xi_1\xi_2+kst\d_3\xi_1\xi_2)\rg \\
 &=& \lg (1,\frac{\d_2+ks\d_3}{\d_2+k\d_3}\xi_2  , \frac{\d_2+ks\d_3}{\d_2+k\d_3}t\xi_1\xi_2)\rg ,\end{array} $$
 where  $k\in \FF_q$  and $\d_2+k\d_3\ne 0$.  Let $b':=\frac{\d_2+ks\d_3 }{\d_2+k\d_3}\xi_2$. Then
$\lg v_2+kv_3\rg =\lg (1, b',b't\xi_1)\Delta_0\rg $. Set  $\Theta=\{ k\in \FF_q\di b'\not\in \FF_q, \d_2+k\d_3\ne 0\}$. Then $|\Theta|\ge q-2$.
In fact, we may set  $g\in \{{\hm g_1}, {\hm g_2}\}$, where $\xi_1=m_1+n_1i$ in both cases.  Set
  $b'=m_2(b+i)$ so that    $b't\xi=m_2t(bm_1+\th n_1+(bn_1+m_1)i).$
Then $\lg v_1\rg , \lg v_2+kv_3\rg \in \mathbf{w}_{(\Delta_0, b+i, \frac{bm_1+n_2\th }{m_1+bn_1}+i)}$, where $m_1+bn_1\ne 0.$

\vskip 3mm

(a) Suppose $d=1$. Remind   $k\in\Theta$. Since there are at most two  values of $k$  such that   $m_1+bn_2=0$,
  we get  at least $\frac{q-4}2$  Baer subplanes  $\mathbf{w}_{(\Delta_0, b+i, \frac{bm_1+n_2\th }{m_1+bn_1}+i)}$  which intersect  a orthegonal frame
$\{ \lg v_1\rg,   \lg v_2+kv_3\rg , \lg v_3'\rg \}$ with   $\mathbf{w}'$, where $v_3'\in \lg v_1, v_2+kv_3\rg ^\perp \cap W'$. Equivalently,   we get  at least $\frac{q-4}2$ subgroups $M_\i$ such that   $M_\i\cap M'\cong D_4$.

 Suppose that $a'=1$ for some of these pairs.  Checking  the proof of  in Lemma 5.8.(2),  we get that $c=\frac{m_1b+m_2}{n_1b+m_1}$ where $n_1b+m_1\ne 0$,  for $g={\hm g_1}$; and    $b\in \FF_q$ and $c=\frac{m_1}{n_1}$ where  $n_1\neq 0$, for $g={\hm g_2}$.   Solving $\frac{m_1b+m_2}{n_1b+m_1}=\frac{bm_1+n_2\th }{m_1+n_1}$ and  $\frac{m_1}{n_1}=\frac{bm_1+n_2\th }{m_1+n_1}$, respectively,
we obtain  at most one solution $b$ (at most two solutions $k$). Therefore,  for cases $b^*,c^*\neq1$, there exist at least  $\frac{q-7}{2}$   pairs $(b, c)\in\FF_q^2$ such that $h(a,b,c)=0$ has three distinct solutions $a$.
\vskip 3mm
(b) Suppose that $d=3$. Then we have to  pick up all the Baer subplanes  $\mathbf{w}_{(\Delta_0, b+i, \frac{bm_1+n_2\th }{m_1+bn_1}+i)}$ which are contained in $\mathcal{B}_0$, meaning
$(\frac{\d_2+ks\d_3}{\d_2+k\d_3})^2t\xi_1\xi_2^2\in\lg \xi^3\rg$, that is
\vspace{-5pt}$$(c_1+c_2i) \frac{\d_2+k\d_3 }{\d_2+ks\d_3}=(x_1+y_1i)^3,$$ for some $x_1,y_1\in\FF_q$,
  where  $t\xi_1\xi_2^2=c_1+c_2i$ for some $c_1,c_2\in\FF_q$, equivalently,
\vspace{-5pt}$$(c_1+c_2i) \frac{\d_2+ks\d_3 }{\d_2+k\d_3}=(x+yi)^3,$$ for some $x,y\in\FF_q$,
  where  $t\xi_1\xi_2^2=c_1+c_2i$ for some $c_1,c_2\in\FF_q$, that is
 \vspace{-5pt}$$(c_1+c_2i)\frac{(m'+sn'k)+(1+sk)i}{(m'+kn')+(1+k)i}=(x^3+3xy^2\th)+(y^3\th+3x^2y)i,
 \,\,\text{and so}$$
 \vspace{-5pt}$$h(z, k):=f(k)(z^3+3z\th)-g(k) (3z^2+\th )=0,$$
 where  $f(k)=c_2d(k)+c_1 e(k),  g(k)=c_1d(k)+c_2\th e(k)$ and $z:=\frac{x}{y}$, where
 \vspace{-5pt}$$d(k)=s(n'^2-\th )k^2+(s+1)(m'n'-\th)k+m'^2-\th \,{\rm and}\, e(k)=(m'-n')(s-1)k.$$
 Clearly, $(d(k), e(k))=1$ so that $(f(k), g(k))=1$, noting that $(\deg(f(k)), \deg(g(k)))\in\{(2,2),(1,2),(2,1)\}$.
   Moreover,  $h(z, k)$ is absolutely irreducible. For the contrary,  suppose   that $h(z, k)$ may be  factorized on some $K\ge \FF_q$, that is
   \vspace{-5pt}$$\begin{array}{ll}
   &h(z,k)=(g_1(k)z+f_1(k))(g_2(k)z^2+u(k)z+f_2(k)),\quad\text{so \, that }\\
   &f(k)z^3-3g(k)z^2+3\th f(k)z -g(k)\th\\ &=g_1(k)g_2(k)z^3+(f_1(k)g_2(k)+g_1(k)u(k))z^2+(g_1(k)f_2(k)+f_1(k)u(k))z+f_1(k)f_2(k),\end{array} $$
\vskip -5pt \f which implies
  \vspace{-5pt}$$\begin{array}{llll}
   &f(k)=g_1(k)g_2(k),
  &-\th g(k)=f_1(k)f_2(k),\\
  &-3g(k)=f_1(k)g_2(k)+g_1(k)u(k),
   &3\th f(k)=g_1(k)f_2(k)+f_1(k)u(k).\end{array} $$
    Since $(f(k), g(k))=1$, we get $(f_\i(k), g_\j(k))=1$ for $\i, \j=1, 2$. Then $f_1(k), g_1(k)\mid u(k)$ so that $f_1(k)g_1(k)\mid u(k)$. Set $u(k)=u_0(k)f_1(k)g_1(k)$. Then
    we have
    \vspace{-5pt}$$\begin{array}{lll}
    &3f_2(k)=g_2(k)\th+u_0(k)g_1(k)^2\th,&
   3\th g_2(k)=f_2(k)+u_0(k)f_1(k)^2.
   \end{array}$$
   Check that  is impossible, as $(\deg(f(k)), \deg(g(k)))\in\{(2,2),(1,2),(2,1)\}$.

   Note that $h(z,k)=0$ has solutions if and only if $l(x,z):=x^2-d(z)=0$ has solutions in $\FF_q^2$, where $d(z)\in\FF_q[z]$ is the discriminant of $h(z,k)$ corresponding to $k$.
   Since $h(z, k)$ is absolutely irreducible and $\deg(d(z))=6$, by Proposition~\ref{weil}, the number $n(z,k)$ of solution of $l(x,z)=0$ is at least $q-5\sqrt{q}$ so that
     the number $n(k)$ of solution $k\in\FF_q$ of  $h(z,k)=0$ is    $\frac 12(q-5\sqrt{q})$.
     Repeated the arguments in (a),  by removing two possible isotropic points in $\PG(v_1^\perp )$ and two values of $k$ such that $a'=1$, we get   at least  $\frac 12(q-5\sqrt{q})-2$  pairs $(c, z)$ such that   $h'(a,c,z)=0$ has three distinct solutions $a$ and so there are at least  $\frac 16(q-5\sqrt{q})-\frac 23$  pairs $(b,c)$ such that   $h(a,b,c)=0$ has three distinct solutions $a$.

 \vskip 3mm
  {\it
  Cases  (ii)  $\rho^{q+1}=-1$, $\tau^{q+1}\ne -1$ and $\t\not\in \FF_q$;  (iv) $\rho^{q+1}=\tau^{q+1}=-1$,  $\t, \t^{-1}\rho\not\in \FF_q$.}
 \vskip 3mm
  In  case $(ii)$, $\lg (0,1, \rho )\Delta_0\rg\in \cal H$ and $\lg (\tau, 1, 0)\Delta_0\rg\in \cal N$, where  $\tau \in\FF_{q^2}\setminus\FF_q$.
 So  $W'$ has an  orthogonal base as follows:
 $$v_1=(\xi_1, 1, 0)\Delta_0,\,v_2=\d_2(t\xi_1, 1, \xi_2)\Delta_0,\,v_3=\d_3(t\xi_1, 1, s\xi_2)\Delta_0,$$
where $\xi_1\in\FF_{q^2}\setminus\FF_q$, $\xi_2\in\FF_{q^2}$, $\d_2 =m'+i$ and $\d_3=n'+i$, where $m',n'\in\FF_q$.

  In case  $(iv)$,  $\lg (0,1, \rho )\Delta_0\rg , \lg (\tau, 1, 0)\Delta_0\rg \in \cal H$,  where $\rho\in \FF_q$,  and  $\tau \in\FF_{q^2}\setminus\FF_q$.
Then  $\lg (\tau, 0, \rho)\rg \in \mathbf{w}'\cap \cal N$ and so $W'$ has an  orthogonal base as follows:
 $$v_1=(\xi_1, 0, 1)\Delta_0,\,v_2=\d_2(t\xi_1,  \xi_2, 1)\Delta_0,\,v_3=\d_3(t\xi_1,  s\xi_2, 1)\Delta_0,$$
where $\xi_1=\tau \rho^{-1} \in\FF_{q^2}\setminus\FF_q$, $\xi_2\in\FF_{q^2}$, $\d_2 =m'+\i$ and $\d_3=n'+\i$, where $m',n'\in\FF_q$.

The arguments for both cases are analogous to case $(i)$, and are thus omitted.\qed

\begin{lem}\label{f2}  ${\bf f}_2\leq \frac 52q^2+4q-\frac 12$.
 \end{lem}
 \demo Since $\mathbf{w}$ has two classes of orthogonal frames with representatives $\mathcal{D}_{0,\i}$, where $\i=1, 2$,  we set $F_{2, \i}=\{g_1 \mathcal{D}_{0,\i}\in F_2\di  g_1\in M\}$.
 It suffices to show that for  $|F_{2,\i}|\leq \frac 12 (\frac 52q^2+4q-\frac 12)$ for each $\i=1,2$.
   Then  $\Phi(\mathcal{D}_{0,\i},\mathbf{w}')=\{ \lg (0, 1, \rho )\Del_{0,\i}\rg , \lg (\tau, 1, 0)\Del_{0,\i}\rg , \lg (\tau, 0, -\rho )\Del_{0,\i}\rg \}$ is one of  cases  $(iii)$, $(v)$ and $(vi)$,  by Eq(\ref{F1F2}). For each  case, since $\lg (0,1, \rho )\Delta_{0,\i}\rg\in \cal H\cap\mathbf{w}'$,  and $\lg (\t, 1, 0)\Delta_{0,\i}\rg\in \mathbf{w}'$, there exists a row vector $\delta\in\FF_{q^2}^3$ such that
$$W'=a_{0,\i}\lg (0,1, \rho )\Del_{0,\i},     (\t, 1, 0)\Del_{0,\i},  \d\Del_{0,\i}\rg_{\FF_q} $$
for $\i=1,2$ and some $a_{0,\i}\in\FF_{q^2}^*$, where $\mathbf{w}'=\PG(W')$.
Define the corresponding set:
  $$L_{\mathcal{D}_{0,\i}}(W'):=a_{0,\i}\lg (0,1, \rho ),(\t, 1, 0),\d \rg_{\FF_q}.$$
  Then we shall  measure the cardinality    of the set $\Psi$ of such frames $\mathcal{D}_{0,\i}$.
For any such $\mathcal{D}_1$, there exists $g_1\in M$ such that $\mathcal{D}_1=g_1\mathcal{D}_{0,\i}$ for some $\i=1,2$. Then
$$W'=a_1\lg (0,1, \rho_1)\Del_1,     (\t_1, 1, 0)\Del_1,  \d_1 \Del_1\rg_{\FF_q}\,  {\rm and}\, L_{\mathcal{D}_1}(W'):=a_1\lg (0,1, \rho_1 ),     (\t_1, 1, 0),  \d_1 \rg_{\FF_q},
$$
where $\Del_1$ is a matrix whose row vectors are from the orthogonal frame $\mathcal{D}_1.$
 Then  $L_{\mathcal{D}_{0,\i}}(W')=L_{\mathcal{D}_1}(W')g_1$ for some $\i=1,2$,
that is,
\begin{equation}\label{g1}
a_{0,\i}\lg (0,1, \rho ),     (\t, 1, 0),  \d\rg_{\FF_q} g_1^{-1}=a_1\lg (0,1, \rho_1 ),     (\t_1, 1, 0),  \d_1 \rg_{\FF_q},
\end{equation}
where $a_{0,\i},a_{1}\in\FF_{q^2}^*$, $g_1\in M$ and $\t$ and $\rho$ are given. Clearly,  $a_{0,\i}a_{1}^{-1}\in\FF_q^*,$ as $g_1\in M.$ Now we discuss  three cases, separately.
\vskip 3mm
 {\it Case  (iii)  $\rho^{q+1}=-1, \tau^2\ne -1,  \tau \in \FF_q.$}
 \vskip 3mm
  First, suppose $\rho\in \FF_q$. Set  $\ell_0=a_{0,\i}\lg (0, 1, \rho ),     (\t, 1, 0)\rg_{\FF_q}$  and  $\ell_1=a_1\lg (0,1, \rho_1 ),     (\t_1, 1, 0)\rg_{\FF_q} $, where  $\rho, \rho_1=\pm \sqrt{-1}$ and $\t, \t_1\in \FF_q\setminus \{ 0, \pm \sqrt{-1}\}$.
  Then both of $\ell_0$ and $\ell_1$ contain two isotropic points.
  Clearly,  $g_1^{-1}$ maps two isotropic points of $\ell_0$ to that of $\ell_1$ and   $\lg \d \rg $ to $\lg \d_1\rg $, where $\d\in\ell_0^\perp$ and $\d_1\in\ell_1^\perp$.   Therefore,   $g_1$  has at most two choices, for given
  $(\t_1, \rho_1)$.   Thus we have $2\cdot (q-3)\cdot 2 $ choices for $g_1$ and so $(q-3)$ choices for desired frames, that is
    $|\Psi|\le q-3$, noting that for any $x\in \{  [1, -1, -1], [-1, 1, -1], [-1, -1, 1]\} $,  $(x g_1)\Del_{0,\i}$  and $g_1\Del_{0,\i}$ give same frame.

   Secondly, suppose $\rho\not\in \FF_q$. Then $a_{1}^{-1}L_{\mathcal{D}_{0,\i}}(W')\cap W\in\lg(\t, 1, 0)\rg$ and so  $a_{1}^{-1}L_{\mathcal{D}_1}(W')\cap W\in\lg (\t_1, 1, 0)\rg$ so that
   $\lg(\t, 1, 0)\rg g_1^{-1}=\lg(\t_1, 1, 0)\rg$. Since $M_{\lg (\t, 1, 0)\rg }\cong D_{2(q\pm 1)}$,   $g_1$ has  at most $2(q+1)$  choices for given $\t_1$ so that it has $2(q+1)\cdot \frac{q-1}2$ choices for all possible $\t_1$ (as $\frac{\t^2+1}{\t_1^2+1}\in (\FF_q^*)^2$). Therefore,  $|\Psi|\le \frac{q^2-1}4$, with the same reason as in last paragraph.

\vskip 3mm {\it Case (v): $ \rho^{q+1}=\tau^{2}=-1,  \rho \not\in \FF_q^*, \tau \in \FF_q.$}
\vskip 3mm
In Case $(v)$, $a_{1}^{-1}L_{\mathcal{D}_{0,\i}}(W')\cap W\in\lg(\t, 1, 0)\rg$ where $\t=\pm \sqrt{-1}$ and   then  $a_{1}^{-1}L_{\mathcal{D}_1}(W')\cap W\in\lg (\t_1, 1, 0)\rg$, so that $\lg(\t, 1, 0)\rg g_1^{-1}=\lg(\t_1, 1, 0)\rg$.
Since $M_{\lg (\t, 1, 0)\rg }\cong \ZZ_p^m:\ZZ_{q-1}$,   $g_1$ has  at most $q(q-1)$  choices for given $\t_1$ so that it has $2q(q-1)$ choices for all possible $\t_1$. Therefore,  $|\Psi|\le \frac{q(q-1)}{2}$.
\vskip 3mm
{\it   Case (vi):  $\rho^{q+1}=\tau^{q+1}=-1,  \t^{-1}\rho \in \FF_q.$}
\vskip 3mm For Case (vi), if $\t, \rho\in \FF_q$, then $a_1^{-1}L_{\mathcal{D}_0,\i}(W')\cap W=\ell_0:=x_0\lg (\t, 1, 0),  (0, 1, \rho)\rg_{\FF_q} $ and so   $a_{0,\i}^{-1}L_{\mathcal{D}_1}(W')\cap W=\ell_1:=x_1\lg (\t_1, 1, 0),  (0, 1, \rho)\rg $ for some $x_0,x_1\in\FF_{q^2}^*$, so that $g_1^{-1}$ maps
$\{x_0(\t, 1, 0),  x_0(0, 1, \rho)\}$ to $\{x_1(\t_1, 1, 0),  x_1(0, 1, \rho_1)\}.$ Then $g_1$   has  at most $8$  choices  for all possible $\t_1$ and $\rho_1$. Therefore,  $|\Psi|\leq2$.

 Suppose that  $\t, \rho\not\in \FF_q$. Write $\rho=\pm\t$. Then $ \Phi(\mathcal{D}_{0,\i},\mathbf{w}')\cap\mathbf{w}=\lg(1,0,\pm1)\Delta_{0,\i}\rg\in\mathbf{w}\cap\mathbf{w}'$.
 Since $\Phi(\mathcal{D}_1,\mathbf{w}')=\{ \lg (0, 1, \rho_1 )\Delta_{1}\rg , \lg (\tau_1, 1, 0)\Delta_{1}\rg , \lg (1, 0, \pm1 )\Delta_{1}\rg \}$, where $\Delta_1=g_1\Delta_{0,\i}$,
  we get that $\lg(1,0,\pm1)\Delta_{0,\i}\rg g_1=\lg(1,0,\pm1)g_1\Delta_{0,\i}\rg \in\mathbf{w}\cap\mathbf{w}'\cap\cal N$.
Since $|\mathbf{w}\cap\mathbf{w}'\cap\cal N|\leq q+2$ and $|M_{\lg(1,0,\pm1)\Delta_{0,\i}\rg}|\leq 2(q+1)$, we have that $|\Psi|\leq 2(q+1)(q+2)/4=\frac{(q+1)(q+2)}{2}.$

\vskip 3mm
Totally, we have at most $ 2[(q-3)+\frac{q^2-1}4+\frac{q(q-1)}{2}+2+\frac{(q+1)(q+2)}{2}]=\frac 52q^2+4q-\frac 12$ frames satisfying one of the above three cases. So our ${\bf f}_2$ has this upper bonud.
  \qed

\begin{lem}\label{n2d}
  Let $n_{3}=\frac{1}{3}|\Pi_3|$ be defiend before Eq(\ref{nd}). Then
  $n_3\leq\frac 16q^2+\frac {11}6q\sqrt{q}+\frac{26}{3}q+\frac{20}{3}\sqrt{q},$ for $q>108.$
\end{lem}
\demo Let $h'(a,c,z)=\sum_{\imath=0}^3 a_{\imath}a^{\imath}$, where $a_{\imath}\in\FF_q[c,z]$, with
the discriminant $d(c,z)$ related to $a$.  By Proposition~\ref{disc},   $(c, z)\in \Pi_3$ if and only if  $h'(a,c,z)$ has only one or two roots for $a$,  only if    $d(c,z)=\th u^2$ for some $u\in\FF_q$.
 Set $l(c,z,u)=\th u^2-d(c,z)$ and $\Theta:=\{(c,z)\in\FF_q^2\mid l(c,z,u)=0,u\in\FF_q\}$.
 Then  $|\Pi_3|\leq|\Theta|$,  which will be  bounded below.

 We recall that $i^2=\th$ where, $\FF_q^*=\lg \th\rg$ and $i\in\FF_{q^2}$.
Suppose that for  some $z\in \FF_q$,  $l(c,z,u)$ is reducible over a field $K\ge \FF_q$.  Then  $\th u^2-d(c,z)=l(c,z,u)=\th(u-l_0(c))(u-l_0'(c))$ over $K$. Consequently, $l_0(c)=-l'_0(c)$, and thus $\th^{-1}d(c,z)=l_0(c)^2$.
  Write $l_0(c)=\sum_{0\leq \jmath\leq 4}k_{\jmath}c^{\jmath}$, as $\deg(d(c,z), c)=8$. Check that  $l_0(c)^2\in\FF_q[c]$ if and only if  $k_{\jmath}\in\FF_q$ or $k_{\jmath}\in i\FF_q$ for all $0\leq \jmath\leq 4$. Therefore, either $\th u^2=d(c,z)=l_1(c)^2$ or   $\th u^2=d(c,z)=\th l_2(c)^2$, that is $u^2=l_2(c)^2$,
  where  $l_1(c), l_2(c)\in \FF_q[c]$.
 Reasonably,  we  introduce the following sets:
 $$\begin{array}{lll}
 \beth_0&=&\{ z\in \FF_q\mid d(c,z)\equiv 0 \, {\it in }\, \FF_q[c]\}; \\
  \beth_1&=&\{ z\in \FF_q\mid  d(c,z)=l_1(c)^2\not\equiv 0 \,  {\it in }\, \FF_q[c]\}; \\
\beth_2&=&\{ z\in \FF_q\mid    d(c,z)=\th l_2(c)^2\not\equiv 0 \,  {\it in }\, \FF_q[c]\}; \\
\beth_3&=&\{ z\in \FF_q\mid   \th u^2-d(c,z)\, {\rm is\, absolutely\, irreducible\, in\, } \FF_q[c,u],  d(c,z)\not\equiv 0\,  {\it in }\, \FF_q[c]\};\\
  \Theta_\i&=&\{(c,z)\in \Theta \mid z\in \beth_\i\}, \i=0, 1, 2, 3.\end{array}$$
  \vskip -5pt
\f Then $\beth_\i\cap \beth_\j=\emptyset$ for any $\i\neq\j$ and consequently $\FF_q=\cup _{\i=0}^3 \beth_\i$, $\Theta=\bigcup_{i=0}^3 \Theta_i$.

By Lemma~\ref{D4Z2},  there exist at least  $\frac 12(q-5\sqrt{q})-1$   pairs $(c, z)\in\FF_q^2$ such that $h'(a,c,z)=0$ has three distinct solutions $a$.  By Proposition~\ref{disc},  for each such pair,  $d(c,z)$ is  a nonzero square element in  $\FF_q^*$.  Then  $0\lvertneqq|\beth_2 |\lvertneqq q$,    provided  $\frac 12(q-5\sqrt{q})-1\gvertneqq 0$, that is $q\ge 37$.

Let $z\in \beth_0$. Then  $d(c,z)\equiv 0$ in $\FF_q[c]$.   Set $d(c,z):=\sum_{\i=0}^8f_\i(z)c^\i$. Then $f_\i(z)=0$ for each $\i$. Since $\deg(f_\i(z), z)\le 12$, we get  $|\beth_0|\le 12$ and so   $|\Theta_0|\le q |\beth_0|$.

Let $z\in \beth_1$. Then $\th u^2=l_1(c)^2$, and thus $u=l_1(c)=0$.  Since $\deg(l_1(c))\le 4$ and $l_1(c)\not\equiv 0$, we have at most 4 solutions $c$. So $|\Theta_1|\le 4|\beth_1|$.

 \vskip 3mm
 {\it Step  1:}  Derive  a general statement:
\vskip 3mm
{\it Statement 1:  Suppose that \vspace{-5pt}$$(u_4(z)c^4+\frac{b_3(z)}{a_3(z)}c^3+
\frac{b_2(z)}{a_2(z)}c^2+\frac{b_1(z)}{a_1(z)}c
+\frac{b_0(z)}{a_0(z)})^2=\sum_{\j=0}^8w_\j(z)c^\j,$$   where $u_4(z), a_\i(z), b_\i(z), w_\j(z)\in \FF_q[z]$ but $a_\i(z)\ne 0$,  for
$0\le \i\le 3$, $0\le \j\le 8$;  $z\in \FF_q\setminus S$, where $S=\cup_{\i=0}^3 S_\i$ and  $S_\i=\{z\in \FF_q\mid  a_\i(z)=0\}$. Then $a_\i(z)\mid b_\i(z)$ where $0\le \i\le 3$,
provided  $q\gvertneqq  6r$, where $r=\max\{\deg(u_4(z)), \deg(a_\i(z)), \deg(b_\i(z)), \deg(w_\j(z))\mid 0\leq\i\leq3, 0\leq\j\leq8\}.$}

\vskip 3mm
{\it Proof of  Statement 1:} Note that if $\frac{f(z)}{g(z)}=h(z)$ holds for any $z\in \FF_q\setminus \Upsilon$, where $\Upsilon=\{z\in \FF_q\mid  g(z)=0\}$ and
 $q\gneqq \deg(g(z))+\max\{ \deg(f(z)), \deg(g(z))+\deg(h(z))\}$,  for some $h(z)$, then $g(z)\mid f(z)$.

 No loss, assume  $(b_\i(z), a_\i(z))=1$, where $\i=0,1, 2, 3$.  Since
\vspace{-5pt}$$w_6(z)=(\frac{b_3(z)}{a_3(z)})^2
+2u_4(z)\frac{b_2(z)}{a_2(z)}
=\frac{a_2(z)b_3(z)^2+2u_4(z)a_3(z)^2b_2(z)}{a_2(z)a_3(z)^2}\in \FF_q[z],$$ by the above argument,
we get $a_3(z)^2\mid a_2(z)$. Similarly, from $w_0(z)\in\FF_q[z]$, $w_2(z)\in \FF_q(z)$ and $w_4(z)\in \FF_q(z)$, we get $a_0(z)\in\FF_q^*$, $a_1(z)^2\mid a_2(z)$ and $a_2(z)^2\mid a_1(z)a_3(z)$, respectively.
Therefore,  $a_3(z)^2a_1(z)^2\mid a_1(z)a_3(z)$. Therefore, $a_1(z), a_3(z)\in \FF_q^*$, forcing $a_2(z)\in \FF_q^*$.

\vskip 3mm
{\it Step 2:}\,  Measure   $|\Theta|$.
 \vskip 3mm
To measure $|\beth_2|$,  write  $ u^2=\th\frac{d(c,z)}{\th ^{2}}=\sum_{\j=0}^8 w_\j(z)c^\j=l_2(c)^2$, where $w_\j(z)\in\FF_q[z]$ and $l_2(c)=\sum_{\i=0}^4u_\i c^\i\in \FF_q[c]$. This yields that
\vspace{-5pt}\begin{equation}\label{w8}
\sum_{\j=0}^8 w_\j(z)c^\j=(\sum_{\i=0}^4u_\i c^\i)^2.\end{equation}
Comparing coefficients, we obtain the following system of equations:
\vspace{-5pt}\begin{equation} \label{system}
  \begin{array}{llll}
&u_4^2=w_8(z), &2u_4u_3=w_7(z),\,  &u_3^2+2u_4u_2=w_6(z),\\
&2u_4u_1+2u_3u_2=w_5(z),  &u_2^2+2u_4u_0+2u_3u_1=w_4(z), &2u_3u_0+2u_2u_1=w_3(z),\\
  &u_1^2+2u_2u_0=w_2(z), &2u_1u_0=w_1(z),  &u_0^2=w_0(z).
  \end{array}\end{equation}
Since  $|\beth_2|\lvertneqq q$,  some of equations in Eq(\ref{w8})   does not hold for some values of $z$. No loss, we assume $w_8(z)\not\equiv 0$, while  we have the same arguments for other cases.
Let $\beth'=\{ z\in K\di   w_8(z)=0\}$.  Then $|\beth'|\leq 12.$

For $z\in \FF_q\setminus \beth'$, by solving the  first five equations  of Eq(\ref{system}) over $\FF_q$  we get
\vspace{-5pt}$$
\begin{array}{ll} &u_4^2=w_8(z), \, u_3=\frac{w_7(z)}{2u_4},\, u_2=\frac{4w_6(z)u_4^2-w_7(z)^2}{8u_4^3},\,
 u_1=\frac{8u_4^4w_5(z)-4w_7(z)w_6(z)u_4^2+w_7(z)^3}{16u_4^5},\\
&u_0=\frac{64u_4^6w_4(z)-16w_6(z)^2u_4^4-5w_7(z)^4+24u_4^2w_6(z)w_7(z)^2-32u_4^4w_5(z)w_7(z)}{128u_4^7}
.\end{array}
$$
 Inserting $u_3$, $u_2$, $u_1$ and $u_0$  to the remaining equations of  Eq(\ref{system})   and using the relation $u_4^2=w_8(z)$, we get four equations about $z$, say $F_\k(z)=0$ where $1\le \k\le 4$, where $F_\k(z)\in \FF_q[z]$ with  $\deg(F_\k(z))\le 96$. We then distinguish the following two cases.
 \vskip 2mm
{\it Case 1:}  Suppose that $F_\i(z)\not\equiv 0$ for some $\i$.
 \vskip 2mm
In this case,  our $z$ has at most 96 possibilities in Eq(\ref{system})  (and so  Eq(\ref{w8}), while some  of $z$  in $\beth'$ may satisfy   Eq(\ref{w8}) too.   Therefore, $|\beth_2|\le 96+12=108$.
For any $z\in \beth_2$, $c$ is free and so $|\Theta_2|=q|\beth_2|$.

For any $z\in \beth_3$,  let $$k_1(z)=|\{ c\in \FF_q\di  (c, z)\in \Theta_3\}| \quad {\rm and}\quad k_1'(z)=|\{ c\in \FF_q\di d(c,z)=0\}|.$$
Then $k_1'(z)\le 8$.   For any $z\in \beth_3$,  since $l(c, z, u)\in \FF_q[c, u]$ is absolutely irreducible,
      by Proposition~\ref{weil} and $\deg(d(c,z), c)\le 8$,  there are at most the number $\d'=q+7\sqrt{q}$ of pairs $(c, u)$ such that $l(c,z,u)=0$  for any $z\in\beth_3$.
  Therefore, $k_1(z)\le \frac 12 (\d'-k_1'(z))+k_1'(z)=\frac 12 \d'+4$.  Note that  $|\beth_3|=q-\sum_{\i=0}^2 |\beth_\i|$.
  Therefore, for all $q> 64$, we get that
 \vspace{-5pt}$$\begin{array}{lll}
|\Theta|&=&\sum_{\i=0}^3 |\Theta_\i|\le  q|\beth_0|+4|\beth_1|+q|\beth_2|+(q-\sum_{\i=0}^2 |\beth_\i|)(\frac 12 \d'+4)\\
&\le & q(\frac 12\d'+4)+(q-\frac 12 \d'-4)(|\beth_0|+|\beth_2|)\le   \frac 12q^2+\frac 72q\sqrt{q}+64q.\\
\end{array}$$ \vspace{-5pt}
{\it Case 2:} Suppose that    $F_\k(z)=0$ for any $z\in \FF_q\setminus \beth'$ with $1\le \k\le 4$.
 \vskip 3mm
In this case, Eq(\ref{w8}) holds  with $u_\i\in K$, for all  $z\in \FF_q\setminus \beth'$, where $K$ is a finite field containing $\FF_q$. This implies that for all such $z$, $u^2-\sum_{\i=0}^8 w_\i(z)c^\i$ is reducible
 over some finite field $K$. So $|\beth_3|\le |\beth'|\le 12 $ and $|\beth_1|+|\beth_2|\ge q-24$. Following the arguments in second paragraph of our proof,  we have the following three subcases:
\vskip 2mm
(2.1)  $w_8(z)=l'(z)^2$ for some $l'(z)\in \FF_q[z]$:  Then   $u_4=u_4(z)$ is a polynomial. By Statement 1, and since the degree of $z$ in  Eq(\ref{system}) is at most $96$ while $|\beth'|\leq 12$, it follows that for $q>108$, all $u_\i$ $(0\le \i\le 4)$ are polynomials.
Consequently, $|\beth_2|=q$, which leads to a contradiction.
  \vskip 2mm

(2.2)   $w_8(z)=\th l'(z)^2$ for some $l'(z)\in \FF_q[z]$:  write
 \vspace{-5pt}$$u^2= \sum_{\j=0}^8 w_\j(z)c^\j=(\sum_{\i=0}^4u_\i c^\i)^2=\frac 1{w_8(z)} (w_8(z)c^4+\sum_{\i=0}^3\frac{f_\i(z)}{g_\i(z)}c^\i)^2,$$
 where $f_\i(z), g_\i(z)\in \FF_q[z]$ and $z\in \FF_q\setminus \beth'$, a contradiction.
 \vskip 2mm
  (2.3)\, Suppose $f(u_4, z)=u_4^2-w_8(z)$  is absolutely  irreducible.
 By Proosition~\ref{weil},   we have  at most $\d''=q+11\sqrt{q}$ solutions $(z, u_4)$ for $f(u_4,z)=0$, as $\deg(w_8(z))\le 12$.
 Let $k_2'(z)=|\{ z\in \FF_q\mid w_8(z)=0\}|.$
  Then $k_2'(z)\le 12$, and consequently $|\beth_2|\le \frac 12 (\d''-k_2'(z))+k_2'(z)\le \frac 12 \d''+6$.  Therefore,
\vspace{-5pt}$$\begin{array}{lll}
|\Theta|&=&\sum_{\i=0}^3 |\Theta_\i|\le  q|\beth_0|+4|\beth_1|+q|\beth_2|+|\beth_3|(\frac 12 \d'+4)\\
&\le &  q|\beth_0|+4(q-|\beth_0|-|\beth_2|-|\beth_3|)+q|\beth_2|+|\beth_3|(\frac 12 \d'+4)\\
&\le &  (q-4)|\beth_2|+(q-4)|\beth_0|+|\beth_3|(\frac 12 \d')+4q\\
&\le &\frac 12q^2+\frac {11}2q\sqrt{q}+26q+20\sqrt{q}.
\end{array}$$

Finally, by taking the maximum possible value for $|\Theta|$, we obtain  $n_3\le \frac 13|\Theta|\le \frac 16q^2+\frac {11}6q\sqrt{q}+\frac{26}{3}q+\frac{20}{3}\sqrt{q}.$
\qed

\begin{lem}\label{n2}  Let $n_{1}=|\Pi_1|$. Then
 $n_1\le \frac{q^2}{2}+\frac{3}{2}q\sqrt{q}+14q-36\sqrt{q},$ for $q\geq80.$
\end{lem}
\demo   Let $d=1$ and $n_{1}=|\Pi_1|$  be defined in Lemma~\ref{L_1(g)3}.  Let $h(a,b,c)=\sum_{\imath=0}^3 a_{\imath}a^{\imath}$, where $a_{\imath}\in\FF_q[b,c]$, with
the discriminant $d(b,c)$, corresponding  to $a$.  By Proposition~\ref{disc}, $(b, c)\in \Pi_1$ if and only if  $h(a,b,c)$ has only one or two roots for $a$,  only if $d(b,c)=\th u^2$ for some $u\in\FF_q$.
 Set $l(b,c,u)=\th u^2-d(b,c)$ and $\Theta:=\{(b,c)\in\FF_q^2\mid l(b,c,u)=0,u\in\FF_q\}$.
 Then  $|\Pi_1|\leq|\Theta|$. In what follows,  we shall get an upper bound for $|\Theta|$.

Let  $i^2=\th$ and $\FF_q^*=\lg \th\rg$.
Suppose that some $b\in \FF_q$,  $l(b,c,u)$ is reducible over a field $K\ge \FF_q$.  Then  $\th u^2-d(b,c)=l(b,c,u)=(iu-l_1(c))(iu-l_1'(c))$ over $K$.   Then  $l_1(c)=-\l'_1(c)$ so that  $d(b,c)=l_1(c)^2$.
  Write $l_1(c)=\sum_{0\leq \jmath\leq 2}k_{\jmath}c^{\jmath}$, as $\deg(d(b,c), c)=4$. Check that  $l_1(c)^2\in\FF_q[c]$ if and only if  $k_{\jmath}\in\FF_q$ or $k_{\jmath}\in i\FF_q$ for all $0\leq \jmath\leq 2$. Then either $\th u^2=d(b,c)=l_1(c)^2$ or   $\th u^2=d(b,c)=(i l_2(c))^2=\th l_2(c)^2$, that is $u^2=l_2(c)^2$,
  where  $l_1(c), l_2(c)\in \FF_q[c]$.
 Therefore, we are now introducing the following sets:
 $$\begin{array}{lll} \beth_0&=&\{ b\in \FF_q\mid  d(b,c)\equiv 0 \, {\it in }\, \FF_q[c]\}; \\
  \beth_1&=&\{ b\in \FF_q\mid  d(b,c)=l_1(c)^2\not\equiv 0 \,  {\it in }\, \FF_q[c]\}; \\
\beth_2&=&\{ b\in \FF_q\mid    d(b,c)=\th l_2(c)^2\not\equiv 0 \,  {\it in }\, \FF_q[c]\}; \\
\beth_3&=&\{ b\in \FF_q\mid   \th u^2-d(b,c)\, {\rm is\, absolutely\, irreducible\, in\, } \FF_q[c,u],  d(b,c)\not\equiv 0\,  {\it in }\, \FF_q[c]\};\\
  \Theta_\i&=&\{(b,c)\in \Theta \mid b\in \beth_\i\}, \i=0, 1, 2, 3.\end{array}$$
Then $\beth_\i\cap \beth_\j=\emptyset$ for any $\i\neq\j$ and so $\FF_q=\cup _{\i=0}^3 \beth_\i$ and $\Theta=\bigcup_{i=0}^3 \Theta_i$.

By Lemma~\ref{D4Z2},  there exist at least  $\frac{(q-5\sqrt{q})}{2d}-\frac{1}{3}$   pairs $(b, c)\in\FF_q^2$ such that $h(a,b,c)=0$ has three distinct solutions $a$.  By Proposition~\ref{disc},  for each such pair,  $d(b,c)$ is  a nonzero square element in  $\FF_q^*$.  Then  $|\beth_2 |\lvertneqq q$,    provided  $\frac 12(q-5\sqrt{q})-1\gvertneqq 0$, that is $q\ge 37$.

Let $b\in \beth_0$. Then  $d(b,c)\equiv 0$ in $\FF_q[c]$.   Set $d(b,c):=\sum_{\i=0}^4f_\i(b)c^\i$. Then $f_\i(b)=0$ for each $\i$. Since $\deg(f_\i(b), b)\le 4$, we get  $|\beth_0|\le 4$ and so   $|\Theta_0|= q |\beth_0|$.

Let $b\in \beth_1$. Then $\th u^2=l_1(c)^2$. Then $u=l_1(c)=0$.  Since $\deg(l_1(c))\le 2$ and $l_1(c)\not\equiv 0$, we have at most 2 solutions $c$. So $|\Theta_1|\le 2|\beth_1|$.

\vskip 3mm
{\it Step 1:}\,  Measure   $|\Theta|$.
 \vskip 3mm
To measure $|\beth_2|$,  write  $ u^2=\th\frac{d(b,c)}{\th ^{2}}=\sum_{\j=0}^2 w_\j(b)c^\j=l_2(c)^2$, where $w_\j(b)\in\FF_q[b]$ and $l_2(c)=\sum_{\i=0}^2u_\i c^\i\in \FF_q[c]$, so that
\vspace{-5pt}\begin{equation}\label{w8'}\sum_{\j=0}^4 w_\j(b)c^\j=(\sum_{\i=0}^2u_\i c^\i)^2.\end{equation}
Then  we get the system
\begin{equation} \label{system'}
\begin{array}{llllll}
   & u_2^2=w_4(b)
   & 2u_1u_2=w_3(b)
&u_1^2+2u_0u_2=w_2(b)
&2u_0u_1=w_1(b)
& u_0^2=w_0(b).
  \end{array}
\end{equation}
Since  $|\beth_2|\lvertneqq q$,  some of equations in the system  does not hold for some values of $b$. No loss, we assume $w_4(b)\not\equiv 0$, while  we have the same arguments for other cases.
Let $\beth'=\{ b\in K\mid   w_4(b)=0\}$.  Then $|\beth'|\leq 4.$

For $b\in \FF_q\setminus \beth'$, by solving the above first five equations (from left to right) in the system Eq(\ref{system'}) over $\FF_q$  we get
$$
\begin{array}{lllll}
&u_2^2=w_4(b),
&u_1=\frac{w_3(b)}{2u_2},
&u_0=\frac{4w_2(b)w_4(b)-w_3(b)^2}{8u_2^3}.
\end{array}
$$
 Inserting $u_1$ and $u_0$  to the remaining equations of the system  and using the relation $u_2^2=w_4(b)$, we get two equations about $b$, say $F_\k(b)=0$ where $1\le \k\le 2$, and $F_\k(b)\in \FF_q[b]$ with  $\deg(F_\k(b))\le 16$. Then we have the following two cases.
 \vskip 3mm
{\it Case 1:}  Suppose that $F_\i(b)\not\equiv 0$ for some $\i$.
 \vskip 3mm
In this case,  our $b$ has at most 16 possibilities in Eq(\ref{system'})  (and so  Eq(\ref{w8'}), while some  of $b$  in $\beth'$ may satisfy   Eq(\ref{w8'}) too.   Therefore, $|\beth_2|\le 16+4=20$.
For any $b\in \beth_2$, $c$ is free and so $|\Theta_2|=q|\beth_2|$.

For any $b\in \Del_3$,  let $$k_1(b)=|\{ c\in \FF_q\di  (b,c)\in \Theta_3\}|, \quad k_1'(b)=|\{ c\in \FF_q\di d(b,c)=0\}|.$$
Then $k_1'(b)\le 4$.   For any $b\in \beth_3$,  since $l(b,c, u)\in \FF_q[c, u]$ is absolutely irreducible,
      by Proposition~\ref{weil} and $\deg(d(b,c), c)\le 4$,  there are at most the number $\d'=q+3\sqrt{q}$ of pairs $(c, u)$ such that $l(b,c,u)=0$  for any $b\in\beth_3$.
  Then $k_1(b)\le \frac 12 (\d'-k_1'(b))+k_1'(b)=\frac 12 \d'+2$.  Note that  $|\beth_3|=q-\sum_{\i=0}^2 |\beth_\i|$.
  Therefore, for all $q> 64$, we get that
$$\begin{array}{lll}
|\Theta|&=&\sum_{\i=0}^3 |\Theta_\i|\le  q|\beth_0|+2|\beth_1|+q|\beth_2|+(q-\sum_{\i=0}^2 |\beth_\i|)(\frac 12 \d'+2)\\
&\le & q(\frac 12\d'+2)+(q-\frac 12 \d'-2)(|\beth_0|+|\beth_2|)\le   \frac 12q^2+\frac 32q\sqrt{q}+14q-36\sqrt{q}-48.\\
\end{array}$$
{\it Case 2:} Suppose that $F_\k(b)=0$ for any $b\in \FF_q\setminus \beth'$ with $1\le \k\le 2$.
 \vskip 3mm
In this case, Eq(\ref{w8'}) holds  with $u_\i\in K$, for all  $b\in \FF_q\setminus \beth'$, where $K$ is a finite field containing $\FF_q$. This implies that for all such $b$, $u^2-\sum_{\i=0}^4 w_\i(b)c^\i$ is reducible
 over some finite field $K$. So $|\beth_3|\le |\beth'|\le 4 $ and $|\beth_1|+|\beth_2|\ge q-8$. Following the arguments in second paragraph of our proof,  we have the following three subcases:
\vskip 3mm
(2.1)  $w_4(b)=l'(b)^2$ for some $l'(b)\in \FF_q[b]$:  Then   $u_2=u_2(b)$ is a polynomial. Then by the above Statement 1,
 all $u_\i$ where $0\le \i\le 2$ are polynomials provided $q>20$, as the degree of $b$ in the system Eq(\ref{system'}) is at most $16$ and $|\beth'|\leq 4$. Then
 $|\beth_2|=q$,  a contradiction.
  \vskip 3mm

(2.2)   $w_4(b)=\th l'(b)^2$ for some $l'(b)\in \FF_q[b]$:  write
 \vspace{-5pt}$$u^2= \sum_{\j=0}^4 w_\j(z)c^\j=(\sum_{\i=0}^2u_\i c^\i)^2=\frac 1{w_4(b)} (w_4(b)c^2+\sum_{\i=0}^1\frac{f_\i(b)}{g_\i(b)}c^\i)^2,$$
 where $f_\i(b), g_\i(b)\in \FF_q[z]$ and $b\in \FF_q\setminus \beth'$, a contradiction.
 \vskip 3mm
  (2.3)\, Suppose $f(u_2, b)=u_2^2-w_4(b)$  is absolutely  irreducible.
 By Proosition~\ref{weil},   we have  at most $\d''=q+3\sqrt{q}$ solutions $(b, u_2)$ for $f(u_2,b)$, as $\deg(w_4(b))\le 4$.
 Let $k_2'(b)=|\{ b\in \FF_q\mid w_4(b)=0\}|.$
  Then $k_2'(b)\le 4$.
  Then $|\beth_2|\le \frac 12 (\d''-k_2'(b))+k_2'(b)\le \frac 12 \d''+2$.  Therefore,

$$\begin{array}{lll}
|\Theta|&=&\sum_{\i=0}^3 |\Theta_i|\le  q|\beth_0|+2|\beth_1|+q|\beth_2|+|\beth_3|(\frac 12 \d'+2)\\
&\le &  q|\beth_0|+2(q-|\beth_0|-|\beth_2|-|\beth_3|)+q|\beth_2|+|\beth_3|(\frac 12 \d'+2)\\
&\le &  (q-2)|\beth_2|+(q-2)|\beth_0|+|\beth_3|(\frac 12 \d')+2q\\
&\le &\frac 12q^2+\frac {3}2q\sqrt{q}+3\sqrt{q}+9q-12.
\end{array}$$

 \f Take the maximal one for $|\Theta|$, we  get $n_1\le \frac 12q^2+\frac 32q\sqrt{q}+14q-36\sqrt{q},$ for $q\geq80$
\qed

\subsubsection{ Upper bounds of $\hm B_3$ and $\hm B_4$}
\begin{lem}\label{C3}
${\hm B_3}\leq q(q^2-1)(q+1)$ and ${\hm B_4}\le \frac{2}{d}q(q+1)(q+2)(q+1-d).$
\end{lem}
\demo (1)  Note that $
{\hm B_3}=(q^2-1)n_{q-1,q-1,1},$ where
\vspace{-5pt}$$ n_{q-1,q-1,1}=|\{M_\imath\in {\cal M}\setminus \{M, M'\}\mid M_\imath\cap M \cong  M_\imath\cap M'\cong \D_{2(q-1)}, M_\imath\cap M\cap M'=1\}|.$$
Let $g=g_A+g_Bi\in \{ {\hm g_1}, {\hm g_2}\}$   be given in Eq(\ref{g1g2}). Using the notation of Eq(\ref{bc}), consider the matrix $ A_{(a', b^*, c^*)}$ given by
$$ A_{(a', b^*, c^*)}=\left\{
               \begin{array}{ll}
                 g_AE_B+g_BE_A, & \hbox{$a'=1$;} \\
                 a(g_AE_B+g_BE_A)+(g_AE_A+g_BE_B\th), & \hbox{$a'=a+i$.}
               \end{array}
             \right.$$
Now, let
\vspace{-5pt} $$E_{A,1}+E_{B,1}i=[1, -b + i, 1],\, E_{A,2}+E_{B,2}i=[1,1, -c+i],\, E_{A,3}+E_{B,3}i=[1, -b + i, -b+i].$$
Correspondingly,  for each $\j=1,2,3$ and each $a'\in \{ 1, a+i\}$,  we define the associated set
$\mathcal{W}_{\j,a'}':=\{(b^*,c^*)\mid \rank(A_{(a', b^*, c^*)})=1\},$ where in the definition of $A_{(a', b^*, c^*)}$, we take $E_A+E_Bi=E_{A,\j}+E_{B,\j}i$. Since there are $\frac{q(q+1)}{2}$ subgroups in $M$ isomorphic to $\D_{2(q-1)}$, and following the same arguments as in Step 2 of the proof of Lemma \ref{hlD}, we obtain
 $$n_{q-1,q-1,1}\leq \frac{q(q+1)}{2} \sum_{j\in \{1, 2, 3\}, a'\in \{1, a+i\}} |\mathcal{W}_{j,a'}'|.$$

Now, we will compute the rank of the matrix $A_{(a',b^*,c^*)}$ in the following.
\vskip 3mm
{\bf case 1: get $\rank(g_AE_{B,1}+g_BE_{A,1})$ for $g_A+g_Bi\in\{g_1,g_2\}$ and $E_{A,1}+E_{B,1}i=[1, -b + i, 1]$}
\vskip 3mm
For this case, then $${\small g_AE_{B,1}+g_BE_{A,1}\in\{\left(
  \begin{array}{ccc}
    0 & 1 & n_1 \\
    1&  t & 0\\
    0& -b & n_2 \\
  \end{array}
\right),\left(
  \begin{array}{ccc}
    0 & 1 & n_1 \\
    0&  1 & 0\\
    1& -bs & n_2 \\
  \end{array}
\right)\}}$$ which are not rank $1$.
So $|\cal W_{1,1}'|=0$.

\vskip 3mm
{\bf case 2: get $\rank(g_AE_{B,2}+g_BE_{A,2})$ for $g_A+g_Bi\in\{g_1,g_2\}$ and $E_{A,2}+E_{B,2}i=[1, 1, -c+i]$}
\vskip 3mm
For this case, then $${\small g_AE_{0,2}+g_BE_{d,2}\in\{\left(
  \begin{array}{ccc}
    0 & 0 & m_1-cn_1 \\
    1& 0 & 0\\
    0& 1 & m_2-cn_2 \\
  \end{array}
\right), \left(
  \begin{array}{ccc}
    0 & 0 & m_1-cn_1 \\
    0& 0 & 0\\
    1& s & m_2-cn_2 \\
  \end{array}
\right)\}.}$$
So for $g=g_1$, we have $|\cal W_{2,1}'|=0$; for $g=g_2$, we have $\rank(g_AE_{B,2}+g_BE_{A,2})=1$ only for $c=\frac{m_1}{n_1}(n_1\neq 0)$, that is, $|\cal W_{2,1}'|\leq1$;
\vskip 3mm

{\bf case 3: get $\rank(g_AE_{B,3}+g_BE_{A,3})$ for $g_A+g_Bi\in\{g_1,g_2\}$ and $E_{A,3}+E_{B,3}i=[1, -b+i, -b+i]$}
\vskip 3mm
For this case, then $${\small g_AE_{0,3}+g_BE_{d,3}\in\{\left(
  \begin{array}{ccc}
    0 & 1 & m_1-bn_1 \\
    1& t & 0\\
    0& -b & m_2-bn_2 \\
  \end{array}
\right), \left(
  \begin{array}{ccc}
    0 & 1 & m_1-bn_1 \\
    0& -1 & 0\\
    1& -bs & m_2-bn_2 \\
  \end{array}
\right)\}.}$$
So for $g=g_1$, we have $|\cal W_{3,1}'|=0$; for $g=g_2$, we also have $|\cal W_{3,1}'|=0$;
\vskip 3mm

{\bf case 4: get $\rank(a(g_AE_{B,j}+g_BE_{A,j})+(g_AE_{A,j}+g_BE_{B,j}\th))$ for $g_A+g_Bi\in\{g_1,g_2\}$ and $E_{A,1}+E_{B,1}i=[1, -b+i, 1]$}
\vskip 3mm
For this case, then $${\small\begin{array}{lll}
&a(g_AE_{B,j}+g_BE_{A,j})+(g_AE_{A,j}+g_BE_{B,j}\th)&\\
&\in\{\left(
  \begin{array}{ccc}
    0 & a-b & m_1+an_1 \\
    e+a& -t(a-b) & 0\\
    -t& \th-ab & m_2+an_2 \\
  \end{array}
\right), \left(
  \begin{array}{ccc}
    0 & a-b & m_1+an_1 \\
    s& -(a-b) & 0\\
    a& (\th-ab)s & m_2+an_2 \\
  \end{array}
\right)\}.&
\end{array}}$$
So for $g=g_1$, we have $\rank(a(g_AE_{B,j}+g_BE_{A,j})+(g_AE_{A,j}+g_BE_{B,j}\th))=1$ only for $b=-e$, that is, $|\cal W_{1,a+i}'|\leq 1$; for $g=g_2$, we have $\rank(g_AE_{B,1}+g_BE_{A,1})\neq1$, that is, $|\cal W_{1,a+i}'|=0$.
\vskip 3mm
{\bf case 5: get $\rank(a(g_AE_{B,2}+g_BE_{A,2})+(g_AE_{A,2}+g_BE_{B,2}\th))$ for $g_A+g_Bi\in\{g_1,g_2\}$ and $E_{A,2}+E_{B,2}i=[1, 1, -c+i]$}
\vskip 3mm
For this case, then $${\small\begin{array}{lll}
&a(g_AE_{B,2}+g_BE_{A,2})+(g_AE_{A,2}+g_BE_{B,2}\th)&\\
&\in\{\left(
  \begin{array}{ccc}
    0 & 1 & (a-c)m_1+(\th-ac)n_1 \\
    e+a& -t & 0\\
    -t& a & (a-c)m_2+(\th-ac)n_2 \\
  \end{array}
\right), \left(
  \begin{array}{ccc}
    0 & 1 & (a-c)m_1+(\th-ac)n_1 \\
    s& -1 & 0\\
    a& as & (a-c)m_2+(\th-ac)n_2 \\
  \end{array}
\right)\}.&
\end{array}}$$
So for $g=g_1$, we have $\rank(a(g_AE_{B,j}+g_BE_{A,j})+(g_AE_{A,j}+g_BE_{B,j}\th))\neq1$ as $t\neq 0$; for $g=g_2$, we have $\rank(a(g_AE_{B,j}+g_BE_{A,j})+(g_AE_{A,j}+g_BE_{B,j}\th))\neq1$, as $s\neq 0$. So we have that $|\cal W_{2,a+i}'|=0$ for $g\in\{g_1,g_2\}$.

\vskip 3mm
{\bf case 6: get $\rank(a(g_AE_{B,3}+g_BE_{A,3})+(g_AE_{A,3}+g_BE_{B,3}\th))$ for $g_A+g_Bi\in\{g_1,g_2\}$ and $E_{A,3}+E_{B,3}i=[1, -b+i, -b+i]$}
\vskip 3mm
For this case, then
$$
{\small\begin{array}{lll}
&a(g_AE_{B,3}+g_BE_{A,3})+(g_AE_{A,3}+g_BE_{B,3}\th)&\\
&\in\{\left(
  \begin{array}{ccc}
    0 & a-b & (a-b)m_1+(\th-ab)n_1 \\
    e+a& (a-b)t & 0\\
    t& \th-ab & (a-b)m_2+(\th-ab)n_2 \\
  \end{array}
\right), \left(
  \begin{array}{ccc}
    0 & a-b & (a-b)m_1+(\th-ab)n_1 \\
    s& b-a & 0\\
    a& s\th-abs & (a-b)m_2+(\th-ab)n_2 \\
  \end{array}
\right)\}.&
\end{array}}
$$
So for $g=g_1$, we have $\rank(a(g_AE_{B,j}+g_BE_{A,j})+(g_AE_{A,j}+g_BE_{B,j}\th))=1$ only for $a=b=-e$; for $g=g_2$, we have $\rank(a(g_AE_{B,j}+g_BE_{A,j})+(g_AE_{A,j}+g_BE_{B,j}\th))=1$ only for $a=b=-e$, as $s\neq 0$. So we have that $|\cal W_{3,a+i}'|=1$ for $g\in\{g_1,g_2\}$.

Then check the rank of the matrix $A_{(a', b^*, c^*)}$ for both ${\hm g_1}$ and ${\hm g_2},$ which is given above, we find that
 $|\cal W_{1,1}'|=|\cal W_{3,1}'|=|\cal W_{2,a+i}'|=0$ and
  $|\cal W_{3,a+i}'|\leq1;$ moreover, $|\cal W_{2,1}'|=0$ and $|\cal W_{1,a+i}'|\leq 1$ for $g={\hm g_1}$, while  $|\cal W_{2,1}'|\leq 1$ and $|\cal W_{1,a+i}'|=0$ for $g={\hm g_2}$. Therefore, we get $n_{q-1,q-1,1}\leq \frac{q(q+1)}{2}\cdot 2=q^2+q$, and consequently, ${\hm B_3}\le (q^2-1)(q^2+q)$.

\vskip 2mm
(2)  By Table \ref{number}, there are at most $\frac{q(q+1)}{2}(\frac{q+1}{d}-1)$  subgroups $M_\imath$ such that $M\cap M_\imath=D_{2(q\pm 1)}$ and so
$n_{q\pm 1, j,1}, n_{j, q\pm 1,1}\le  \frac{q(q+1)}{2}(\frac{q+1}{d}-1)$ for any possible $j$. Thus
$$\begin{array}{lll}
{\hm B_4}&=&2[(q+2)n_{3,q+1,1}+qn_{3,q-1,1}
  +n_{q+1,3,1}+n_{q-1,3,1}]\\
  &\leq& 2(q+2+q+1+1)\frac{q(q+1)}{2}(\frac{q+1}{d}-1)
  =\frac{2}{d}q(q+1)(q+2)(q+1-d).  \hskip 2cm \Box
\end{array} $$
\subsection{Upper bounds of ${\hm C}$ and ${\hm D}$}
\begin{lem}\label{D}\,
${\hm C}\leq  3.85q^4$ for $d=1$; and $3.26 q^4$ for $d=3$.
\end{lem}
\demo
Recall that
${\hm C}=\sum\limits_{4, q\not\in \{j, k\},l\ne 1}n_{j,k,l}(|\{\cal E(M_\imath): M_\imath\in{\cal W}_{j,k,l}\}|-1).$
By the definition of ${\cal W}_{j,k,l}$, we get that $1\neq l\lneqq \min\{j,k\}.$
This implies that
$$\begin{array}{lll}
{\hm C}&=&3n_{3,3,2}+[2(q+1)-1](n_{3,q+1,2}+n_{q+1,3,2})+[2(q-1)-1](n_{3,q-1,2}+n_{q-1,3,2})\\
&&+[(q+1)^2-1]n_{q+1,q+1,2}+[(q-1)^2-1]n_{q-1,q-1,2}\\
&&+[(q+1)(q-1)-1](n_{q+1,q-1,2}+n_{q-1,q+1,2})+[(q-1)^2-1]n_{q+1,q+1,3}\\
&&+[(q-3)^2-1]n_{q-1,q-1,3}+[(q-1)(q-3)-1](n_{q+1,q-1,3}+n_{q-1,q+1,3})\\
&\leq&3n_{3,3,2}+(2q+1)(n_{3,q+1,2}+n_{q+1,3,2}+n_{3,q-1,2}+n_{q-1,3,2})\\
&&+(q^2+2q)(n_{q+1,q+1,2}+n_{q-1,q-1,2}+n_{q+1,q-1,2}+n_{q-1,q+1,2})\\
&&+(q^2-2q)(n_{q+1,q+1,3}+n_{q-1,q-1,3}+n_{q+1,q-1,3}+n_{q-1,q+1,3}).\\
\end{array}$$
Inserting the upper bounds for $n_{j,k,l}$  that will be established in Lemma \ref{n,jkl}, we get the upper bound for ${\hm C}$ stated in this lemma.
\qed

\begin{lem}\label{used2}
In $M\cong \SO(3,q)$, each involution is contained in at most $\frac{q+1}{2}$ subgroups $\D_4$ and $q+1$ subgroups $\D_{2(q\pm 1)}$; and each $\D_4$ is contained in at most $3$ subgroups $D_{2(q\pm1)}$.
(A routine checking)
\end{lem}

\begin{lem}\label{n,jkl}
With the notations as above, we get
\vskip -3mm $$\begin{array}{lll}
&(1)\,
n_{3,3,2}\leq \frac{q+1}{d}(q+2)(\frac{q+1}{2})^2;&\\
&(2)\, n_{3,q+1,2}+n_{3,q-1,2},n_{q+1,3,2}+n_{q-1,3,2}\leq (q+2)\frac{q+1}{2}(q+1);&\\
&(3)\, n_{q+1,q-1,3}+n_{q-1,q+1,3}+n_{q+1,q+1,3}+n_{q+1,q+1,3}\leq \frac{9(q+1)}{2};&\\
&(4)\, n_{q+1,q+1,2}+n_{q+1,q-1,2}+n_{q-1,q+1,2}+n_{q-1,q-1,2}\leq (q+2)(q+1).&
\end{array}$$
\end{lem}
\demo
Set $A:=M_\imath\cap M$ and $B:=M_\imath\cap M'$ for any $M_\imath\in{\cal W}_{j,k,l}$ with $4,q\notin\{j,k\}$ and $l\neq 1$.
\vskip 3mm
(1) Suppose  $A\cong B\cong\D_4$ and $A\cap B\cong\ZZ_2$. Then $\D_{2s}\cong \lg A,B\rg\leq M_\imath$ with $s\di (q\pm 1)$ and $s>2$. By Table \ref{number}, each subgroup $\D_s$ is contained in at most $\frac{q+1}{d}$ subgroups $\SO(3,q)$. Moreover, $M \cap M'$ contains at most $q+2$ involutions. Therefore, by Lemma \ref{used2}, we have
\vspace{-5pt}$$\begin{array}{l}
n_{3,3,2}\leq \frac{q+1}{d}|\{(t,A,B)\mid t\in A\cap B, A\cong B\cong \D_4, \D_{2s}\cong \lg A,B\rg\leq M_\imath \}|
 \leq\frac{q+1}{d}(q+2)(\frac{q+1}{2})^2.
\end{array}$$\vspace{-20pt}
\vskip 3mm
(2) For $M_\imath\in{\cal W}_{3,q\pm1,2}\cup{\cal W}_{q\pm1,3,2}$ with $A\cap B\cong \ZZ_2$, we have $\lg A,B\rg= M_\imath$. Therefore by Lemma \ref{used2}, we have
$$\begin{array}{lll}
&&n_{3,q+1,2}+n_{3,q-1,2},n_{q+1,3,2}+n_{q-1,3,2}\\
&\leq & |\{(t,A,B)\mid t\in A\cap B, A\cong \D_4, B\cong \D_{2(q\pm1)}, \lg A,B\rg= M_\imath \}|
 \leq(q+2)\frac{q+1}{2}(q+1).
\end{array}$$\vspace{-20pt}
\vskip 3mm
(3) Suppsoe that $A,B\cong \D_{2(q\pm 1)}$ and $A\cap B\cong \D_4$. Then  $\lg A,B\rg=M_\imath$. Note that  each involution in $M\cap M'$ is contained  at most $\frac{q+1}{2}$ subgroups $\D_4$. Therefore by Lemma \ref{used2},
\vspace{-5pt}$$\begin{array}{lll}
&&n_{q+1,q-1,3}+n_{q-1,q+1,3}+n_{q+1,q+1,3}+n_{q+1,q+1,3}\\
&\leq &|\{(C,A,B)\mid \D_4\cong C\leq A\cap B, A, B\cong \D_{2(q\pm1)}, \lg A,B\rg= M_\imath\}|\le \frac{q+1}{2}\times 3^2.
\end{array}$$ \vspace{-20pt}
\vskip 3mm
(4) For any given $t\in M\cap M'$ and  subgroup $A\leq M$ with $A\cong\D_{2(q\pm 1)}$, define
$${\small\begin{array}{l}
\Upsilon=\{M_\imath\mid t\in A\cap B, M\cap M_\imath=A, M'\cap M_\imath=B\cong \D_{2(q\pm1)}, M_\imath\in{\cal W}_{q+1,q\pm1,2}\cup{\cal W}_{q-1,q\pm1,2}\}.
 \end{array}}$$
\f  We shall   show   $|\Upsilon|\le 1$.   Since there are at most $q+2$ involutions in $M\cap M'$, Lemma \ref{used2} yields $n_{q+1,q+1,2}+n_{q+1,q-1,2}+n_{q-1,q+1,2}+n_{q-1,q-1,2}\leq (q+2)(q+1).$

We now prove that $|\Upsilon|\le 1$.
Let $\{\lg \a\rg , \lg \b \rg , \lg \g\rg \}$ be an fixed orthogonal frame of $\mathbf{w}$ and followed from Lemma \ref{lem-inv-bij}, we may set $\D_4\cong\lg\tau_{\a},\tau_{\b},\tau_{\g}\rg\leq A$. Without loss of generality, we assume  $t:=\tau_{(1,0,0)}\in M\cap M'$, where the coordinate triple $(1,0,0)$  corresponds to the vector $\a$ under the basis $\{\a,\b,\g\}$. Then $\mathbf{ w}'=\{\lg Yg(\a,\b,\g)^T\rg\mid Y\in\FF_q^3\}$, where
$$g={\small{{\left(
\begin{array}{ccc}
1&0&0 \\
0&1&c+di \\
0&i&c'+d'i\\
\end{array}
\right)  }}}$$
with $c'd-cd'=1$ and $c,c',d,d'\in\FF_q$.
Following a parallel argument to the one given in the proof of Lemma \ref{hlD}, we may assume that \vspace{-5pt}$$\mathbf{w}_\i=\{\lg XD(\a,\b,\g)^T\rg\mid X\in\FF_q^3\},\, {\rm with}\,
 D\in\{[1,b+i,1],[1,1,b+i],[b+i,1,1]\}.$$
Then $\mathbf{w}_\i\cap \mathbf{w}'\cap {\cal N}\neq \emptyset$ only if there exist row vectors $X,Y\in\FF_q^3$ such that
 \vspace{-5pt}$$XD(\a,\b,\g)^T=a'Yg(\a,\b,\g)^T
$$
with $a'\in\{1,a+i\}$ for some $a\in\FF_q$.
In other words, there exist two nonzero vectors $X,Y\in\FF_q^3$ and a permutation matrix $P$ satisfying the equivalent conditions:
$$XP[1,1,b+i]=a'Yg \Leftrightarrow XP=a'Yg[1,1,(b+i)^{-1}]
\Leftrightarrow Y A_{a'}=0,$$
where $A_{1}=g_AD_B+g_BD_A$ and $A_{a+i}=a(g_AD_B+g_BD_A)+(g_AD_A+g_BD_B\theta)$ for $D_A=[1,1,\frac{b}{b^2-\theta}],\,\,D_B=[0,0,\frac{-1}{b^2-\theta}],$
$g_A={\small{{\left(
\begin{array}{ccc}
1&0&0 \\
0&1&c \\
0&0&c'\\
\end{array}
\right)  }}}$
and
$g_B={\small{{\left(
\begin{array}{ccc}
0&0&0 \\
0&0&d \\
0&1&d'\\
\end{array}
\right)  }}}.$
Obviously
$$\begin{array}{ll}
&\mathbf{w}_\i\cap \mathbf{w}'=\{\lg Yg(\a,\b,\g)^T \rg \mid YA_{a'}=0, Y\in \FF_q^3\,\,\,\,\text{for some}\,\,\,\, a'\in \{1,a+i\mid a\in\FF_q\}\},\\
&\begin{array}{lll}
\mathbf{w}_{\i}\cap \mathbf{w}'\neq \emptyset&\Longleftrightarrow& |A_{a'}|=0\,\,\,\,\text{for some}\,\,\,\, a'\in\{1,a+i\mid a\in\FF_q\}\\
&\Longleftrightarrow &1\leq\rank(A_{a'})\leq 2\,\,\,\,\text{for some}\,\,\,\,a'\in\{1,a+i\mid a\in\FF_q\}.
\end{array} \end{array}$$
Now, by the same arguments as in Lemma~\ref{hlD}, we get that $M_\imath\cap M\cong\D_{2(q\pm 1)}$ only for the case $\rank(A_{a'})=1$. Then we have $|\Upsilon|\leq |\mathcal{A}|$, where
$$\begin{array}{ll}
 \mathcal{A}:=|\{b\mid \rank (A_{a'})=1,\,\,\text{for some}\,\, a'\in\{1,a+i\}\,\,\text{and some permutation matrix}\,\, P\}|.
\end{array}$$
For cases $a'=1$ and $a'=a+i$, we have
\vspace{-5pt}$$A_{1}={\small{{\left(
\begin{array}{ccc}
0&0&0 \\
0&0&\frac{bd-c}{b^2-\theta} \\
0&1&\frac{bd'-c'}{b^2-\theta}\\
\end{array}
\right)  }}} \quad {\rm and}\quad  A_{a+i}={\small {{\left(
\begin{array}{ccc}
1&0&0 \\
0&1&\frac{bc-d\theta+abd-ac}{b^2-\theta} \\
0&a&\frac{bc'-d'\theta+abd'-ac'}{b^2-\theta}\\
\end{array}
\right)  }}.}$$
\f respectively.
 Check that   $\rank(A_{1})=1$ if and only if $bd-c=0,$ which implies that either $b=\frac{c}{d}$ for $d\neq0$ or there is no solution for $b$, where $c'd-cd'=1$;
and  $\rank(A_{a+i})\neq 1$. For both cases,  $|\Upsilon|\leq |\mathcal{A}|\leq1$, as desired.
\qed

\begin{lem}\label{DD}
${\hm D}\leq 2.02 q^4$ for $d=1;$ and $2.503 q^4$ for $d=3$.
\end{lem}
\vskip -2mm \demo
By the definition, we have
\vspace{-5pt}$$\begin{array}{lll}
{\hm D}&=&\sum\limits_{\{4, q\}\cap  \{j, k\}\ne \emptyset} n_{j,k, l}(|\{\cal E(M_\imath): M_\imath\in{\cal W}_{j,k,l}\}|-1)\\
&\leq&\sum\limits_{4\in \{j, k\}} n_{j,k,l}(|\{\cal E(M_\imath): M_\imath\in{\cal W}_{j,k,l}\}|-1)+\sum\limits_{q\in \{j, k\}} n_{j,k,l}(|\{\cal E(M_\imath): M_\imath\in{\cal W}_{j,k,l}\}|-1)\\
&\leq&\sum\limits_{4\in \{j, k\}} n_{j,k,l}[3(q+2)-1]+\sum\limits_{q\in \{j, k\}} n_{j,k,l}[q(q+2)-1]\\
&\leq&2|\{M_\imath\in\O\setminus M\mid M_\imath\cap M\cong A_4\}|[3(q+2)-1]\quad\text{(as $4=j$ or $4=k$)}\\
&&+2|\{M_\imath\in\O\setminus M\mid M_\imath\cap M\cong \ZZ_p^m.\ZZ_2\}|[q(q+2)-1]\quad\text{(as $q=j$ or $q=k$)}\\
&\leq&2\times\frac{q(q^2-1)}{12}\times x_8\times[3(q+2)-1]+2\times\frac{q^2-1}{2}\times x_5\times[q(q+2)-1],\\
\end{array}$$
 where   $x_5$ and $x_8$ were given  in Table~\ref{number}, so  we get the desired upper bound of ${\hm D}$.\qed

\vskip 3mm
\begin{center}{\large\bf Acknowledgements}\end{center}
\vskip 2mm
The first author thanks the supports of the National Natural Science Foundation of China (12301446, 12571362).
The second author thanks the supports of the National Natural Science Foundation of China (12471332).
The third author thanks the supports of the National Natural Science Foundation of China (12301422).


\end{document}